\documentclass[12pt]{amsart}

\usepackage{amssymb,amsmath,amsthm,amsfonts,amscd,latexsym, color, xcolor, mathtools}
\usepackage[hidelinks]{hyperref} 
\usepackage{graphics}
\usepackage{subcaption}
\usepackage{wrapfig}
\usepackage{pst-node}
\usepackage{tikz-cd} 
\usepackage{enumitem}
\usepackage{mathrsfs}
\usepackage{svg}
\usepackage{cleveref}
\usepackage{todonotes}

\newcommand{\bbR}{\mathbb R}

\newcommand{\bbZ}{\mathbb Z}

\newcommand{\bbH}{\mathbb H}

\newcommand{\Lip}{\mathrm{Lip}}

\DeclareMathOperator{\arccosh}{acosh}
\DeclareMathOperator{\arcsinh}{asinh}
\DeclareMathOperator{\arctanh}{atanh}

\DeclareFontFamily{U}{mathx}{\hyphenchar\font45}
\DeclareFontShape{U}{mathx}{m}{n}{<-> mathx10}{}
\DeclareSymbolFont{mathx}{U}{mathx}{m}{n}
\DeclareMathAccent{\widebar}{0}{mathx}{"73}
	
\crefformat{section}{\S#2#1#3} 
\crefformat{subsection}{\S#2#1#3}
\crefformat{subsubsection}{\S#2#1#3}
	\allowdisplaybreaks

\theoremstyle{definition} 

\usepackage{thmtools, thm-restate}
\theoremstyle{plain}

\newtheorem{theorem}{Theorem}[section]

\newtheorem{proposition}[theorem]{Proposition}
\newtheorem{lemma}[theorem]{Lemma}
\newtheorem{corollary}[theorem]{Corollary}

\newtheorem{definition}[theorem]{Definition}
\newtheorem{remark}[theorem]{Remark}

\numberwithin{equation}{section}
\newtheorem*{convention}{Convention}

\newtheorem*{theorem*}{Theorem}

\newcommand{\nocontentsline}[3]{}
\let\origcontentsline\addcontentsline
\newcommand\stoptoc{\let\addcontentsline\nocontentsline}
\newcommand\resumetoc{\let\addcontentsline\origcontentsline}

\newcounter{foo}

\newtheorem{theo}[foo]{Theorem}
\newtheorem{cor}[foo]{Corollary}

\newcommand{\len}{\mathrm{len}}
	
\usepackage[margin=1in]{geometry}
\usepackage{blindtext}
\usepackage{subfiles} 

  \newcounter{Cconstant} 
  \newcommand{\newC}[1]{\refstepcounter{Cconstant}\label{#1}} 
  \newcommand{\useC}[1]{C_{\ref{#1}}}

  \newcounter{Dconstant} 
  \newcommand{\newD}[1]{\refstepcounter{Dconstant}\label{#1}} 
  \newcommand{\useD}[1]{D_{\ref{#1}}}

  \newcounter{Kconstant} 
   
  \newcounter{dconstant} 
   
  \newcounter{econstant} 
  \newcommand{\newe}[1]{\refstepcounter{econstant}\label{#1}} 

\begin{document}
\title{Flexibility and Degeneracy around a theorem of Thurston}
\author[A. Nolte]{Alexander Nolte}
\date{}
\begin{abstract}
We give two flexible and degenerate constructions related to a theorem of Thurston.
First, we produce geodesic segments for Thurston's asymmetric metric on Teichm\"uller space $\mathcal{T}(S_g)$ that remain geodesics after adding arbitrary $\varepsilon$-Lipschitz noise to all but one Fenchel-Nielsen coordinate.
Then, for all $2 < n \leq 3g-3$ we construct open sets in $\mathcal{T}(S_g)^n$ for which the limit cones of the corresponding representations in $\mathrm{PSL}_2(\bbR)^n$ are cones over explicit finite-sided polyhedra.
Each construction is as degenerate as possible and has applications to the basic structure and local non-rigidity of the involved objects.
\end{abstract}

\maketitle

\allowdisplaybreaks


\stepcounter{section}

In this paper, we identify and study degenerate phenomena in two constructions that compare lengths of geodesics in hyperbolic surfaces.
The two constructions are as follows.

First, for $X_1$ and $X_2$ marked hyperbolic metrics on a closed orientable surface $S$, consider the supremum $\mathrm{K}(X_1, X_2)$ of the ratios $\ell_{X_2}(\gamma)/\ell_{X_1}(\gamma)$ of lengths $\ell_{X_i}(\gamma)$ of geodesics in free homotopy classes of closed curves  $\gamma$ in $S$ in the two metrics.
This quantity was studied in depth by Thurston in 1986 \cite{thurston1998minimal}, who proved that $\mathrm{K}(X_1, X_2)$ is equal to the optimal Lipschitz constant of identity-isotopic maps $X_1 \to X_2$ and that $\log \mathrm{K}(X_1,X_2)$ defines an asymmetric metric $d_{\mathrm{Th}}$ on the Teichm\"uller space $\mathcal{T}(S)$ of isotopy classes of hyperbolic metrics on $S$.

Second, consider the following related construction for $n$ hyperbolic surfaces $X_i$ corresponding to pairwise non-conjugate representations $\rho_i: \pi_1 S \to \mathrm{PSL}_2(\bbR)$ ($i=1, ..., n$).
For $\gamma \in \pi_1S$, denote by $\lambda(\gamma) \in \bbR^n$ the list $(\ell_{X_1}(\gamma), ..., \ell_{X_n}(\gamma))$ of lengths of $\gamma$ in $X_i$.
Then the \textit{limit cone} $\mathcal{L}(\rho(\pi_1S))$ of $\rho = (\rho_1, ..., \rho_n)$ is the closure of the union of the rays $\bbR^+\lambda(\gamma)$ $(\gamma \in \pi_1S)$.
This is a case of a fundamental construction of Benoist \cite{benoist1997proprietes} for discrete subgroups of Lie groups (\S \ref{ss-limit-cones-results}), and \cite[\S 1.2]{benoist1997proprietes} shows $\mathcal{L}(\rho(\pi_1S))$ is convex and has nonempty interior.

There are not many convex cones in $\bbR^2$.
In particular, $\mathcal{L}(\rho(\pi_1S))$ is determined by its boundary slopes when $n =2$, which are exactly $\mathrm{K}(X_1, X_2)$ and $\mathrm{K}(X_2, X_1)^{-1}$.
When $n > 2$, there are far more cones in $\bbR^n$ and the shapes of limit cones are rich and not well-understood.

\medskip

We construct two families of oddities: geodesic segments for $d_\mathrm{Th}$ that remain geodesics after adding arbitrary $\varepsilon$-Lipschitz noise to all but one Fenchel-Nielsen coordinate (Thm.~\ref{thm:main-geos}) and open sets of $n$-tuples $(\rho_1, ..., \rho_n) \in \mathcal{T}(S)^n$ whose limit cones are cones over explicit finite-sided polyhedra (Thm. ~\ref{theo:main-cones}).
Both are related to the following peculiar theorem of Thurston:

\begin{theo}[Thurston, {\cite[Thm. 10.7]{thurston1998minimal}}]\label{thm:thurston}
On the complement of a closed set of Hausdorff codimension $1$ in $\mathcal{T}(S) \times \mathcal{T}(S)$, the supremum $\mathrm{K}(X_1, X_2)$ is realized by a simple curve $\gamma$ that remains optimal for all pairs in a neighborhood of $(X_1, X_2)$.
\end{theo} 

The generic stability of optimizers in Thm. \ref{thm:thurston} is atypical of optimization problems: usually, optimality is fragile and delicate to establish (e.g. \cite{bessonCourtoisGallot1995entropies,hutchingsMorganRitpreRos2002proof,viazovska2017sphere}).
For an adjacent example, no such stability is present when optimizing extremal length ratios on surfaces \cite[p.36]{kerckhoff1980asymptotic}.
We mention that Thm. \ref{thm:thurston} is not isolated, and now fits into a family of results in ergodic theory known as {Hunt-Ott phenomena} (e.g. \cite{bochi2018ergodic,contreras2016ground,huntOtt1996optimal}).

We prove our results by developing methods that allow one to conclude, among other things, that a given simple closed curve optimizes $\mathrm{K}(X_1, X_2)$ in classes of examples.
The stability underlying Thm.~\ref{thm:thurston} appears for us as leniency in the estimates in our proofs, and allows us to carry out quite direct arguments that would otherwise be hopeless.
Our constructions are as degenerate as possible in precise senses.
The degeneracies are rather dramatic; for instance, we prove that there are isometric embeddings of the $\mathrm{L}^\infty$ metric on $[0,1]^2$ into $(\mathcal{T}(S_g), d_{\mathrm{Th}})$ for all $g \geq 3$ (Cor. \ref{cor-l-infty}).

\medskip

One of the points of this paper is to contribute some new understanding of what is happening in the generic case of Thurston's theorem.
The idea of both our constructions is to see the stability of optimizers in Thm. \ref{thm:thurston} in the geometric isolation of simple closed geodesics in hyperbolic surfaces.
This is different from Thurston's proof of Thm. \ref{thm:thurston}, which is based on an analysis of the non-generic set.
We present our results together with the hope of bringing out their common underlying idea and their relation to Thurston's work.

Our two methods are different, and each has settings where we can use it to prove more than the other.
Broadly speaking, our construction of geodesics for $d_{\mathrm{Th}}$ is analytic and is useful in analyzing very close together hyperbolic structures.
On the other hand, our method for limit cones uses coarser estimates from Gromov hyperbolicity.
It can be used to prove stronger conclusions, but requries more separation between the involved surfaces.

In the rest of the introduction, we discuss geodesics for $d_{\mathrm{Th}}$ in \S \ref{ss-thurston-geos-results}, limit cones in \S \ref{ss-limit-cones-results}, estimates on lengths of geodesics in \S \ref{ss-non-rot-len}, further applications in \S\ref{sss-cors}, and proofs in \S \ref{ss-proof-remarks}.

\begin{remark} We prefer to work with closed surfaces, but all results in \S\ref{ss-thurston-geos-results}, \S\ref{ss-limit-cones-results}, and \S\ref{sss-cors} have straightforward analogues for hyperbolic surfaces of finite type (i.e. with finitely many punctures and geodesic boundary components).
We briefly explain this in \S\ref{ss-finite-type-thurston} and \S\ref{ss-finite-type-limit-cones}.
\end{remark}

\stoptoc

\subsection{Thurston Geodesics}\label{ss-thurston-geos-results}
Thurston showed that $d_{\mathrm{Th}}$ is geodesic in \cite{thurston1998minimal}, and geodesics for $d_{\mathrm{Th}}$ have been of interest since.
The long-term \cite{dumasLenzhenRafiTao2020coarse,lenzhenRafiTao2012bounded,lenzhenRafiTao2015shadow} and infinitesimal \cite{barnatanOhshikaPapadopoulos2025convex,huangOhshikaPapadopoulos2025infinitesimal,pan2023local} structure of general geodesics for $d_{\mathrm{Th}}$ has seen considerable study.
On short time-scales, the \textit{envelopes} of all geodesics with fixed start and end points have seen recent study \cite{dumasLenzhenRafiTao2020coarse,panWolf2024envelopes}.
On the other hand, the local structure of general geodesics for $d_{\mathrm{Th}}$ has remained unclear.

The reason for the lack of local results on general geodesics seems to be the rigidity of prior examples.
Beyond Thurston's examples and a variation in \cite{alessandriniDisarlo2022generalizing}, Pan-Wolf have a construction of distinguished geodesics for $d_{\mathrm{Th}}$ based on limits of harmonic maps \cite{panWolf2025ray} (c.f. \cite{daskalopoulosUhlenbeck2024transverse}), and Papadopoulos-Th\'eret \cite{papadopoulosTheret2012someLipschitz} and Papadopoulos-Yamada \cite{papadopoulosYamada2017deforming} have real-analytic examples based on explicit deformations of certain hexagons.
These are all rather rigid constructions, and are in particular all concatenations of real-analytic segments.

We show that general geodesics for ${d}_{\mathrm{Th}}$ are as non-rigid and irregular as possible.
The main result (Thm. \ref{thm:main-geos}) is a bit technical; here is a corollary that illustrates its flavor.

\begin{cor}[$\mathrm{L}^\infty$ regions]\label{cor-l-infty}
	Let $T > 0$ and $k \geq 1$. Then there is a real-analytic isometric embedding of $[0,T]^k$ with the $\mathrm{L}^\infty$ metric $d_{\infty}$ into $(\mathcal{T}(S_g), d_{\mathrm{Th}})$ for all $g$ so that $3g-3 \geq 2k$. 
\end{cor}

We emphasize that the embeddings in Cor. \ref{cor-l-infty} are genuine metric isometries.
So geodesics in $([0,T]^k, d_\infty)$ are mapped to geodesics in $d_{\mathrm{Th}}$, and $d_{\mathrm{Th}}$ is symmetric on the images of these embeddings.
This makes Cor. \ref{cor-l-infty} new even for $k =1$, see the discussion of Cor. \ref{thm:double-geos} below.

Cor. \ref{cor-l-infty} is an example of a theme of this paper: phenomena around Thm. \ref{thm:thurston} make statements that one expects to be approximately or asymptotically true (e.g. \cite{minsky1996extremal}) actually hold on the nose on regions inside Teichm\"uller space. See also Thm. \ref{theo:main-cones} and Cor. \ref{cor-designer-asymm}.

To contextualize the general result below, there are two fundamental constraints on Thurston geodesics.
First, from the definition of $d_{\mathrm{Th}}$ in terms of $\mathrm{K}$, all length functions are Lipschitz on geodesic segments for ${d}_{\mathrm{Th}}$.
Second, Thurston proved \cite[Cor. of Thm. 8.2]{thurston1998minimal} that along any geodesic segment $X_t : [0,T] \to \mathcal{T}(S)$ in $\mathrm{d}_{Th}$ there is a nonempty (chain-recurrent) lamination $\lambda$ that is maximally stretched in a strong sense.\footnote{For $t \in [0, T]$ let $\lambda(t)$ be the geodesic realization of $\lambda$ on $X_t$.
	Thurston proves any optimal Lipschitz map $f: X_{t} \to X_{t'}$ for $0 \leq t < t' \leq T$ sends $\lambda(t)$ to $\lambda(t')$ and expands every leaf of $\lambda(t)$ by a factor of $e^{t' - t}$.}
Thm. \ref{thm:main-geos} below shows that these are the only qualitative constraints on Thurston geodesics in general. 

To state the result, recall Fenchel-Nielsen coordinates for $\mathcal{T}(S)$.
That is, fixing a pants decompositon $\mathscr{P}= \{\gamma_1, ..., \gamma_{3g-3}\}$ for $S$, one has $\mathcal{T}(S) \cong (0,\infty)^{3g-3} \times \bbR^{3g-3}$ with the first $3g-3$ coordinates $\ell_i$ given by the lengths of the geodesic representatives of $\gamma_i$ and the final $3g-3$ coordinates given by \textit{twist parameters} $t_i$ across the curves $\gamma_i$ (see e.g. \cite{fathiLaudenbachPoenaru2012thurston}).
To get optimal constants, it is useful to replace length coordinates with their logarithms, parameterizing $\mathcal{T}(S)$ by $\bbR^{6g-6}$.
We call these \textit{log-Fenchel-Nielsen coordinates}.

\begin{theo}[Noisy Geodesics]\label{thm:main-geos}
	Let $D, T \in (0, \infty)$ be given.
	Then there is a constant $c \in \bbR$ with the following property.
	For any $X_0 = ( \ell_1, ...,  \ell_{3g-3}, t_1, ..., t_{3g-3} ) \in \mathcal{T}(S_g)$ in log-Fenchel-Nielsen coordinates with $\ell_i < c$ for all $i$, for any coordinate-wise $1$-Lipschitz $F : [0,T] \to \bbR^{3g-4}$ and any $D$-Lipschitz $G: [0,T] \to \bbR^{3g-3}$ with $F(0) = 0$ and $G(0) = 0$, the path $X_t$ $(0 \leq t \leq T)$ given in log-Fenchel-Nielsen coordinates by $$X_t = (\ell_1, \ell_2 , ..., \ell_{3g-3} , t_1 , ... t_{3g-3} ) + (t,F(t), G(t))$$ is a geodesic for $d_{\mathrm{Th}}$.
\end{theo}

The bounds of coordinate-wise $1$-Lipschitz $F$ and $D$-Lipschitz $G$ for arbitrary $D < \infty$ are optimal among Lipschitz conditions.\footnote{A stronger statement for an asymmetrization of the Lipschitz condition is given in Thm. \ref{thm:main-technical-geos}.}
One has an analogous statement in Fenchel-Nielsen coordinates, with $\varepsilon(X_0, T)$-Lipschitz $F$ in place of coordinatewise $1$-Lipschitz $F$.

Thm. \ref{thm:main-geos} has a number of corollaries due to its contrast to the rigidity of known examples of geodesics for $d_{\mathrm{Th}}$ (see \S \ref{sss-cors}).
For one such result, the rigidity of prior examples left it open whether smooth geodesics for $d_{\mathrm{Th}}$ were determined in a neighborhood of their initial point by their tangent vectors.
An immediate consequence of Thm. \ref{thm:main-geos} is a negative resolution:

\begin{cor}[Geodesic Germs Non-Rigid]\label{thm-germs-nonrigid}
	For any $k > 0$, there exist $X \in \mathcal{T}(S)$ and real-analytic Thurston geodesic segments $\gamma_1, \gamma_2$ through $X$ whose $k$-jets agree at $0$ but that agree on no neighborhood of $0$.
\end{cor}

\subsection{Limit Cones}\label{ss-limit-cones-results}

Our main results on limit cones give a class of $n$-tuples of Fuchsian representations with the surprising property that their limit cone can be computed explicitly.
These examples have the further counter-intuitive property that their limit cones are cones over {finite-sided polyhedra}, despite being defined with infinitely many rays.

\begin{theo}[Polyhedral Limit Cones]\label{theo:main-cones}
	Let $n \geq 3$ and let $n \leq k \leq 3g-3$.
	Then there is a nonempty open set $U$ of $\mathcal{T}(S_g)^n$ so that for any $(\rho_1, ..., \rho_n) \in U$ the limit cone $\mathcal{L}(\rho(\pi_1S_g))$ of the sum representation $\rho: \pi_1 S_g \to \mathrm{PSL}_2(\bbR)^{n}$ is the cone over a $k$-vertex polyhedron $P_\rho$.
	
	The polyhedra $P_\rho$ and open sets $U$ are explicitly described in Fenchel-Nielsen coordinates.
\end{theo}

Polyhedra are the most degenerate convex bodies in the sense that they are completely determined by finitely many points.
Their appearance here is a manifestation of a degeneracy phenomenon in which finitely many curves dominate all others in families of length-ratio optimization problems, as in Thm. \ref{thm:thurston}.
This finiteness is responsible for the tractability of a full description of these limit cones.

\begin{remark}
	After we proved Thm. \ref{theo:main-cones}, we learned that Danciger-Gu\'eritaud-Kassel independently and simultaneously found similar examples as part of a broader investigation into limit cones \cite{dancigerGueritaudKassel2025limit}.
	Their method is different from ours: \cite{dancigerGueritaudKassel2025limit} is based on an elegant surgery construction, and computes limit cones through proving that certain hyperplanes that support the convex hull of the projections $\lambda(\gamma)$ of {\rm{simple}} curves $\gamma$ actually support the full limit cone.
\end{remark}

Limit cones are defined in more generality\footnote{The natural setting is Zariski-dense subgroups of semisimple real linear algebraic groups \cite{benoist1997proprietes}, where the Jordan projection to a Weyl chamber plays the role of listing lengths of hyperbolic elements of $\mathrm{PSL}(2,\bbR)^n$.}
 and play a major role in the study of discrete subgroups of Lie groups (e.g. \cite{benoist2004convexesI,davaloRiestenberg2024finitesided,deyOh2025deformations,guichard2005regularity,kasselTholozan2024sharpness,sambarino2015orbital,tsouvalas2023holder,zhangZimmer2024regularity}).
Theorem \ref{theo:main-cones} and \cite{dancigerGueritaudKassel2025limit} are the first examples of finitely generated Zariski-dense subgroups of semisimple Lie groups of real rank at least $3$ with explicitly computed limit cones.
C.f.  \cite{cantrellCertReyes2024jointtranslation} for an appearance of polyhedrality in an analogous construction for finite-valence metric graphs.

We note that limit cones of general subgroups with no further hypotheses are extremely flexible.
In our setting, there are infinitely generated free subgroups $\Lambda$ of $\mathrm{SL}(2,\bbR)^n$ with limit cone any closed convex cone in $(\bbR^+)^n$ with nonempty interior \cite[\S 1.2]{benoist1997proprietes}.
However, it is far more challenging to compute the limit cone of a given subgroup $\Lambda$ of interest, and the problem becomes richer when $n \geq 3$ and there are many cones in $\mathbb{R}^n$.

\medskip

Our method always produces polyhedral cones containing $(1, ..., 1) \in \mathbb{R}^n$ and with between $n$ and $3g-3$ verticies.
We have given sharp enough proofs to produce {all} such cones.

\begin{theo}\label{theo-designer-cone}
	Let $P \subset \bbR^n$ be the cone over a convex polyhedron with $n \leq k \leq 3g-3$ vertices and so that $(1, ..., 1) \in \mathrm{Int}(P)$.
	Then there is a representation $\rho \in \mathcal{T}(S_g)^n$ with limit cone $P$.
\end{theo}

In fact, our main technical result on limit cones (Thm. \ref{thm-cones}) shows polyhedrality for a rather flexible family of representations, as follows.
Take Fenchel-Nielsen coordinates for a pants decomposition $\mathscr{P}$ satisfying a technical hypothesis (Def. \ref{def-convenient}).
Let $\rho \in \mathcal{T}(S)^n$ have all twist parameters $0$ and suppose the cone $\mathcal{C}$ over the convex hull of $\lambda(\gamma)$ for $\gamma \in \mathscr{P}$ contains $(1, ..., 1)$ in its interior.
Let $\rho_t$ $(t \in (0,1))$ be the deformations of $\rho$ given by simultaneously scaling all length coordinates by $t$. 
Then we prove that for all $t$ sufficiently small, $\mathcal{L}(\rho_t(\pi_1S)) = \mathcal{C}$.

We conclude by remarking that a corollary of our main technical Theorem \ref{thm-cones} is that the limit cone of a Zariski-dense representation $\rho \in \mathcal{T}(S)^n$ does not locally determine $\rho$.

\begin{cor}[Not Locally Rigid]\label{cor-cones-nonrigid}
	For $n \geq 3g-3$, there are submanifolds $M$ with boundary of $\mathcal{T}(S_g)^n$ of dimension $n(6g-6) - nk + 1$ so that $\mathcal{L}(\rho(\pi_1 S_g))$ is the same $k$-vertex polyhedron with non-empty interior for all $\rho \in M$.
\end{cor}

\subsection{Non-rotating length}\label{ss-non-rot-len} Our limit cones are verified through a uniform estimate on how lengths of arbitrary closed curves distort between pinched and untwisted hyperbolic surfaces (Thm. \ref{theo-distortion}).
We hope this estimate will be of broader use, and explain it here.

Roughly, we use combinatorial data to correct lengths of closed curves for contributions from twisting, with a simple scheme that is asymptotically accurate for curves with large spiraling.
We then show that the distortion of the remaining length between well-pinched surfaces is controlled by the ratios of \textit{logarithms} of lengths of the pants curves, in a uniform way with well-controlled constants.

To be precise, let $\mathscr{P}$ be a pants decomposition of $S$ and let $\mathscr{U}(\mathscr{P}) \subset \mathcal{T}(S)$ be the subset with zero twisting parameters for $\mathscr{P}$.
For $X \in \mathscr{U(P)}$, the seams of the pants decomposition glue up into a set $\mathscr{H}$ of simple closed curves.
Let $\mathscr{U}_{\varepsilon}(\mathscr{P})$ be the subset of $\mathscr{U(P)}$ in which all pants curves have length no more than $\varepsilon$.
For $\gamma_i \in \mathscr{P}$, let $\mathcal{A}(\gamma_i)$ be the annular cover corresponding to $\gamma_i$.
We call a pants decomposition $\mathscr{P}$ \textit{convenient} if it is embedded and every curve in $\mathscr{H}$ has connected intersection with each pant in $\mathscr{P}$.\footnote{Convenient pants decompositions exist for all $g \geq 2$. We expect that convenience is not necessary for our results, but prefer to work with the assumption because it simplifies combinatorial book-keeping.}

For any closed geodesic $\gamma$ in $X \in \mathscr{U(P)}$, the circularly ordered set $\mathscr{C_H}(\gamma)$ of transverse intersections of $\gamma$ with $\mathscr{H}$ is independent of the choice of $X \in \mathscr{U(P)}$.

\begin{definition}[Rotation number]\label{def-combo-rot}
	Let $\gamma$ be any closed curve and let $\gamma_i \in \mathscr{P}$.
	\begin{enumerate}
		\item \label{condition-counting} Let $n(\gamma, \gamma_i)$ be the number of pairs of consecutive $w_j, w_{j+1} \in \mathscr{C_H}(\gamma)$ so that there is a lift $\widetilde{\gamma}$ of $\gamma$ to the annular cover $\mathcal{A}(\gamma_i)$ corresponding to $\gamma_i$ so that the lifts of the curves in $\mathscr{H}$ corresponding to $w_j$ and $w_{j+1}$ intersect the core curve of $\mathcal{A}(\gamma_i)$.
		\item The {\rm{combinatorial rotation}} $r(\gamma, \gamma_i)$ of $\gamma$ about $\gamma_i$ is $\frac{1}{2} n(\gamma, \gamma_i)$.
	\end{enumerate}
\end{definition}

\begin{figure}
	\begin{center}
		\includegraphics[scale=0.43]{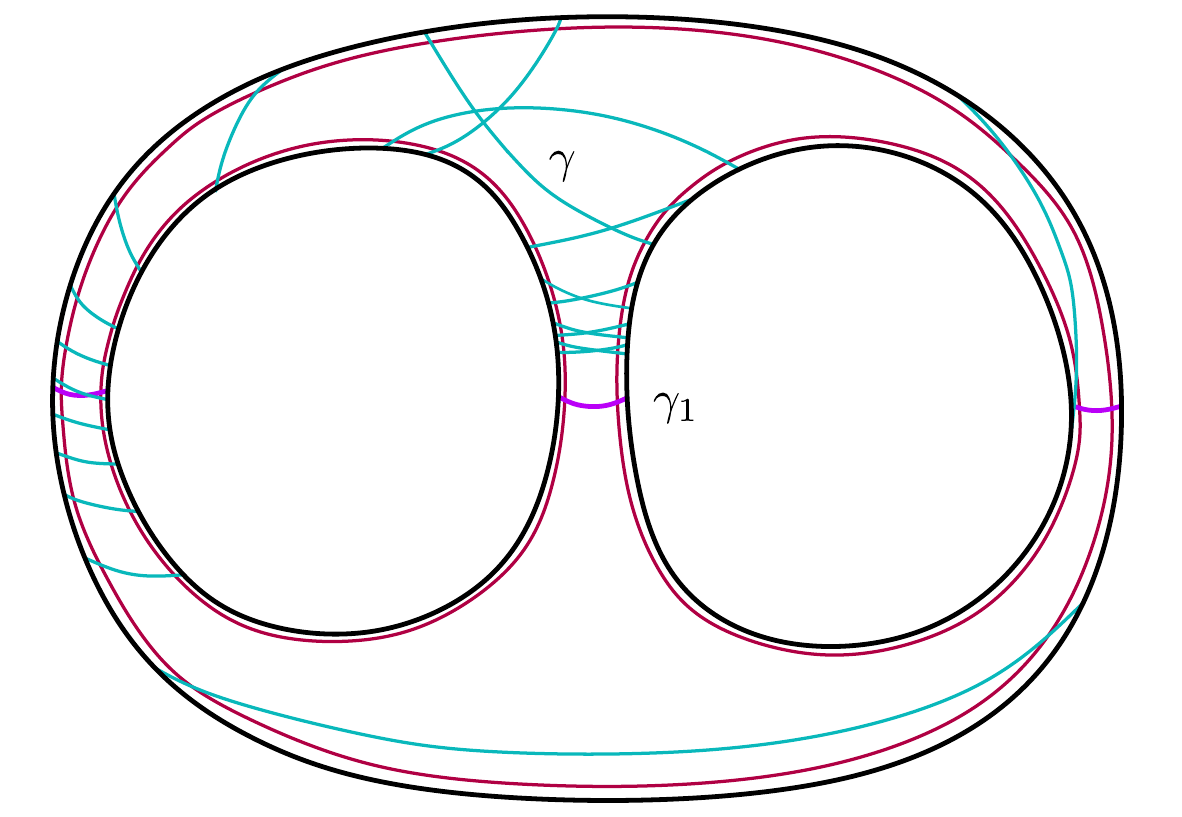} \hspace{1.5cm}
		\includegraphics[scale=0.36]{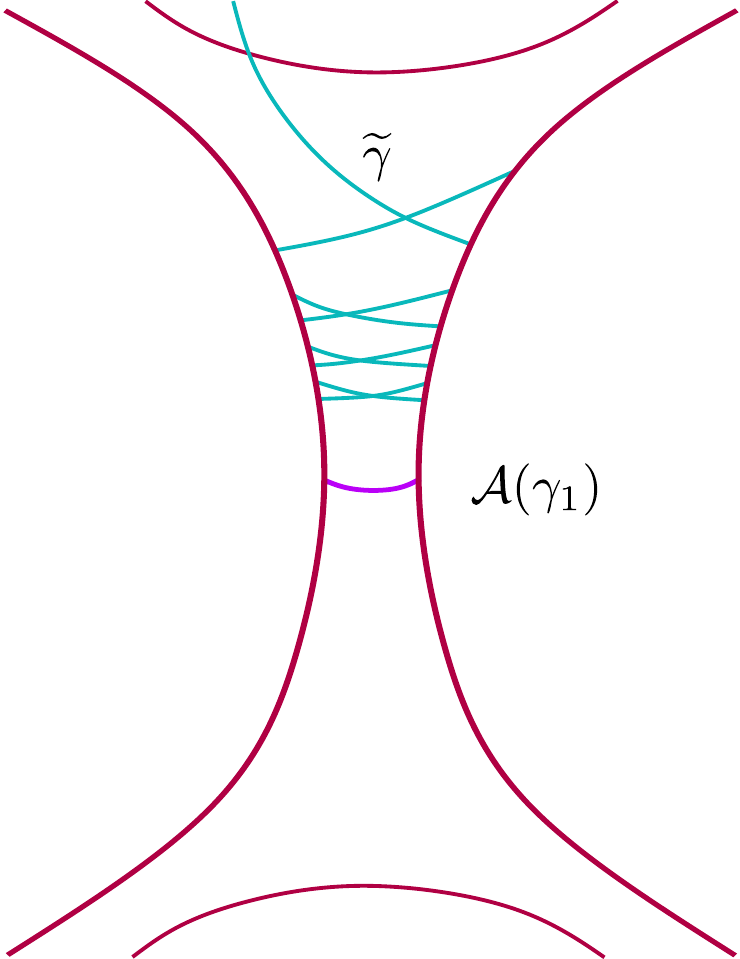}
	\end{center}
	\caption{Left: a genus two surface, with a convenient pants decomposition in purple, seams in red, and a teal closed curve $\gamma$. Right: lifts of $\gamma$ and seams curves to $\mathcal{A}(\gamma_1)$ for the middle pants curve $\gamma_1$ that demonstrates the spiraling in the middle contributes to the combinatorial rotation $r(\gamma, \gamma_1)$.}\label{fig-combo-rot}
\end{figure}

See Figure \ref{fig-combo-rot}.
Well-definition of $r(\gamma, \gamma_i)$ and independence of $X \in \mathscr{U(P)}$ is proved in \S \ref{s-len-combo}.
We correct lengths for contributions from twisting using combinatorial rotation:

\begin{definition}[Non-rotating length]\label{def-non-rotating-length}
	For $\gamma \in \pi_1(S) - \{e\}$ and $X \in \mathscr{U}(P)$, define the {\rm{rotational projection}} $\mathscr{R}(\gamma,X)$ and {\rm{non-rotating length}} $\mathscr{L}(\gamma,X)$ of $\gamma$ by \begin{align}
	\mathscr{R}(\gamma,X) = \sum_{i=1}^{3g-3} r(\gamma, \gamma_i) \ell_X(\gamma_i), \qquad \qquad
	\mathscr{L}(\gamma,X) = \ell_X(\gamma) - \mathscr{R}(\gamma,X). \label{eq-rot}
	\end{align}
\end{definition}

For small $\varepsilon > 0$ and a fixed $X \in \mathscr{U}_\varepsilon(\mathscr{P})$, the rotational projection $\mathscr{R}(\gamma,X)$ is typically small compared to $\ell_X(\gamma)$, but can be arbitrarily multiplicatively close in size to $\ell_X(\gamma)$ for curves $\gamma$ that spiral large amounts around pants curves.
We prove:

	\begin{theo}\label{theo-distortion}
		Fix a convenient pants decomposition $\mathscr{P}$.
		For every $C > 1$ there is an $\varepsilon > 0$ so that
		for all $X,Y \in \mathscr{U}_{\varepsilon}(\mathscr{P})$ and for all closed curves $\gamma$ not in $\mathscr{P}$,
			\begin{align}
				\min_{\gamma_i \in \mathscr{P}, \sigma \in \{\pm 1\}} C^{-1} \left(\frac{\log(\ell_X(\gamma_i))}{\log(\ell_Y(\gamma_i))}\right)^{\sigma}  \leq   \frac{\mathscr{L}(\gamma,X)}{\mathscr{L}(\gamma,Y)} \leq \max_{\gamma_i \in \mathscr{P}, \sigma \in \{\pm 1\}} C \left(\frac{\log(\ell_X(\gamma_i))}{\log(\ell_Y(\gamma_i))}\right)^{\sigma}. \label{eq-intro-comparison}
			\end{align}
	\end{theo}

The point of Thm. \ref{theo-distortion} is that the estimate is uniform across all (including non-simple) closed curves $\gamma$ and that $C>1$ may be chosen freely.
The logarithm terms come from ratios of the heights of collars around the pants curves.
The proof relies on a combinatorialization scheme for arbitrary closed geodesics; we discuss it in \S \ref{ss-proof-remarks}.

\subsection{Further applications on Thurston geodesics}\label{sss-cors}
We now return to Thurston geodesics to note that the flexibility and sharpness of Thm. \ref{thm:main-geos} allows some further corollaries on general geodesics and symmetrizations of $d_{\mathrm{Th}}$ to be deduced.

One direct consequence for the basic structure of geodesics is:
\begin{cor}[Irregular Geodesics]\label{thm-irregular-geodesics}
    There are smooth Thurston geodesic segments that are nowhere real-analytic, and there are nowhere smooth Thurston geodesic segments.
\end{cor}

We also remark our method constructs a family of explicit geodesic segments for ${d}_{\mathrm{Th}}$ that are also {parameterized} geodesic segments for the reversed metric $d_{\mathrm{Th}}'$ defined by $d_{\mathrm{Th}}'(X,Y)= d_{\mathrm{Th}}(Y,X)$.
These are the first such examples; note that Papadopoulos-Th\'eret \cite{papadopoulosTheret2012someLipschitz} give a family of real-analytic geodesic rays $d_{\mathrm{Th}}$ that are unparameterized geodesics for $d_{\mathrm{Th}}'$.

\begin{cor}[Noisy Symmetric Geodesics]\label{thm:double-geos}
	Let $D, T > 0$ be given.
	Then there exists a constant $c \in \bbR$ so that for any $X_0 = (\ell_1, ..., \ell_{3g-3}, t_1, ..., t_{3g-3} ) \in \mathcal{T}(S)$ in log-Fenchel-Nielsen coordinates with $\ell_i < c$ for all $i$, for any coordinatewise $1$-Lipschitz $F : [0,T] \to \bbR^{3g-5}$ and $D$-Lipschitz $G: [0,T] \to \bbR^{3g-3}$ with $G(0) = 0$ and $F(0) =0$, the path $X_t$ $(0 \leq t \leq T_0)$ given in log-Fenchel-Nielsen coordinates by $$X_t = (\ell_1, \ell_2, \ell_3,  ..., \ell_{3g-3} , t_1 , ... t_{3g-3} ) + (t,-t,F(t), G(t))$$ is a geodesic for both $d_{\mathrm{Th}}$ and the reversed metric $d'_{\mathrm{Th}}$.
\end{cor}

In particular, the subset of $\mathcal{T}(S) \times \mathcal{T}(S)$ on which $d_{\mathrm{Th}}$ is symmetric contains geodesic segments of arbitrarily long finite length.
These examples are, as a consequence, geodesics for all standard symmetrizations\footnote{E.g. $D(X,Y) = \max \{d_{\mathrm{Th}}(X,Y), d_{\mathrm{Th}}(Y,X) \} $ and $D'(X,Y) = d_{\mathrm{Th}}(X,Y) + d_{\mathrm{Th}}(Y,X)$ ($X,Y \in \mathcal{T}(S)$).} of $d_{\mathrm{Th}}$.
This pair of consequences contributes towards the open-ended Problems 6.1 and 6.2 in the problem list \cite{su2016problems}.
Note there are no full geodesic rays for both $d_{\mathrm{Th}}$ and $d'_{\mathrm{Th}}$ whose maximally strectched lamination is a multicurve (\cite[p. 2452]{theret2010divergence}), so the restriction to geodesic segments here is essential.

While $d_{\mathrm{Th}}$ is symmetric on the geodesics in Cor. \ref{thm:double-geos}, we can actually essentially prescribe the asymmetry of a geodesic for $d_{\mathrm{Th}}$ that is an unparameterized geodesic for $d'_{\mathrm{Th}}$ using our main technical result (Thm. \ref{thm:main-technical-geos}). 

\begin{cor}[Specified Asymmetry]\label{cor-designer-asymm}
	Let $T< \infty$ and $f: [0,T] \to \bbR$ be a decreasing bilipschitz function.
	Then there is a path $\gamma: [0,T] \to \mathcal{T}(S)$ with the following properties:
	\begin{enumerate}
		\item $\gamma$ is a geodesic segment for ${d}_{\mathrm{Th}}$ and is a geodesic for $d'_{\mathrm{Th}}$ after reparameterization,
		\item For all $0 \leq t \leq t' \leq T$, $d_{\mathrm{Th}}(\gamma(t'), \gamma(t)) = f(t) - f(t')$.
	\end{enumerate}
\end{cor}

\subsection{Remarks on Proofs}\label{ss-proof-remarks}
We conclude by discussing proofs of Thms. \ref{thm:main-geos} and \ref{theo:main-cones}.
Of the two, the harder proof is of Thm. \ref{theo:main-cones}.
The way that both proofs were found was through the intuition that Thurston's Theorem \ref{thm:thurston} says that when simple closed curves win optimization problems for length-ratios on hyperbolic surfaces, they typically win {with room to spare}.
A practical theme of our proofs is that this phenomenon leads to the viability of rather direct methods of verifying optimizers in related problems.
One manifestation of this is that our proofs are almost entirely based on classical hyperbolic geometry.

\subsubsection{Thurston Geodesics}
On Thurston geodesics, the conceptual point behind Thm. \ref{thm:main-geos} is an elementary construction of explicit Lipschitz-optimal maps on hyperbolic cylinders that have their Lipschitz constants dampen out as one moves away from the core.
The key contrast between these model maps and Thurston's construction \cite{thurston1998minimal} is that our model maps {only} maximally stretch the core curve, while in Thurston's examples the maximally stretched set contains no isolated curve.
We note that essentially these model maps are used near boundaries of hexagons in the geodesics of \cite{papadopoulosTheret2012someLipschitz} and \cite{papadopoulosYamada2017deforming}.

The remainder of the proof is to actually construct Lipschitz maps between the relevant hyperbolic surfaces that extend these model maps and do not increase the Lipschitz constant.
This is made realistically feasible by Lang-Schroeder's Theorem \cite{LangSchroeder1997kirszbraun} on Lipschitz extensions in negative curvature (see also \cite{gueritaudKassel2017maximally}), which reduces the construction of such maps to getting sufficient control on derivatives of geometric quantities in well-pinched hyperbolic surfaces.
Everything involved is expressed in terms of hyperbolic trigonometric functions, and the necessary estimates are neither unexpected nor pleasant.
The heart of this part is in getting estimates that work when the involved surfaces are very close; for less close-together surfaces more flexible techniques are useful, c.f. for instance \cite{calderonTao2025deflating} or Part II of this paper.

The use of Lang-Schroeder's theorem \cite{LangSchroeder1997kirszbraun} in our construction of geodesics is the main technical distinction from the methods of \cite{panWolf2025ray,papadopoulosTheret2012someLipschitz,papadopoulosYamada2017deforming,thurston1998minimal} and allows us to obtain much more flexibility in our geodesics by removing the need to be completely explicit in our optimal Lipschitz maps away from the maximally stretched lamination.

\subsubsection{Limit Cones}
The limit cones of Thm. \ref{theo:main-cones} are the convex hulls of Jordan projections of short curves in a pants decomposition, and the computation of limit cones is deduced directly from Thm. \ref{theo-distortion}, which we discuss the proof of here.

We emphasize that the main difficulty of Thm. \ref{theo-distortion} is working with all closed curves.
Indeed, the simple case of Thm. \ref{theo-distortion} is quite similar to a proposition of Mirzakhani \cite[Prop. 3.5]{mirzakhani2008growth}, which has a considerably simpler proof.

The first step in proving Thm. \ref{theo-distortion} is to get combinatorial data out of closed geodesics and bound the shape of geodesics in well-pinched and untwisted hyperbolic metrics with this data.
The basic idea is to encode the shape of geodesics in the orders of endpoints of curves they intersect in the Gromov boundary of $\pi_1S$, which is independent of the Fuchsian representation.
This basic method is known, e.g. \cite{leeZhang2017collar}.
What appears to be new here is that we consider intersections with two families of \textit{intersecting} curves, and use the structure encoded in the \textit{difference} between the orders of the corresponding sets on the top and bottom of the circle at infinity.

Our systems of curves are the pants curves in $\mathscr{P}$ and the seams curves $\mathscr{H}$; their union cuts $X \in \mathcal{T}(S)$ into hexagons.
This choice of curves has a few pleasant features.
First, the data we obtain fairly cleanly separates the combinatorics of geodesics as they cross between pants from their behavior within pairs of pants.
Next, the data within pairs of pants from seams intersections is convenient for encoding the structure of segments of closed geodesics in individual pairs of pants.
There, our basic combinatorial scheme is inspired by the set-up of the analysis of Chas-McMullen-Phillips \cite{chasMcMullenPhillips2019almost} of closed geodesics in pairs of pants.

With the structure of geodesics constrained, we split up an arbitrary closed geodesic $\gamma$ into three regimes of behavior based on a threshold $T$ for combinatorial twisting data: segments crossing between pants, segments within a single pair of pants that twist $T$ times around one pants curve, and segments in a single pair of pants without $T$-twisting segments around a pants curve.
Each type of segment can then be analyzed.
At its heart, the key point of this method is that\footnote{When $\gamma$ is not entirely contained in a well-pinched pair of pants. Otherwise, a modified proof works.} after sufficient pinching we can guarantee a lower bound on the non-rotational contributions to the length of $\gamma$ that grows linearly in the number of segments in our decomposition.
These length guarantees allow us to convert additive uncertainty in lengths of segments to multiplicative uncertainty, as appears in the statement of Thm. \ref{theo-distortion}.

\subsection{Outline}
The paper is divided into two parts: the first (\S \ref{s-models}-\ref{s-optimal-Lip-maps}) on geodesics in $d_{\mathrm{Th}}$ and the second (\S\ref{s-len-ests}-\ref{s-limit-cones}) on limit cones.
Each part begins with an outline of the flow of ideas throughout that part, and the part's sections.
Before this, we begin with a section (\S \ref{s-basics}) collecting conventions, reminders, and elementary fundamental results.
In \S \ref{s-basics}, the subsection \S\ref{ss-used-in-both} is used in both parts, \S\ref{ss-lipschitz-fundamentals} is used only in Part I, and \S\ref{ss-zariski-lemma} is used only in Part II.

The parts are almost independent, but not entirely so: our discussion of limit cones relies on the discussion of Thurston geodesics for a few technical points, most notably Lemma \ref{lemma-compact-part-estimate}.

\medskip

\par \noindent \textbf{Acknowledgements.} I want to thank Mike Wolf for suggesting I study Thurston's paper \cite{thurston1998minimal} and for introducing me to the question of infinitesimal rigidity of Thurston geodesics.
I want to thank Max Riestenberg for asking me if limit cones of Zariski dense maximal surface groups in $\mathrm{PSL}_2(\bbR)^3$ can ever be polyhedral, and particularly thank Aaron Calderon for discussions on Thurston's asymmetric metric.
I want to thank Jeff Danciger, Subhadip Dey, Fran\c cois Gu\'eritaud, Fanny Kassel, Timoth\'ee Lemistre, Yair Minsky, Hee Oh, Eduardo Reyes, \c Ca\v gr{\i} Sert, Rom\'eo Troubat, and Anna Wienhard for discussions, and Parker Evans, Fanny Kassel, and Max Riestenberg for comments.
This material is based on work supported by the National Science Foundation under DMS Grants No. 1842494, 2502952, and 2005551.

\tableofcontents
\resumetoc

\section{Conventions and Fundamentals}\label{s-basics}

\subsection{Notational Conventions}\label{ss-used-in-both}
We begin by setting notation for fundamental or repeatedly used constructions.

We write the fundamental group of a surface $S$ of genus $g\geq 2$ by $\Gamma$, denote the cardinality of a set $B$ by $\# B$, and label curves with the following convention:

\begin{convention}
	Curves are always written with Greek letters. The length of a curve is denoted by the corresponding Latin letter, or by $\ell(\cdot)$ unless otherwise specified.
	For instance, the length of a curve $\sigma_{ij}$ may be written as $s_{ij}$ or $\ell(\sigma_{ij})$.
\end{convention}	

\subsubsection{Annuli and Collars}\label{sss-collars-notation}
For any $a > 0$ there is a unique complete hyperbolic annulus $\mathcal{A}_a'$ with a core geodesic of length $a$, up to isometry.
Splitting $\mathcal{A}_a'$ along the core geodesic gives two isometric copies of a one-sided hyperbolic cylinder, which we denote by $\mathcal{C}'_a$.
The cylinder $\mathcal{C}'_a$ is foliated by curves $\eta_r$ of constant distance $r$ from the core curve, whose arc-length is given by $a \cosh(r)$.
We denote $\eta_r$ by $\eta_{a, r}$ when the core length is not clear from context.

Recall that a curve of constant geodesic curvature in $[0,1)$ in a hyperbolic surface is called a \textit{hypercycle}.
The curves $\eta_r$ are hypercycles of geodesic curvature $\tanh(r)$.

We often cut off $\mathcal{C}'_a$ along a curve $\eta_r$.
Typically, we have a fixed parameter $\delta_* > a$ and normalize so that the arc-length of the non-geodesic boundary curve $\eta_r$ has arc-length $\delta_*$. We then write the cylinder so obtained as $\mathcal{C}_{a, \delta_*}$, or suppress the $\delta_*$ and write it as $\mathcal{C}_a$.
Alternatively, at times we wish to record instead the distance $R$ from the core curve to the boundary component, and abuse notation to write the cylinder so obtained as $\mathcal{C}_{a,R}$.
When needed, we write the thickness of such a cylinder as $R_a$, suppressing the $\delta_*$-dependence.
Which notation is in use will be clear from context.

\subsubsection{Fenchel-Nielsen Coordinates}
As the phenomena we study in both parts of this paper are manifestations of the shapes of collars of geodesics, our preferred coordinate system on Teichm\"uller space is Fenchel-Nielsen coordinates.
See e.g. \cite{fathiLaudenbachPoenaru2012thurston} for their defintion and basic properties; we set notation below.

Let $\mathscr{P}$ be a pants decomposition of $S$ with pants curves $\gamma_1, ..., \gamma_{3g-3}$.
Let us write the Fenchel-Nielsen coordinates for $\mathcal{T}(S)$ as $(\ell_1, ..., \ell_{3g-3}, t_1, ..., t_{3g-3}) \in \bbR_{+}^{3g-3} \times \bbR^{3g-3}$, with $\ell_i = \ell_X(\gamma_i)$ the length of $\gamma_i$ on $X \in \mathcal{T}(S)$ and $t_i$ the twist parameter of $\gamma_i$ on $X$.
It will at times be useful to us to replace the length parameters $\ell_X(\gamma_i)$ by their logarithms, in which case the corresponding coordinate system for $\mathcal{T}(S)$ is globally $\bbR^{6g-6}$.
Let us call such coordinates \textit{log-Fenchel-Nielsen coordinates}.

When considering an individual pair of pants $P$, we write the boundary curves as $\alpha_i$ $(i = 1, ..., 3)$.
The pair of pants $P$ is split uniquely into two isometric all right hexagons $H^\pm$ by three \textit{seams curves} between $\alpha_i$ and $\alpha_j$ with $i, j \in \{1,2,3\}$ distinct.

\subsection{Lipschitz maps}\label{ss-lipschitz-fundamentals}
We now discuss the fundamental properties of Lipschitz maps and their extensions that we shall use in constructing optimal Lipschitz maps.

The main result that we shall use is a generalization of a classical result of Kirszbraun: if $C$ is any subset of $\bbR^n$ and $f: C \to \bbR^n$ is Lipschitz then $f$ admits a global extension $\overline{f}: \bbR^n \to \bbR^n$ with the same Lipschitz constant \cite{Kirszbraun1934Uber,Valentine1945lipschitz}.
This result was later extended by Lang-Schroeder \cite{LangSchroeder1997kirszbraun} to the setting of spaces of bounded curvature.
There is an essential asymmetry in Lang-Schroeder's result between Lipschitz constants of size greater than $1$ or less than $1$, depending on the signs of the curvature.
The Lipschitz constants for which the extension theorem holds are well-adapted to our uses:

\begin{theorem}[{Corollary of \cite[Thm.~A]{LangSchroeder1997kirszbraun} or \cite[Thm. 5.1]{alexanderKapovitchPetrunin2011Alexandrov}}]\label{thm:kirszbraun}
    Let $C \subset \bbH^2$ and $f: C \to \bbH^2$ be $L$-Lipschitz with $L \geq 1$.
    Then there exists an $L$-Lipschitz extension $\overline{f} : \bbH^2 \to \bbH^2$ of $f$.
\end{theorem}

This corollary is obtained by scaling metrics to arrange for $L \geq 1$ in \cite{LangSchroeder1997kirszbraun} or \cite{alexanderKapovitchPetrunin2011Alexandrov}.
This also follows immediately from \cite[Thm.~ 1.6.(1)]{gueritaudKassel2017maximally} applied with the trivial group.

\subsubsection{Gluing} We collect some elementary local-to-global and gluing phenomena for Lipschitz maps in forms that we shall use later.
The first lemma is a variation on the standard observation that a locally $L$-Lipschitz map on a geodesically convex set is $L$-Lipschitz:

\begin{lemma}[Localization]\label{lemma-elem-local-global}
    Let $X, Y$ be metric spaces and let $K \subset O \subset X$ be subsets with the property that for every $x,y \in K$ there is a path $\gamma_{xy}$ in $O$ from $x$ to $y$ with $\ell(\gamma_{xy}) = d(x,y)$.
    Let $f: O \to Y$ be locally $L$-Lipschitz.
    Then the restriction of $f$ to $K$ is $L$-Lipschitz.
    
    If $X$ and $Y$ are Riemannian manifolds, $O$ is open, $f$ is $C^1$, and $||Df(v)||/||v|| \leq L$ for all nonzero tangent vectors $v$ to $O$, then $f$ is $L$-Lipschitz on $K$. 
\end{lemma}

\begin{proof}
    For $x,y \in K$, the length of $f(\gamma_{xy})$ is at most $Ld(x,y)$, so $d(f(x),f(y))\leq L d(x,y)$.
    In the manifold case, the length of $f(\gamma_{xy})$ is similarly bounded by $L d(x,y)$.
\end{proof}

On gluing, we have:

\begin{lemma}[Gluing]\label{lemma-elem-gluing}
    Let $X$ be a geodesic metric space, let $\mathcal{K} = \{K_1, ..., K_n\}$ be a covering of $X$ by compact subsets with the property that any distance-minimizing geodesic $\gamma$ in $X$ intersects $\partial K_i$ ($i=1, ..., n$) only finitely many times, and let $f_i : K_i \to Y$ be $L$-Lipschitz maps so that $f_i = f_j$ on $K_i \cap K_j$ for all $i, j$.
    Then the induced map $f : X \to Y$ is $L$-Lipschitz.
\end{lemma}

\begin{proof}
    For a distance-minimizing geodesic $\gamma_{xy}$ from $x$ to $y$ in $X$ and $I \subset [0,t]$ with $\gamma_{xy}(I) \subset K_i$, the length of $f\circ \gamma_{xy}(I)$ is no more than $L \ell(\gamma_{xy}(I))$.
    The hypotheses ensure that $\gamma_{xy}$ may be partitioned into finitely many such segments, so that  $d(f(x),f(y)) \leq L d(x,y)$ for all $x,y \in X$.
\end{proof}

\subsection{Sums of Fuchsian Representations}\label{ss-zariski-lemma}
The following Zariski density condition is surely known to experts, and follows from a straightforward modification of Dal'bo-Kim's method in \cite{dalBoKim2000criterion}, which includes the $n = 2$ case.
We present a proof as we need this to apply Benoist's limit cone theorem in \S \ref{s-limit-cones} and could not find a source.

\begin{proposition}\label{prop-Z-closures}
	Let $\rho_1, ..., \rho_n \in \mathcal{T}(S)$ be Fuchsian and pairwise non-conjugate.
	Then the sum representation $\rho$ defined by $\rho(\gamma) = (\rho_1(\gamma), ..., \rho_n(\gamma) )$  $(\gamma \in \Gamma)$ is Zariski-dense.
\end{proposition}

\begin{proof}
	Do induction on $n$.
	The case $n = 1$ is classical and is implied by the Zariski-density of lattices \cite{borel1960density}, for instance.
	Suppose that the claim holds for $n -1$ and let $\rho$ be given by $\rho(\gamma) = (\rho_1(\gamma), ..., \rho_n(\gamma))$ with $\rho_i$ pairwise non-conjugate and Fuchsian.
	
	Consider the Zariski closure $G = \mathrm{Cl}_Z(\rho(\Gamma))$.
	Let $\pi_i : \mathrm{PSL}_2(\bbR)^n \to \mathrm{PSL}_2(\bbR)^{n-1}$ be the projection $(g_1, ..., g_n ) \to (g_1, ..., g_{i-1}, g_{i+1}, ..., g_n)$ ($i = 1, ..., n$).
	By the inductive hypothesis $\pi_i(G) = \mathrm{PSL}_2(\bbR)^{n-1}$ for all $n$.
	Let $K_i = \ker (\pi_i|_G)$.
	Note that any element $(g_1, ..., g_n)$ of $K_i$ has $g_j = e$ for all $j \neq i$ by definition of $\pi_i$, and so identify $K_i$ with a subgroup of $\mathrm{PSL}_2(\bbR)$.
	From this constraint, we see that the Zariski-dense group $\rho_i(\Gamma)$ normalizes $K_i$, and so $K_i$ is trivial or $\mathrm{PSL}_2(\bbR)$.
	Note that if $K_i = \mathrm{PSL}_2(\bbR)$ then surjectivitiy of $\pi_i$ implies $G = \mathrm{PSL}_2(\bbR)^n$.
	 
	We are left to rule out the case all $K_i$ are trivial, so suppose for contradiction that this occurs.
	Then the restrictions to $G$ of $\pi_1$ and $\pi_2$ are isomorphisms onto $\mathrm{PSL}_2(\bbR)^{n-1}$, and so the map $\Phi: \mathrm{PSL}_2(\bbR)^{n-1} \to \mathrm{PSL}_2(\bbR)^{n-1}$ given by $g \mapsto \pi_2(\pi_1^{-1}(g))$ is an isomorphism that extends $(g, g_3, ..., g_n) \mapsto (\rho_2 (\rho_1^{-1}(g)),g_3, ..., g_n)$ on the Zariski-dense image of $\pi_2$. 
	This implies that $\Phi$ may be written as $\Phi(g, g_3 ..., g_n) = (\phi(g), g_3, ..., g_n)$ with $\phi : \mathrm{PSL}_2(\bbR) \to \mathrm{PSL}_2(\bbR)$ an automorphism extending $\rho_2 \circ \rho_1^{-1}$ on $\rho_1(\Gamma)$.
	We conclude that $\rho_1$ and $\rho_2$ are conjugate, which is a contradiction.
\end{proof}

\medskip
\begin{center}
    \sc{Part I: Thurston Geodesics}
\end{center}

We first construct geodesics for $d_{\mathrm{Th}}$.
We verify the involved paths are geodesics by constructing optimal Lipschitz maps for nearby hyperbolic surfaces $X,Y$ with as many short curves as possible.
These optimal Lipschitz maps are defined peicewise, after cutting $X$ and $Y$ up into collar neighborhoods of short curves, \textit{cuff} neighborhoods of the long boundary components of these collars, and \textit{shorts} with fixed length hypercycle boundary components.

The maps on collars are explicit and the cores of collars are the only places where optimal Lipschitz constants are realized.
These are the key conceptual point.
Away from collars, our task is to extend the maps to maps from $X$ to $Y$ while not using more distortion than the Lipschitz constants of model maps on the collars.
Then maps between shorts are defined using Lang-Schroeder's theorem.
Twisting coordinates are handled on cuff regions, using a quantitative dampening out Lipschitz constants of our model maps away from cylinders' cores.
The relevant technical points are about the asymptotic behavior of shapes of hyperbolic surfaces while curves are pinched, which can be understood with hyperbolic trigonometry.

\medskip

Section \S\ref{s-models} discusses model maps on annuli: \S\ref{ss-models} defines our model maps and proves their optimality properties, and \S\ref{ss-dampening} discusses how the local Lipschitz constants of these model maps behave away from the core geodesics.
Section \S \ref{s-optimal-Lip-maps} then discusses how to extend these model maps to the entire surface.
Within \S\ref{s-optimal-Lip-maps}, in \S\ref{ss-shorts} we construct Lipschitz optimal maps between shorts, in \S \ref{ss-cuffs} we show how the dampening of Lipschitz constants from \S\ref{ss-dampening} can be used to fit twisting into cuff regions, and then in \S \ref{ss-assembly} we combine these maps to prove Thm. \ref{thm:main-geos}.
The concluding subsection \S\ref{ss-corollaries} deduces the remaining corollaries claimed on Thurston geodesics in the introduction.

\section{Model Maps on Cylinders}\label{s-models}

\subsection{Models}\label{ss-models}
    For any $a_2> a_1 > 0$ and hyperbolic cylinders $\mathcal{A}_{a_1, R_1}$ and $\mathcal{A}_{a_2, R_2}$ with $R_1 \geq R_2$ define the \textit{model map} $f_{a_1,a_2}: \mathcal{A}_{a_1, R_1} \to \mathcal{A}_{a_2, R_2}$ by sending curves of of constant distance from core curves to each other, linearly in distance from core curves and preserving crossing geodesics orthogonal to the boundaries, arranged so that $f_{a_1,a_2}$ is orientation-preserving.

Variations on this map are used throughout our construction of geodesics for Thurston's metric, both as parts of Lipschitz-optimal maps and as a technical tool.
The following collects the ways in which these model maps are Lipschitz optimal.
These maps have the curious property that both they {and their inverses} are Lipschitz optimal.

\begin{proposition}\label{prop:doubly-optimal}
For $a_2 > a_1 >0$ and $R_1 \geq R_2$, the maps $f_{a_1,a_2}: \mathcal{A}_{a_1, R_1} \to \mathcal{A}_{a_2, R_2}$ satisfy:
\begin{enumerate}
    \item $f_{a_1,a_2}$ is Lipschitz optimal: if $g : \mathcal{A}_{a_1, R_1} \to \mathcal{A}_{a_2, R_2}$ is any identity-homotopic Lipschitz map then $a_2 / a_1 = \Lip(f_{a_1,a_2}) \leq \Lip(g)$,
    \item $f^{-1}_{a_1,a_2}$ is Lipschitz optimal among maps $(\mathcal{A}_{a_2,R_2} , \partial \mathcal{A}_{a_2,R_2}) \to (\mathcal{A}_{a_1, R_1}, \partial \mathcal{A}_{a_1, R_1})$: if $g : \mathcal{A}_{a_2, R_2} \to \mathcal{A}_{a_1, R_1}$ is a map that takes boundaries to boundaries and is homotopic to the identity then $\Lip(f^{-1}_{a_1,a_2}) \leq \Lip(g)$.
    The Lipschitz constant of $f^{-1}_{a_1, a_2}$ is given by
\begin{align}
	  \max \left\{\frac{R_1}{R_2}, \, \frac{a_1\cosh(R_1)}{a_2\cosh(R_2)} \right\} \label{eq-inv-bound}.
\end{align}  
   
    \item \label{claim-inverse-bounds} Suppose additionally that $R_1/R_2 \geq (a_1 \cosh(R_1))/(a_2\cosh(R_2))$.
    Then if $h$ is any other identity-homotopic homeomorphism $\mathcal{A}_{a_1, R_1}\to \mathcal{A}_{a_2,R_2}$ so that $h$ is Lipschitz optimal and $h^{-1}$ is Lipschitz optimal among maps $(\mathcal{A}_{a_2,R_2} , \partial \mathcal{A}_{a_2,R_2}) \to (\mathcal{A}_{a_1, R_1}, \partial \mathcal{A}_{a_1, R_1})$, then $h$ differs from $f$ by a rotation.
\end{enumerate}
\end{proposition}

We remark that when the boundary components of $\mathcal{A}_{a_1, R_1}$ and $\mathcal{A}_{a_2, R_2}$ have the same length, Eq.~ (\ref{eq-inv-bound}) is given by $R_1/R_2$ and the additional hypothesis of (\ref{claim-inverse-bounds}) is satisfied.

\begin{proof}
We have $f(\eta_{a_1,r}) = \eta_{a_2, (R_2/R_1) r}$ for all $r \in [-R_1, R_1]$, with the mapping linear in arc length on these curves.
The target curve $\eta_{a_2, (R_2/R_1) r}$ has arc-length $\cosh(r R_2/R_1)a_2$.
We claim that the ratios of the lengths of any such curves is no more than at the cores.

The function of interest is then $F(t) = \cosh(Ct)/\cosh(t)$ for $C = (R_2/R_1) \leq 1$.
One may verify using calculus the desired claim that $F(t)$ is nonincreasing on $[0, \infty)$.
Indeed, let $\varphi \in C^1(\bbR^+)$ and $C \in \bbR - \{0\}$.
Then $\frac{d}{dt} (\varphi(Ct)/\varphi(t)) \leq 0$ if and only if $C(\varphi'(Ct)/\varphi(Ct)) \leq \varphi'(t)/\varphi(t)$ for all $t$.
Now, $C \cosh'(Ct)/\cosh(Ct) = C \tanh (Ct)$ and $\cosh'(t)/\cosh(t) = \tanh(t)$.
As  $0 < C < 1$, the claim follows from the monotonicity of $\tanh$.

At the endpoint $t = R_1$ we have $ F(t) = (a_2\cosh(R_2))/(a_1 \cosh(R_1))$.
As $R_2 \leq R_1$, the above shows that the derivative of the restriction of $f_{a_1,a_2}$ to any curve $\eta_{a_1, r}$ is bounded above by $a_2/a_1$ and the derivative of the restriction of $f^{-1}_{a_1,a_2}$ to any curve $\eta_{a_2,r}$ is bounded above by $(a_1\cosh(R_1))/(a_2\cosh(R_2))$.
Using this computation and the fact that transversal geodesics orthogonal to boundary components are orthogonal to curves parallel to the core, $||D f_{a_1,a_2}(v)||/||v|| \leq a_2 / a_1$ for all nonzero $v \in \mathrm{T} \mathcal{A}_{a_1,R_1}$ and \begin{align}
	\frac{||Df_{a_1,a_2}^{-1}(v)||}{||v||} \leq \max \left\{\frac{R_1}{R_2}, \, \frac{a_1\cosh(R_1)}{a_2\cosh(R_2)} \right\}
\end{align} for all nonzero $v \in \mathrm{T}\mathcal{A}_{a_1, R_1}$.
As collars are geodesically convex, this implies that $\Lip(f_{a_1,a_2}) \leq a_2/a_1$ and $\Lip(f_{a_1,a_2}^{-1})$ is bounded by Eq.~ (\ref{eq-inv-bound}).

On the other hand, if $g$ is any Lipschitz map $\mathcal{A}_{a_1, R_1} \to \mathcal{A}_{a_2, R_2}$ that is not null-homotopic then examination of the restriction of $g$ to the core curve shows $\Lip(g) \geq a_2/a_1$, and for equality to hold $g$ must map core curves to each other by arc length.
So $f_{a_1,a_2}$ is Lipschitz optimal.

Similarly, if $g$ is any Lipschitz map $(\mathcal{A}_{a_2,R_2}, \partial  \mathcal{A}_{a_2,R_2}) \to (\mathcal{A}_{a_1, R_1}, \partial \mathcal{A}_{a_1, R_1})$ then we obtain a bound $\Lip(g) \geq (a_1\cosh(R_1))/(a_2\cosh(R_2))$ by examining the outer boundary curve.
By examining the restriction of $g$ to the crossing, boundary orthogonal geodesics in $\mathcal{A}_{a_2, R_2}$ we obtain the bound $\Lip(g) \geq R_1/R_2$.
So $f_{a_1,a_2}^{-1}$ is also Lipschitz optimal, with Lipschitz constant given by Eq.~ (\ref{eq-inv-bound}).
Finally, if Eq.~ (\ref{eq-inv-bound}) is given by $R_1/R_2$, then equality implies $g$ maps boundary orthogonal geodesics to each other linearly in arc length.
Combining the equality case considerations now gives (\ref{claim-inverse-bounds}).
\end{proof}

\subsection{Away from the core}\label{ss-dampening}
A key property of these model maps is that their Lipschitz constants ``dampen'' away from the core, as follows.
Fix a constant $\delta_* > 0$, and take cylinders to have boundary hypercycle of length $\delta_*$.
For $0 < a_1, a_2 < \delta_*$ and $0 < D < \min(R_{a_1}, R_{a_2})$, let $L(a_1, a_2, D)$ denote the optimal Lipschitz constant of the restriction of the stretch map $f_{a_1,a_2}: \mathcal{C}_{a_1} \to \mathcal{C}_{a_2}$ to a $D$-neighborhood of the long boundary component of $\mathcal{C}_{a_1}$.
In this subsection we prove:

\begin{proposition}[Lipschitz Damping]\label{prop:damping}
	Fix $\delta_* > 0$ and let $D_0 > 0$, $T_0 > 0$ and $\varepsilon > 0$ be given.
	Then there is a $a_0>0$ so that for any  $0 \leq l < a_0$ and $|t| < T_0$, 
\begin{align}\label{eq:damping}
		0 \leq \log \mathrm{L}(l, e^t l, D_0) \leq \varepsilon t.
\end{align}
\end{proposition}

\subsubsection{Height Distortion}
When contracting rather than expanding the length of core curves and holding boundary lengths constant, we must control the distortion of the heights of cylinders.
To this end, for $l, t> 0$ let $R_{l}(t) = \arccosh(\delta_* l^{-1} e^{t})$ be the height of the one-sided collar with core of length $l e^{-t}$ and non-geodesic boundary curve of length $\delta_*$.
Our normalization is so that $R_{l}(t)$ is increasing in $t$.
The basic estimate here is: 

\begin{lemma}[Height Distortion]\label{lem:thickness-distortion} Let $T_0 > 0$ be fixed. 
	Then for any $\varepsilon > 0$ there is a $a_0 > 0$ so that for all $a < a_0$ and $0 \leq t \leq T_0$ we have
\begin{align}\label{eq-opp-direction-bound}
	0 \leq \log \frac{R_{a}(t)}{R_{a}(0)} \leq \varepsilon t.
\end{align}
Furthermore, we may also arrange for $0 \leq \frac{d}{dt} \frac{R_{a}(t)}{R_{a}(0)} \leq \varepsilon$ for all $a, t$ as above. 
\end{lemma}

The following observation from calculus, which follows from the bound $\log(1 + x) \leq x$ for $x \geq 0$ is helpful for proving inequalities of this form.

\begin{lemma}\label{lem:easy-est} 
	Let $f \in C^1[0,t]$ be positive and $f'(x) \in [0, C]$ for $x \in [0,t]$. Then $\log \frac{f(t)}{f(0)} \leq \frac{|t| C}{f(0)}.$
\end{lemma}

\begin{proof}[Proof of Lemma {\ref{lem:thickness-distortion}}]
	The lower bound is trivial.
	For the upper bound, note that the denominator $R_{a}(0) =  \arccosh(\delta_* a^{-1})$ grows unboundedly large as $a \to 0$ while the numerator's derivative remains uniformly bounded as $a$ grows small:
	$$ \lim_{a \to 0} \sup_{t \in [0,T_0]} \frac{d}{dt} R_{a}(t) = \lim_{a \to 0} \sup_{t \in [0,T_0]}  \frac{d}{dt} \arccosh (\delta_* a^{-1} e^{t} ) = \lim_{a \to 0} \sup_{t \in [0,T_0]}  \frac{\delta_* a^{-1} e^t}{\sqrt{(\delta_* a^{-1} e^t)^2 - 1 } } = 1.$$
	Then apply Lemma \ref{lem:easy-est}.
	The final remark is proved with the same estimates.
\end{proof}
	 
	 \subsubsection{Offset distortion}
	 
	 Let $\gamma_0$ be a geodesic in $\bbH$ and fix $\delta > 0$.
	 Take two perpendicular geodesics $\sigma_i$ and $\sigma_f$ to $\gamma_0$ at distance $x < \delta$ apart, let $\eta_{0,x}$ be the hypercycle segment parallel to $\gamma_0$ of arc-length $\delta$ bounded by $\sigma_i$ and $\sigma_f$.
	 For $t \in \bbR$ let $\eta_{t,x}$ be the hypercycle parallel to $\gamma_0$ at oriented distance $t$ from $\eta_{0,x}$ along $\sigma_i$ away from $\gamma_0$.
	 It is helpful to have an expression for the arc-length $\ell(\eta_{t,x})$ in terms of $\delta$, $t$, and $x$:
	 
	 \begin{lemma}[Offset Distortion]\label{lemma-offset}
	 	For fixed $\delta > 0$, there is a real-analytic function $g_{\delta}(x)$ defined for $x \in (-\delta, \delta)$ with $g(0) = g'(0) = 0$ so that 
	 	\begin{align}\label{eq-offset}
	 		\ell(\eta_{t,x}) &= \delta e^{t} + \sinh(t)  g_{\delta}(x) 
	 	\end{align}
	 \end{lemma}
	 
	 Note that the expression on the left of (\ref{eq-offset}) is only defined for $x > 0$ but the expression on the right is defined for $x \in (-\delta, \delta)$.
	 
	 \begin{proof}
	 	For $x > 0$, the curve $\eta_{t,x}$ has constant distance $\arccosh\left(\frac{ \delta}{x}\right) + t $ from $\gamma_0$, with length
	 	\begin{align*}
	 		x \cosh \left(\arccosh\left(\frac{ \delta}{x}\right) + t \right) &= \delta \cosh(t) + x \sinh \left( \arccosh \left(\frac{\delta}{x}\right) \right) \sinh t\\
	 		&= \delta \cosh (t) + \sinh (t) \left( \delta - x \right)^{1/2} (\delta + x)^{1/2}.
	 	\end{align*} 
	 	The first few terms of the Taylor expansions in $x$ for the square root terms are
	 	\begin{align*}
	 		\left(\delta - x  \right)^{1/2} \left(\delta + x \right)^{1/2} &= \left( \delta^{1/2} - \frac{x}{2\delta^{1/2}} - \frac{x^2}{8 \delta^{3/2}} \right) \left(\delta^{1/2} + \frac{x}{2\delta^{1/2}} - \frac{x^2}{8 \delta^{3/2}} \right) + O(x^3) \\
	 		&=  {\delta} - \frac{x^2}{2\delta} + O(x^3). 
	 	\end{align*}
	 	Note the order $1$ cancellation.
	 	The $-x^2/(2\delta) + O(x^3)$ term here is real-analytic in a neighborhood of $0$; write it as $g_\delta(x)$.
	We see that this $g_{\delta}$ satisfies the desired properties and
	 	\begin{align*}
	 		\ell(\eta_{t,x}) &= \delta \cosh (t) + \sinh (t) (\delta + g_{\delta}(x)) = \delta e^{t}  + \sinh(t) g_{\delta}(x).
	 	\end{align*}
	 \end{proof}
	 
	 \subsubsection{Lipschitz Damping}
	 We now prove the Lipschitz Damping Proposition \ref{prop:damping}.

\begin{proof}[Proof of Prop. \ref{prop:damping}]
		The estimates break into cases on whether $t \geq 0$ or $t < 0$.
		For $t < 0$, the proof of Prop.~\ref{prop:doubly-optimal} shows that $f_{a, e^t a}$ does not increase the length of any curve of constant distance to the core curve, and so Lemma \ref{lem:thickness-distortion} shows that for sufficiently small $a_0$ we have $$ 0 \leq \log \mathrm{L}(a, e^t a, D_0) \leq \log \frac{R_{a}(-t)}{R_{a}(0)} \leq \varepsilon t.$$
		
		For $t > 0$, the map $f_{a,e^ta}$ has $||Df_{a,e^ta}(v)|| \leq ||v||$ for $v$ tangent to geodesics orthogonal to the core curve.
		The proof of Prop.~\ref{prop:doubly-optimal} also shows that for tangent vectors $v$ to curves $\eta_r$ of constant distance $r$ from the core curve that $||Df_{a,e^ta}(v)||/||v||$ is monotonically decreasing in $r$.
		As any length-minimizing geodesic between two points in the $D_0$-neighborhood of the long boundary component is contained in the $(D_0 + \delta/2)$-neighborhood of the boundary component, let $v$ be a tangent vector of norm $1$ to the curve of distance $D_0 + \delta/2$ from the long boundary of $\mathcal{C}_{a}$.
		From the Offset Distortion Lemma \ref{lemma-offset}, to deduce the desired bound from Lemma \ref{lemma-elem-local-global}, the quantity to estimate to control $||Df_{a,e^ta}(v)||$ is:
		\begin{align*}
		\log \mathrm{L}(a, e^t a, D_0)\leq  \log \frac{\delta\exp(-\frac{R_a(-t)}{R_a(0)} (D_0 + \delta/2)) + \sinh( -\frac{R_{a}(-t)}{R_a(0)} (D_0 + \delta/2)) g_{\delta}(e^ta)}{\delta \exp(-(D_0 + \delta/2)) + \sinh(-(D_0 + \delta/2)) g_{\delta}(a) }.
		\end{align*}
		Lemma \ref{lem:easy-est} reduces this to noting that for $t \in [0, T_0]$, the function $$F(t, a) = \delta \exp \left( -\frac{R_a(-t)}{R_a(0)}(D_0 + \delta/2)  \right) + \sinh\left( -\frac{R_a(-t)}{R_a(0)} (D_0 + \delta/2)\right)g_{\delta}(e^ta)$$ has a uniform lower bound on $F(t,a)$ for $a$ in a sufficiently small neighborhood of $0$ and showing that $\frac{d}{dt} F(t,a)$ tends to $0$ uniformly on $[0, T_0]$ as $a \to 0$.
		This is seen directly by taking the derivative, applying Lemma \ref{lem:thickness-distortion}, and using that $g_{\delta}(0) = g'_{\delta}(0) = 0$.	
\end{proof}

\section{Maps between Polygons}\label{s-optimal-Lip-maps}

The other component of our construction is of Lipschitz extensions away from core curves.
The general method is to cut up the complements of collars into hexagons with alternating geodesic and hypercycle arc boundary segments, to map these hexagons to each other in a controlled way, then to use Lang-Schroeder's theorem to extend to interiors.

Throughout this section, fix a $\delta_* > 0$ small enough to be comfortably less than the $\delta$-hyperbolicity constant of $\bbH^2$.

\subsection{Shorts}\label{ss-shorts}
Let $P$ be a pair of pants with geodesic boundary components $\alpha_1, \alpha_2, \alpha_3$ with lengths $2a_1, 2a_2, 2a_3$.
If hypercycles of length $\delta^*$ about the boundary components bound a topological pair of pants $P_{\delta_*}$, call $P_{\delta_*}$ the \textit{$\delta_*$-shorts} corresponding to $P$.
Note that $\delta_*$-shorts are well-defined when $a_1, a_2, a_3 < \delta_*/2$.

The pair of pants $P$ can be split uniquely into two isometric all right hexagons. Denote the intersection of such a hexagon with $P_{\delta_*}$ by $H_{\delta_*}(a_1,a_2, a_3)$ and call $H_{\delta_*}(a_1, a_2, a_3)$ the \textit{$\delta_*$-short hexagons} in $P_{\delta_*}$.
For $a_i, a_i'< \delta_*/2$ ($i = 1,2,3$), denote by $\psi : \partial H_{\delta_*}(a_1, a_2, a_3) \to \partial H_{\delta_*}(a_1', a_2', a_3')$ the map that takes verticies to verticies and edges to edges linearly in arc length.
The goal of this subsection is to prove the following quantitative geometric convergence of the shorts $P_{\delta_*}(a_1,a_2,a_3)$ to $P_{\delta_*}(0,0,0)$ as $(a_1,a_2,a_3) \to (0,0,0)$.

\begin{proposition}[Shorts Maps]\label{prop:shorts}
	Let $\varepsilon > 0$ and $T_0 < \infty$. Then there is a $c > 0$ so that if $a_1,a_2,a_3 < c$, and $|t_1|, |t_2|, |t_3| \leq t  \leq T_0$, then there exists an identity-isotopic map $\Psi : P_{\delta_*}(a_1,a_2,a_3) \to P(e^{t_1} a_1, e^{t_2}a_2, e^{t_3}a_3)$ satisfying the following properties: 
	\begin{enumerate}
		\item $\Psi$ has logarithm of its optimal Lipschitz constant bounded by $\varepsilon |t| $,
		\item $\Psi$ extends $\psi$ on the boundary of each $\delta_*$-short hexagon in $P_{\delta_*}$.
	\end{enumerate}
\end{proposition}

A use of Lang-Schroder's extension Theorem will reduce the proof to controlling the optimal Lipschitz constant of $\psi$, which we denote by $\mathrm{L}(\psi)$.

\subsubsection{Preliminaries and Pentagons}
\begin{figure}
	\includegraphics[scale=0.43]{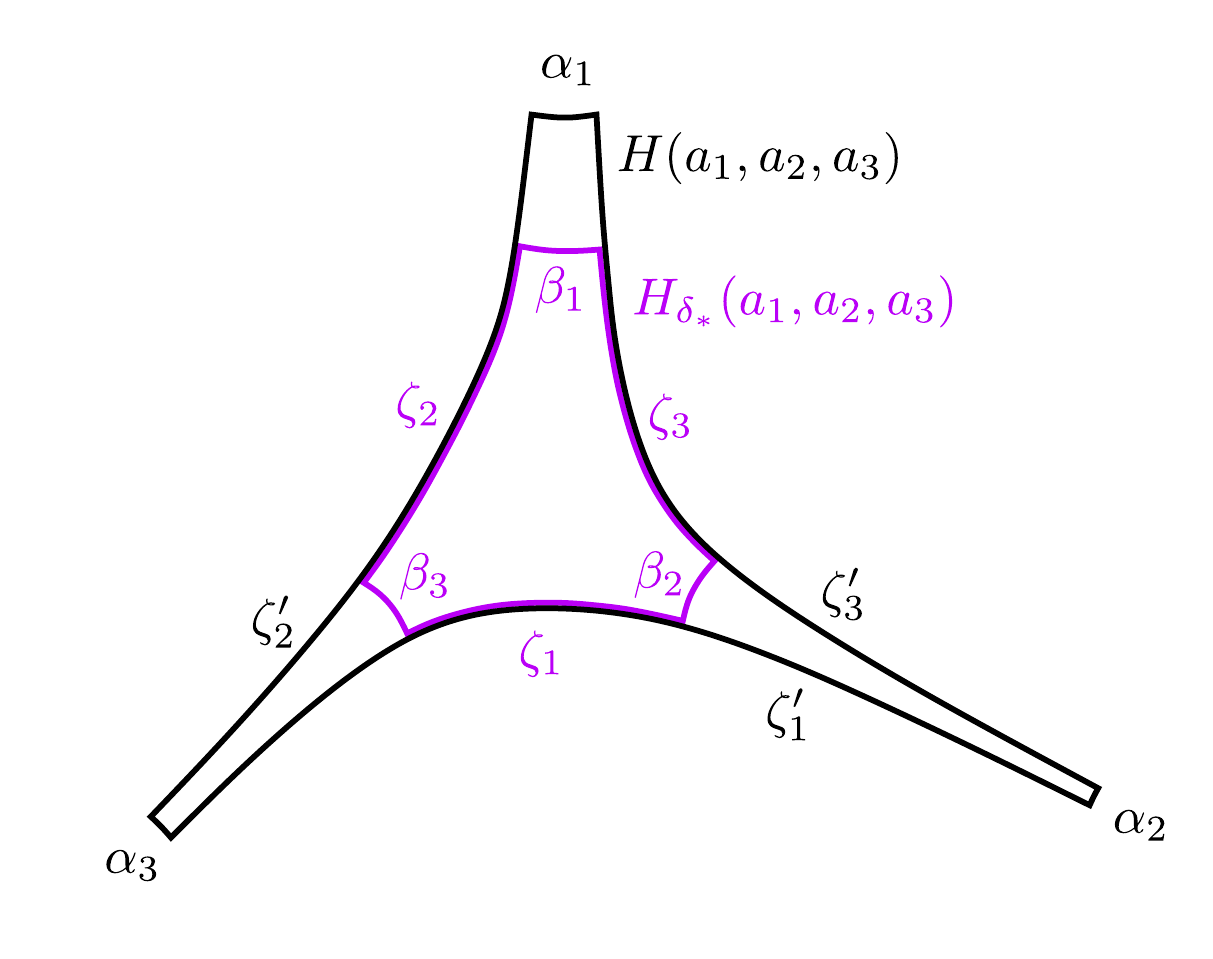} \hspace{0.7cm} \includegraphics[scale=0.43]{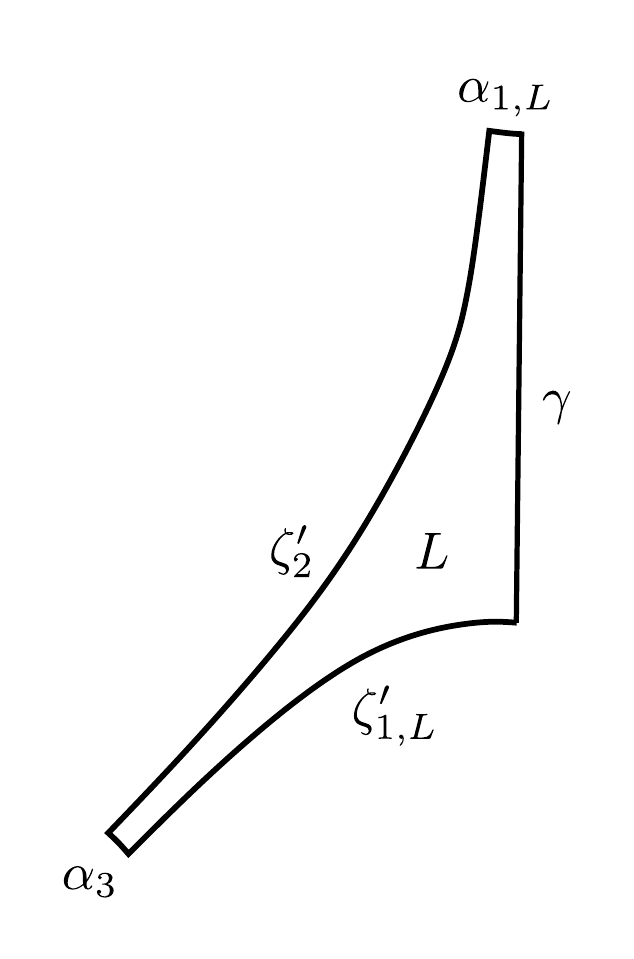}
	\caption{Labels on a hexagon and one of the pentagons it is split into.} \label{figure-pentagons}
\end{figure}

Label the sides of $H_{\delta_*}(a_1, a_2, a_3)$ by $\beta_1, \zeta_3, \beta_2, \zeta_1, \beta_3, \zeta_2$ in cyclic order so that $\zeta_i$ is a geodesic segment opposite a hypercycle arc $\beta_i$ for all $i$.
See Figure \ref{figure-pentagons}.
Write the lengths as $b_i$ for $\beta_i$ and $c_i$ for $\zeta_i$ for $i = 1,2,3$.
In particular, $\delta_*/2 = b_i$ for $i = 1, 2,3$.
Write the sides of the full hexagon that correspond to the seams of the pants $P$ as $\zeta_i'$ and write their lengths as $c_i'$ for $i=1,2,3$. 

In this paragraph we prove:
\begin{lemma}\label{prop-shorts-geos}
	As functions of $(a_1, a_2, a_3)$, the lengths $c_i$ $(i=1,2,3)$ admit real-analytic extensions to a neighborhood of $(0,0,0)$ in $\bbR^3$.
\end{lemma}

The point for our uses later is that this uniformly bounds the partial derivatives of $c_i$ on a neighborhood of $(0,0,0)$.

Lemma \ref{prop-shorts-geos} is proved with elementary hyperbolic trigonometry and calculus by finding a formula for $c_i$ that makes the claim clear.
We explain the details as there is some delicacy: the obvious formulas for $c_i$ are as sums of positive and negative terms that grow unboundedly as $(a_1, a_2, a_3)$ tends to $0$.
So showing extension across $0$ requires seeing cancellation.

We focus on $c_1$, as the roles of $c_i$ are symmetric.
The argument here is greatly simplified by splitting hexagons into right pentagons.
Let us set notation.
The hexagon $H$ may be split uniquely along a geodesic $\gamma$ perpendicular to $\alpha_1$ into two all right pentagons $K$ and $L$.

Let us take $K$ to have sides $\gamma, \alpha_{1,K}, \zeta_3',  \alpha_{2}, \zeta_{1,K}'$ and $L$ to have sides $\gamma, \zeta_{1,L}', \alpha_{3}, \zeta_2', \alpha_{1,L}$, with parallel conventions on lengths and write $\ell(\gamma) = t$.
See Figure \ref{figure-pentagons}.
In particular, $c_1' = c_{1,K}' + c_{1,L}'$ and $a_1 = a_{1,K} + a_{1,L}$.

From uniqueness of this splitting, the values of $c_{1,K}', c_{1,L}', a_{1,K},$ and $a_{1,L}$ are well-defined functions of $a_1, a_2, a_3$.
We focus on $c_{1,K}'$; the analysis of $c_{1,L}'$ is similar.
Hyperbolic trigonometry \cite[\S 7.18]{beardon1995geometry} gives the identities $a_{1,K} + a_{1,L} = a_1$, $\sinh a_{1,K} \sinh t = \cosh a_2,$ $\sinh \alpha_{1,L} \sinh t = \cosh a_3$, $\sinh c'_{1,K} = \cosh a_{1,K}/\sinh a_2$, and $\sinh c'_{1,L} = \cosh a_{1,L} / \sinh a_3$.
From these identities and hyperbolic function addition formulas, one derives the identity $F(a_{1,K}, a_1, a_2, a_3) = 0$ with
\begin{align}\label{v-identity}
	F(a_{1,K},a_1, a_2, a_3) = \sinh(a_{1,K}) \left( 1 + \frac{\cosh a_1 \cosh a_2 }{\cosh a_3} \right) - \cosh (a_{1,K}) \left( \frac{\sinh a_1 \cosh a_2}{\cosh a_3} \right),
\end{align} with $F$ well-defined on all of $\bbR^4$.
Note that the $a_{1,K}$-derivative of $F$ is non-vanishing at $0$, so that by the implicit function theorem, $a_{1,K}$ is the restriction of a real-analytic function of $(a_1, a_2, a_3)$ defined in a neighborhood of $0$.
We remark that one may explicitly compute the derivatives of $a_{1,K}$ using Eq. (\ref{v-identity}) if desired, though we shall not need these.

As $\sinh c_{1,K}'  = \cosh a_{1,K} / \sinh a_2$, we see $c_{1,K}'$ is also real-analytic in $(a_1, a_2, a_3)$ for $a_2 > 0$.

Write $c_{1,K} = c_{1,K}' - \arccosh (\delta_*/2a_2)$ and $c_{1,L} = c_{1,L}' - \arccosh (\delta_*/2a_3)$, so that $c_1 = c_{1,K} + c_{1,L}$.
Proposition \ref{prop-shorts-geos} now follows from:
\begin{lemma}\label{lem:cut-good}
	$c_{1,K}$ is the restriction of a real-analytic function of $(a_1, a_2, a_3)$ defined on a neighborhood of $0$.
\end{lemma} 

\begin{proof}
	It suffices to prove the claim for $\sinh c_{1,K}$.
	Using addition formulas for hyperbolic trigonometric functions,
	\begin{align*}
		\sinh c_{1,K} &= \sinh ( c_{1,K}' - \arccosh (\delta_* / 2a_2)) \\
				&= \frac{\delta_*\sinh(c_{1,K}')}{2a_2} - \cosh(c_{1,K}') \left(\frac{\delta_*}{2a_2}  -1\right)^{1/2} \left(\frac{\delta_*}{2a_2} + 1\right)^{1/2}.
	\end{align*}
	The first few terms of the Taylor polynomial for the square root terms are
	\begin{align*}
		\left(\frac{\delta_*}{2a_2}  -1\right)^{1/2} \left(\frac{\delta_*}{2a_2} + 1\right)^{1/2} &= \frac{1}{2a_2}  \left( \delta_*^{1/2} - \frac{a_2}{\delta_*^{1/2}} - \frac{a_2^2}{2 \delta_*^{3/2}} \right) \left(\delta_*^{1/2} + \frac{a_2}{\delta_*^{1/2}} - \frac{a_2^2}{2 \delta_*^{3/2}} \right) + O(a_2^2) \\
		&=  \frac{\delta_*}{2a_2} - \frac{a_2}{\delta_*} + O(a_2^2) . 
	\end{align*}
	The $O(a_2^2)$ term here is real-analytic in a neighborhood of $0$; write it as $g(a_2)$.
	Combining this with the identity $\sinh c_{1,K}' = \cosh a_{1,K}/\sinh a_2$, and then using the formula $\arcsinh x = \log(x + \sqrt{x^2 +1})$ we see
	\begin{align*}
		\sinh c_{1,K} &= \frac{\delta_* \sinh c_{1,K}'}{2a_2} - (\sinh c_{1,K}' + e^{-c_{1,K}'})\left( \frac{\delta_*}{2a_2} - \frac{a_2}{\delta_*} + g(a_2) \right) \\
		&= \cosh a_{1,K} \left( \frac{a_2 \delta_*^{-1} + g(a_2)}{\sinh a_2} \right) + e^{-c_{1,K}'} \left(\frac{\delta_*}{2a_2} + \frac{a_2}{\delta_*} + g(a_2)\right) \\
		&= \cosh a_{1,K} \left( \frac{a_2 \delta_*^{-1} + g(a_2)}{ \sinh a_2} \right) \\
		&\qquad + \frac{\sinh a_2}{ \cosh a_{1,K} + \sqrt{\cosh^2 a_{1,K} + \sinh^2 a_2 } }\left(\frac{\delta_*}{2a_2} + \frac{a_2}{\delta_*} + g(a_2)\right).
	\end{align*}
The claim now follows from the analyticity of $\sinh(x)/x$ and $a_{1,K}$.
\end{proof}


\subsubsection{Non-Geodesic Boundaries and Corners}
We next address the distortion of the union of two adjacent sides of $H(a_1, a_2, a_3)$.

\begin{lemma}[Corners]\label{lem:corners}
	Let $\delta_* > 0$ and $T_0 < \infty$. Then there is a $c > 0$ so that if $a_1, a_2, a_3 < c$, and $|t_1|, |t_2|, |t_3| < t  \leq T_0$, then the restriction $\psi'$ of $\psi: \partial H_{\delta_*}(a_1, a_2, a_3) \mapsto \partial H_{\delta_*}(e^{t_1}a_1, e^{t_2}a_2, e^{t_3}a_3)$ two any two consecutive sides of $\partial H_{\delta_*}$ has the logarithm $L'(\psi)$ of its optimal Lipschitz constant bounded above by $\varepsilon t$.
\end{lemma}


\begin{proof}
	Write $A = (a_1, a_2, a_3)$ and $A' = (e^{t_1}a_1, e^{t_2}a_2, e^{t_3}a_3)$.
	The domain of $\psi'$ consists of a hypercycle $\beta_i$ of length $\delta_*/2$ and a geodesic segment $\zeta_j$.
	These isometrically embed inside of the complete hyperbolic collars $\mathcal{C}'_{a_i}$ and $\mathcal{C}'_{a_i'}$.
	We write the distance from the core curves to the constant-distance curves of length $\delta_*$ in the cylinder $\mathcal{C}'_a$ as $R_a$.
	 Put
	 $$ \Omega =  \bigcup_{t \in [R_{a_i} - \delta_*/2, R_{a_i} + c_i(A)]} \eta_t \subset \mathcal{C}'_{a_i}.$$
	 
	Consider the following modification $h: \Omega \to \mathcal{C}_a'$ of the model map $f_{a_i,a_i'}$.
	On $\Omega \cap \mathcal{C}_{a_i}$ define $h$ by $f_{a_i,a_i'}$.
	On the remainder, map the curves $\eta_t$ to curves of constant distance in $\mathcal{C}_a'$, affine linearly with rate $c_i(A')/c_i(A)$.
	The Lipschitz Damping Proposition \ref{prop:damping} controls the local Lipschitz constant of $h$ on the part of $\Omega$ inside $\beta_i$.
	The verification of an appropriate bound for the local Lipschitz constant on the outside of $C_{a_i}$ is directly analogous to the proof of the Lipschitz Damping Proposition, with the modification that $c_{i}(A')/c_i(A)$ replaces $R_l(t)/R_l(0)$ and Lemma \ref{prop-shorts-geos} plays the role of Lemma \ref{lem:thickness-distortion}.
\end{proof}

\subsubsection{Embedding}\label{sss-emb}

The final basic input to Prop.~\ref{prop:shorts} is control on distances between points on non-adjacent sides of $\partial H_{\delta_*}$.

For $A = (a_1, a_2, a_3) \in \bbR_+^{3}$ Let us define \textit{normalized hexagon inclusions} $\iota_{A} : H_{\delta_*}(A) \to \bbH^2$ by requiring the corner between $\beta_{1}$ and $\zeta_3$ to go to $i$ in the upper-half-plane model of $\bbH$ with $\zeta_3$ tangent to the positive imaginary axis.
 In this paragraph, we prove:

\begin{lemma}[Embeddings Similar]\label{lem:embeddings-sim}
	Let $\varepsilon > 0$ and $T_0 < \infty$. Then there is a $c > 0$ so that if $a_1, a_2, a_3< c$, and $|t_1|, |t_2|, |t_3| < t  \leq T_0$, then with $A = (a_1, a_2, a_3)$ and $A' = (e^{t_1}a_1, e^{t_2}a_2, e^{t_3}a_3)$, for all $p \in \partial H_{\delta_*}(A)$, we have $d(\iota_{A}(p), \iota_{A'}(\psi(p))) \leq \varepsilon t$.
\end{lemma}

The remaining prelimary task is to control the embeddings of non-geodesic sides of shorts.

\begin{definition}
	For $a < \delta_*$, denote the hypercycle of length $\delta_*$ in $\mathcal{C}_a'$ by $\nu_a$.
	Take isometric identifications of the universal cover of $\mathcal{C}'_a$ and the upper-half plane model of $\bbH^2$ so that $\nu_a$ is tangent to $(i,1)$.
	Let $Q(a,t)$ be the point on $\nu_a$ of arc-length $t$ from $i$, positively oriented. 
\end{definition}

\begin{lemma}\label{lem:const-geo-curv-fine}
	For any $\delta_* < 1$, the function $Q(a,t): [0, \delta_*] \times[0, \delta_*] \to\bbR$ is the restriction of a real-analytic function defined on a neighborhood of $\{0\} \times [0,\delta_*]$ in $\bbR^2$.
\end{lemma}

The requirement $\delta_* < 1$ is an artifact of being elementary and efficient in our treatment of matrix logarithms in the proof, and makes no difference to the main results.

\begin{proof}
	Computing the holonomy $M(a)$ of the core curve of $\mathcal{C}_{a}$ in its uniformizing representation with this normalization gives 	\begin{align}
			M(a) &= \begin{bmatrix}
				\cosh(a/2) & -\exp (\arccosh(\delta_*/a)) \sinh (a/2) \\
				-\exp (-\arccosh(\delta_*/a)) \sinh(a/2) & \cosh(a/2)
				\end{bmatrix} \nonumber \\
				 &= \begin{bmatrix}
					\cosh(a/2) &   -{\sinh(a/2)} \left(\frac{\delta_*}{a} + \sqrt{ \left( \frac{\delta_*}{a} \right)^2 - 1  } \right)			\\
					-\frac{\sinh(a/2)}{\frac{\delta_*}{a} + \sqrt{ \left( \frac{\delta_*}{a} \right)^2 - 1  }} & \cosh(a/2)
				\end{bmatrix} \nonumber \\
				&= \begin{bmatrix}
					\cosh(a/2) &   -\frac{\sinh(a/2) \left({\delta_*} + \sqrt{ \delta_*^2 - a^2 } \right)		}{a}	\\
					-\frac{a \sinh(a/2)}{\delta_* + \sqrt{ \delta_*^2 - a^2  }} & \cosh(a/2)
				\end{bmatrix}. \label{eq-holonomy-good}
	\end{align}
	From the analyticity of $a^{-1}\sinh(a/2)$, the final expression in Eq. (\ref{eq-holonomy-good}) is a real-analytic function of $a$ that is defined for $|a| < \delta_*$.
	Note $M(0) = \exp( b(0))$ with $b(0) = \begin{bmatrix} 0 & \delta_* \\ 0 & 0 \end{bmatrix}$.
	As $\delta_* < 1$, a neighborhood of $M(0)$ is contained in the domain of convergence for the standard Taylor series for the matrix logarithm centered at the identity, so write $b(a) = \log M(a)$ for $a$ in a neighborhood of $0$.
	Then the desired points are $Q(a,t) = \exp [\delta_*^{-1}t  b(a)]\cdot i$, which is real-analytic and well-defined on a neighborhood of	$\{0\} \times [0,\delta_*]$ in $\bbR^2$.
\end{proof}

We now prove Lemma \ref{lem:embeddings-sim}.

\begin{proof}[Proof of Lemma \ref{lem:embeddings-sim}]
	Cyclically order the verticies and tangent data of $\partial H_{\delta_*}(A)$ as $p_i \in T^1 \bbH^2$ $(i=1, ..., 6) $ starting at $p_1 = (i, 1)$.
	Proposition \ref{prop-shorts-geos} and Lemma \ref{lem:const-geo-curv-fine} show that each $p_i$ and point $\iota_A(\psi(q))$ for $q$ in the side after $p_i$ is a real-analytic function of $p_{i-1}$ and $A$ with uniformly bounded derivative for $p_{i-1}$ in a pre-selected compact subset $K$ and $A$ in a fixed neighborhood of $0$ in $\mathbb{R}^3$ .
	For $A$ of entries sufficiently small, the map $[0,T]^3 \to \mathbb{R}^3$ given by $(t_1, t_2,  t_3) \mapsto (e^{t_i}a_1, e^{t_2}a_2, e^{t_3}a_3)$ has arbitrarily small derivatives, and the lemma follows from the chain rule.
	\end{proof}

\subsubsection{Hexagons and Extensions} We now bound the Lipschitz constants of the maps $\psi$.

\begin{proposition}[Hexagon Maps]\label{prop:hex-map}
	Let $\varepsilon > 0$ and $T_0 < \infty$. Then there is a $c > 0$ so that if $0 < a_1, a_2, a_3 < c$, and $|t_1|, |t_2|, |t_3| < t  \leq T_0$, then $\psi: \partial H_{\delta_*}(a_1, a_2, a_3) \mapsto \partial H_{\delta_*}(e^{t_1}a_1, e^{t_2}a_2, e^{t_3}a_3)$ has $\log \mathrm{L}(\psi) \leq \varepsilon t$.
\end{proposition}

\begin{proof}
	Adopt the notations in the proposition statement, and write $A = (a_1, a_2, a_3)$ and $A' = (e^{t_1}a_1, e^{t_2}a_2, e^{t_3}a_3)$.
	For $p, q$ on the same or adjacent sides of $H_{\delta_*}(A)$, the relevant estimate follows from the Corners Lemma \ref{lem:corners}.
	 For $p, q$ not on adjacent sides of $H_{\delta_*}(A)$ begin by noting that there is a uniform $\delta > 0$ so that $d(p, q) \geq \delta$ for all $A$ in a neighborhood of $0$, as $\delta_*$ is fixed.
	 
	For points on non-adjacent sides, we view $H_{\delta_*}(A)$ and $H_{\delta_*}(A')$ as (isometrically) embedded in $\bbH^2$ by $\iota_A$ and $\iota_{A'}$.
	In Lemma \ref{lem:embeddings-sim}, take $T_0$ so that in the conclusion we obtain a bound of $t\varepsilon\delta /2$.
	Using the triangle inequality and Lemma \ref{lem:embeddings-sim}, we get
		\begin{align*}
		\log \frac{d(\psi(p), \psi(q))}{d(p,q)} &\leq \log \left( 1 +  \frac{ d(\psi(p), p) + d(q, \psi(q))}{d(p,q)} 
		\right) \leq \frac{d(p, \psi(p)) + d(q, \psi(q))}{\delta} \leq t \varepsilon.
	\end{align*}
\end{proof}

We now have everything that we need to prove Prop. \ref{prop:shorts}.

\begin{proof}[Proof of Prop.~\ref{prop:shorts}]
	Denote the two hexagons in $P_{\delta_*}(A)$ and $P_{\delta_*}(A')$ as $H_{\delta_*}^\pm(A)$ and $H_{\delta_*}^{\pm}(A')$, respectively.
	Isometrically embed $H(A)$ and $H(A')$ inside $\bbH^2$.
	After these embeddings, $\psi : \partial H_{\delta_*}^\pm(A) \to \partial H^\pm_{\delta_*}(A')$ is $\exp(|t|\varepsilon)$-Lipschitz by Prop.~\ref{prop:hex-map}.
	Thm.~\ref{thm:kirszbraun} now produces a $\exp(|t|\varepsilon)$-Lipschitz extension $\Psi^{\pm}_0 : H^\pm_{\delta_*}(A) \to \bbH^2$, and we define $\Psi^{\pm}$ as the post-composition of $\Psi^\pm_0$ with closest-point projection to $H(A')$.
	As $H(A')$ is convex and $\bbH^2$ is nonpositively curved, closest point projection to $H(A')$ is $1$-Lipschitz, and so $\Psi^{\pm}$ is $\exp(|t|\varepsilon)$-Lipschitz.
	Define $\Psi$ by gluing together $\Psi^+$ and $\Psi^-$ along $\partial H_{\delta_*}^\pm(A)$, noting that the maps agree on $\partial H_{\delta_*}^\pm(A)$ and are equal to $\psi$ there.
	Then $\Psi$ satisfies all desired properties, and is isotopic to the identity from the contractibility of hexagons.
\end{proof}

\subsubsection{Uniform Estimates on Shorts}
The following slight modification of Prop.~\ref{prop:shorts} will be used in the second part of the paper, but is not immediately needed.
We record it here because it is conveniently accessible with what has been done so far.

\begin{lemma}[Uniform Shorts Maps]\label{lemma-compact-part-estimate}
	Let $\delta > 0$ and $C > 1$ be given.
	Then there is an $\varepsilon \geq 0$ so that for any $P, P' \in \mathcal{P}_\varepsilon$ there is a $C$-Lipschitz identity-isotopic map $f : P_{\delta} \to P'$ extending $\psi$ on each $\delta$-short hexagon in $P$.
\end{lemma}

\begin{proof}
	The same proof as Prop.~\ref{prop:hex-map} can be used to show that there is an $\varepsilon > 0$ sufficiently small so that for all $a_1, a_2, a_3 < \varepsilon$, there are $C^{1/2}$-Lipschitz maps $\partial H_{\delta}(a_1, a_2, a_3) \to \partial H_{\delta}(0,0,0)$ and $\partial H_{\delta}(0,0,0) \to \partial H_{\delta}(a_1, a_2, a_3)$.
	Compose these maps together and apply Thm.~\ref{thm:kirszbraun}, as in the proof of Prop.~\ref{prop:shorts}.
\end{proof}

\subsection{Cuffs}\label{ss-cuffs}
We next turn to maps on cuffs.
Let us define a \textit{cuff} to be a region in a hyperbolic cylinder bounded by two hypercycles and not containing a closed geodesic.

Cuffs $O(a, \delta_1, \delta_2)$ are determined by three data points: the length $a$ of the core curve of the hyperbolic cylinder containing $O(a,\delta_1, \delta_2)$ and the lengths $\delta_1 < \delta_2$ of the boundary components.
We also write $O'(a,L,\delta_2)$ to be the cuff in a cylinder of core length $a$, long boundary component of length $\delta_2$ and height $L$.
Let $F(a,\delta_1, \delta_2)$ be the distance between the top and bottom of $O(a,\delta_1, \delta_2)$, so that $F(a, \delta_1, \delta_2) = \arccosh(\delta_2/a) - \arccosh(\delta_1/a )$.

\begin{lemma}\label{lemma-twisting-aux-control}
	For any $\delta_2 > \delta_1 > 0$, the function $F$ admits a strictly positive real-analytic extension to a neighborhood of $(0,\delta_1, \delta_2)$.
	Furthermore, this extension satisfies $\frac{\partial}{\partial a} F(0, x, y) = 0$ for all $x < y$. 
\end{lemma}

\begin{proof}
	By now, the proof is routine.
	We prove the claim for $\cosh(F)$.
	Expanding Taylor series in $a$, using hyperbolic trigonometry, and noting the same order-$1$ cancellations as in Lemmas~\ref{lemma-offset} and \ref{lem:cut-good}:
		\begin{align*}
		\cosh\left( \arccosh\left( y/a \right) - \arccosh\left( x/a \right) \right) &= \frac{xy}{a^2} - \sinh (\arccosh (y/a)) \sinh (\arccosh (x/a)) \\
		&=  \frac{xy}{a^2} - \frac{1}{a^2} \left( (x - a)^{1/2}(x +a)^{1/2} (y -a)^{1/2}(y + a )^{1/2} \right) \\
		&= \frac{x}{2 y} + \frac{y}{2 x} + G(a,x, y), 
	\end{align*}where $G$ is a real-analytic function of $a$ defined in a neighborhood of $(0, \delta_1, \delta_2)$ and so that $G(0, \delta_1,\delta_2) = 0$.	
	The claim on the $a$-derivative vanishing is because $G$ is even in $a$, from the intermediary expression in terms of square roots obtained above.
\end{proof}

We will find it technically useful to uniformize cuffs $O(a, \delta_1, \delta_2)$ to flat cylinders.
Recall that every annulus $A$ is conformally equivalent to a unique flat cylinder $C_M = S^1 \times [0,M]$, and that the \textit{modulus} $\mathrm{Mod}(A)$ of $A$ is defined as $M$.

From Lemma \ref{lemma-twisting-aux-control} and Lemma \ref{lemma-offset} we see that the model maps between annuli $O'(a, L, \delta_2)$ for fixed $L$ and $\delta_2$, for $a$ in a neighborhood of $0$ are uniformly $K$-quasiconformal.
So from basic properties of the modulus and that $O'(0, L, \delta_2)$ grows unboundedly as $L \to \infty$:
\begin{corollary}[Modulus Bound]\label{cor:mod-bd}
	For any $M < \infty$ and $\delta_2 >0$, there is an $L < \infty$ and $a_0 > 0$ so that $\mathrm{Mod}(O'(a,L,\delta_2)) \geq M$ for all $a < a_0$.  
\end{corollary}

\subsubsection{Twisting}
For a cuff $O'$ of modulus $M$, let $C_M$ denote the flat cylinder $S^1 \times [0,M]$, let $\varphi: O' \to C_M$ be a conformal equivalence, and let $\lambda(p)$ be the conformal factor at $p$.
For $s > 0$, let $r_s : C_M \to C_M$ be the standard rotation of slope $s$ given by $(\theta,r) \mapsto (\theta + rs,r)$, and let $\mathrm{tw}_s : O' \to O'$ be $\varphi^{-1} \circ r_s \circ \varphi$.

Take oriented orthonormal bases for tangent spaces $(b_o,b_l)_p, (c_o,c_l)_q$ of points $p \in O'$ and $q \in C_M$ so that $b_o$ and $c_o$ are tangent to the orbits of the isometric rotations of $O'$ and $C_M$.
Next, note that $\varphi$ conjugates conformal equivalences of $O'$ to conformal equivalences of $C_M$ and so takes curves of constant distance from the core of the cylinder $O'$ to curves of constant height in $C_M$.
As $r_s$ preserves curves of constant height, it follows that the conformal factors $\lambda$ are rotation-invariant and that $D\varphi_{p}$ is diagonal with respect to the distinguished bases.
With respect to these bases, as $D\varphi_{p}$ is then diagonal and conformal, hence $\lambda(p) \mathrm{Id}$.
From this discussion, we also see that the total twisting of $\mathrm{tw}_s$ is $Ms$.

Note that the differential of $r_s$ with respect to an orthonormal basis adapted to the product structure of the flat metric on $C_M$ is a shear by $s$.
Computing with the chain rule for $p \in O'$, with respect to the basis $(b_o,b_l)_p$ the differential $(D \mathrm{tw}_s)_p$ is also a shear by $s$.
So, the top singular value of $(D \mathrm{tw}_s)_p$ may be computed to be $\sigma(s) = (s^2/2 + s(s^2  + 4)^{1/2}/2 + 1)^{1/2}$, which is in particular a smooth function in a neighborhood of $0$ with $\sigma(0) = 1$.
So $\mathrm{tw}_s$ extends to a $\sigma(s)$-Lipschitz homeomorphism of the half-cylinder containing $O'$ with total twisting $Ms$.

Twisting is now handled by the following: 
\begin{proposition}[Room to Twist]\label{prop:room-to-twist}
	Let $C_0, T_0< \infty$.
	There exists an $a_0 > 0$ so that if $a < a_0$, $t < T_0$, and $|s| < C_0 t$ then for all $k \in [e^{-C_0t}, e^t]$, there is an $e^t$-Lipschitz homeomorphism $f_s : \mathcal{C}_a \to \mathcal{C}_{ k a}$  of twisting $s$ that takes boundary components to boundary components by arc length.
\end{proposition}

\begin{proof}
	Adopt the strategy of pre-composing $f_{\nu, k\nu}$ with the map that is $\mathrm{tw}_s$ on $L$-collars of the boundaries of $C_\nu$ and the identity otherwise.
	Prop.~\ref{prop:doubly-optimal} and Lemma \ref{lem:thickness-distortion} show that for sufficiently small $a$, maps $f'_s : \mathcal{C}_a \to \mathcal{C}_{ka}$ exist with the desired Lipschitz constants and no twisting.
	Prop.~\ref{prop:damping} ensures that on cuffs these maps $f'_s$ have a linear lower bound of log-Lipschitz constants left to spend on twisting while remaining $e^t$-Lipschitz. Smoothness of $\sigma(s)$ together with the uniform lower bound on the module of $L$-collars of boundary components (Cor.~\ref{cor:mod-bd}) ensures that a linear lower bound of twisting can be done within this freedom in log-Lipschitz constant.
	Taking $a_0$ to be small enough to gaurantee large enough moduli of cuffs through Cor.~\ref{cor:mod-bd} for all $t < T_0$, the claim follows.
	\end{proof}

\subsection{Assembly}\label{ss-assembly}
We now have everything used to certify Lipschitz optimality for enough maps to prove Thm.~\ref{thm:main-geos}. 
To give an optimal statement, we make the following definition.
\begin{definition}
	A function $f: [0, T] \to \bbR$ is {\rm{$(F,B)$-Lipschitz}} if $|f(x) - f(y)| \leq F|x -y|$ for all $x \leq y$ and $|f(x) - f(y)| \leq B|x-y|$ for all $x \geq y$.  
\end{definition}
Then we have:

\begin{theorem}\label{thm:main-technical-geos}
	Let $D, T \in (0, \infty)$ be given.
	Then there exists a constant $c \in \bbR$ with the following property.
	For any $X_0 = ( \ell_1, ..., \ell_{3g-3}, t_1, ..., t_{3g-3} ) \in \mathcal{T}(S_g)$ with $ \ell_i < c$ for all $i$, for any coordinatewise $(1, D)$-Lipschitz $F : [0,T] \to \bbR^{3g-4}$ and any $D$-Lipschitz $G: [0,T] \to \bbR^{3g-3}$ with $F(0) = G(0) = 0$, the path $X_t$ $(0 \leq t \leq T)$ given in log-Fenchel-Nielsen coordinates by $$X_t = (\ell_1,  \ell_2 , ...,  \ell_{3g-3} , t_1 , ... t_{3g-3} ) + (t,F(t), G(t))$$ is a geodesic for $d_{\mathrm{Th}}$.
\end{theorem}

\begin{proof}
	Let us initially take $c$ sufficiently negative to be able to decompose all $X_t$ into long collars around core curves of length $e^{\ell_i}$ ($i = 1, ..., 3g-3$) and $\delta_*$-shorts for small fixed $\delta_*$.
	Write $G = (g_1, ..., g_{3g-3})$.
	We will take $c$ sufficiently small to apply Props.~\ref{prop:shorts} and \ref{prop:room-to-twist} as needed below.
		We must prove that for all $t_1 < t_2 \in [0,T]$ that there is a $e^{t_2 - t_1}$-Lipschitz map $\Psi_{t_1, t_2}$ in the appropriate isotopy class of maps $X_{t_1} \to X_{t_2}$.
	
	We define $\Psi_{t_1, t_2}$ peicewise on the $\delta_*$-shorts and collars.
	On $\delta_*$-shorts, we take $\Psi_{t_1,t_2}$ to be the map $\Psi$ of Prop.~\ref{prop:shorts}; it has $\mathrm{L}(\Psi) \leq \varepsilon' |t_1 - t_2| $ where we may arrange for any $\varepsilon'$ of our choosing with appropriate choice of $c$, though all we shall need is that $\varepsilon' < 1$.
	On collars, we define $\Psi_{t_1,t_2}$ to be the stretch-and-twist map of Prop.~\ref{prop:room-to-twist} to have twisting $g_k(t_2) - g_k(t_1)$ on the $k$-th cylinder segment for all $k$.
	By taking $c$ sufficiently small at the start, Prop.~\ref{prop:room-to-twist} shows $\Psi_{t_1,t_2}$ has its optimal Lipschitz constant bounded by $e^{t_1 - t_2}$ for all $0 \leq t_1 <  t_2 \leq T$.
	The asymmetry on the requirements of $F$ is used here, through the asymmetry of the bounds in the conclusion of Prop.~\ref{prop:room-to-twist}.
	The maps on collars and cylinders glue together to be maps $X_{t_1} \to X_{t_2}$ in the correct isotopy class by the compatibility of the two constructions on the boundaries of shorts and cylinders, and the correct twisting being arranged to glue with changes in twist parameters given by $G$. 
	
	Furthermore, the decomposition of $X_{t_i}$ into shorts and cuffs is geometrically acceptable to glue together Lipschitz constants with Lemma \ref{lemma-elem-gluing} and obtain that $L(\Psi_{t_1, t_2}) \leq t_2 - t_1$.
	Because the core curve $\gamma_1$ has expanded by a factor of exactly $e^{t_2- t_1}$, we have that $L(\Psi_{t_1, t_2}) = t_2 - t_1$ and $d_{\mathrm{Th}}(X_{t_1}, X_{t_2}) = t_2 - t_1$, as desired.
\end{proof} 

\subsection{Corollaries}\label{ss-corollaries}
Let us explain how to deduce the remaining points on Thurston geodesics from the introduction.

Cor.~\ref{thm:double-geos} follows from Thm.~\ref{thm:main-technical-geos} because $f(t) = -t$ is $1$-Lipschitz and Cors.~\ref{thm-germs-nonrigid} and \ref{thm-irregular-geodesics} follow from basic facts on Lipschitz functions.
Let us prove Cor.~\ref{cor-l-infty}.

\begin{proof}[Proof of Cor. {~\ref{cor-l-infty}}]
	Let $T < \infty$ be given and $2k \leq 3g-3$.
	For $X_0 = (\ell_1, ..., \ell_{3g-3}, t_1, ..., t_{3g-3})$ in log-Fenchel-Nielsen coordinates, let $\Psi_{X_0} : [0,T]^k \to \mathcal{T}(S)$ be given for $(x_1, ..., x_k) \in [0,T]^k$ by
	$$\Psi_{X_0}(x_1, ..., x_k) = X_0 + (x_1, -x_1, x_2, -x_2, ..., x_k, -x_k, 0, ..., 0).$$  
	Then Thm.~\ref{thm:main-technical-geos} shows that for $X_0$ with sufficiently negative $\ell_1, ..., \ell_{3g-3}$ that $$d_{\mathrm{Th}}(\Psi_{X_0}(x_1, ..., x_k), \Psi_{X_0}(y_1, ..., y_k)) = \max_{i=1, ..., k}  |x_i - y_i| \qquad ((x_1, ..., x_k), (y_1, ..., y_k) \in [0,T]^k).$$
	This proves Cor.~\ref{cor-l-infty}.
\end{proof}

The only remaining point on $d_{\mathrm{Th}}$ claimed in the introduction is now Cor.~\ref{cor-designer-asymm}.

\begin{proof}[Proof of Cor. {~\ref{cor-designer-asymm}}]
Let $f(x) : [0,T] \to \bbR$ be decreasing, $B$-bilipschitz, and normalize $f$ so $f(0) = 0$.
For the initial value $X_0 : (\ell_1, ..., \ell_{3g-3}, t_1, ..., t_{3g-3})$, take $\gamma_{X_0}(t) : [0,T] \to \mathcal{T}(S)$ to be given in log-Fenchel-Nielsen coordinates by $X_0 + (t, f(t), 0, ..., 0)$, and let $\eta_{X_0} : [0,-f(T)] \to \mathcal{T}(S)$ be the reparameterization $s \mapsto X_0 + (f^{-1}(-s), -s, 0, ..., 0)$.
Then using that $f$ is decreasing, Thm.~\ref{thm:main-technical-geos} shows that for $X_0$ with sufficiently negative log-length coordinates, $\gamma_{X_0}$ is a geodesic for $d_{\mathrm{Th}}$, that $\eta_{X_0}$ is a geodesic for $d'_{\mathrm{Th}}$ with the same image, and $d'(\gamma_{X_0}(t), \gamma_{X_0}(t')) = f(t) - f(t')$ for all $0 \leq t \leq t' \leq T$.	
\end{proof}

\subsection{Finite type surfaces}\label{ss-finite-type-thurston}
We briefly explain how to deduce analogous results in the general finite-type case, i.e. in which $S$ may have punctures and geodesic boundary components.
The statements are modified by replacing $3g-3$ with the number $\kappa(S)$ of curves in a pants decomposition of $S$, and excluding $0$-dimensional Teichm\"uller spaces from Cor.~\ref{thm-germs-nonrigid} and \ref{cor-designer-asymm}.

The case of surfaces $S'$ with geodesic boundary components and no punctures is obtained by doubling $S'$ to a closed surface $S$ and then restricting optimal Lipschitz maps to $S'$.
This works because the optimal Lipschitz maps of Thm.~\ref{thm:main-technical-geos} map each pair of pants in the decomposition induced by $\mathscr{P}$ into itself, and so define maps of $S'$.

To also handle punctures, one uses that cusp neighborhoods with horocycle boundary are determined up to isometry by the length of their boundary.
So modify the maps of Thm.~\ref{thm:main-technical-geos} to be by these isometries on cusp neighborhoods.
The only point of note is that extension across $0$ in Lemmas \ref{prop-shorts-geos} and \ref{lem:const-geo-curv-fine} allows for the same proof as Prop. \ref{prop:shorts} to go through.

	\begin{center}
	{\textsc{Part II: Limit Cones}}
\end{center}

Almost all of this part is spent proving Thm.~\ref{theo-distortion}, and the final relevant statement is Prop.~\ref{prop-combo}. 
Then polyhedrality of the limit cones we consider is deduced quickly.

We begin with a brief outline.
The basic idea is to break up a closed geodesic into its intersections with pairs of pants, then analyze contributions to its length from a few different factors.
It is useful to group these contributions to length into three classes, by the order of how they distort between distinct pinched, barely twisted, marked hyperbolic surfaces $X,Y$:
\begin{enumerate}
	\item\label{contrib-type-bad} The spiraling of geodesics around core curves $\gamma_i$ of long collars.
	The ratios of length contributions from segments of this type behave on the order of $[\ell_X(\gamma_i)/\ell_Y(\gamma_i)]^{\pm1}$.
	\item\label{contrib-excursion} Portions of geodesics that correspond to excursions deep into collars or crossings of collars, corrected to remove spiraling contributions.
	The ratios of length contributions from this type of behavior are no greater in order than $[\log(\ell_X(\gamma_i))/\log(\ell_Y(\gamma_i))]^{\pm1}$.
	\item\label{contrib-thick} Portions of geodesics that stay in thick parts of pants.
	The ratios of length contributions from this type of behavior are bounded by uniform constants.
\end{enumerate}

We mention that properly identifying and partitioning into these cases and setting thresholds for consideration are major points.
A crucial feature of this decomposition is that we can arrange that a curve picks up a definite amount of length whenever it is considered for contributions of the form (\ref{contrib-type-bad}).
This allows us to convert \textit{additive} errors in estimates coming from $\delta$-hyperbolicity into uniformly bounded \textit{multiplicative} errors in the contributions of type (\ref{contrib-excursion}) and (\ref{contrib-thick}).

In order to exhibit the desired phenomena, we need tight enough estimates so that this behavior is separated out.
The estimation procedure is also complicated by the need for our estimates to be uniform across \textit{all} isotopy classes of closed (not necessarily simple) curves and for the estimates to hold across a fairly large family of hyperbolic structures.
Note that though working only with simple closed curves is possible in the case of $n=2$ due to work of Thurston \cite[\S 3]{thurston1998minimal}, Thurston's proof of this breaks when $n > 2$.

\medskip

Section \S\ref{s-len-ests} gives the basic estimates on hyperbolic lengths and geodesics in annuli used in our combinatorialization.
Then \S\ref{s-len-combo} combinatorializes the lengths of geodesics on pinched surfaces: this is by far the longest section in this part.
In \S\ref{ss-combo-basics} we set basic notation and recall fundamentals.
In \S\ref{ss-combo-combo} we collect invariant data and use this data to constrain the shape of geodesics in pinched surfaces: \S\ref{sss-crossing-data} describes the data and its basic properties, \S\ref{sss-twisting-data} and \S \ref{sss-rot-track} use this data to encode twisting, and \S\ref{sss-init-decomposition} uses this data to give well-behaved partitions of geodesics.
Then \S\ref{ss-combo-length-ests} finishes the combinatorialization scheme and deduces the relevant length estimates.
Combinatorialization of lengths of crossing segments is done in \S\ref{sss-crossing-combo}, the decomposition of segments of our partition that stay in a pair of pants and its basic properties is established in \S\ref{sss-internal-curve-control}, and the main estimate is finally deduced in \S\ref{sss-general-combinatorialization}.
In \S\ref{s-limit-cones} we use our combinatorialization to construct sums of Fuchsian representations with polyhedral limit cones.

\section{Basic Length Estimates}\label{s-len-ests}

We begin by giving definitions and basic estimates from hyperbolic geometry of annuli that will be used throughout this part.

\subsection{Twisting in Annuli}
Let us begin with a notion of twisting.
For $a > 0$, let $\mathcal{A}_a'$ be as in \S \ref{sss-collars-notation} above, and write the core geodesic of $\mathcal{A}_a'$ as $\beta_a$.
There is a canonical closest-point projection to $\beta_a$ from $\mathcal{A}'_a$.
For any curve $\gamma$ in $\mathcal{A}'_a$, define the \textit{basic rotation} $\bar r(\gamma)$ as the net (real, positive) number of rotations around $\beta_a$ that the projection of $\gamma$ to $\beta_a$ makes.

A notion that is less well-defined but has properties that are useful for invariance arguments is as follows.
Let $\sigma, \sigma'$ be orthogonal geodesics to $\beta_a$ in $\mathcal{A}_{a}'$ spaced exactly half-way apart and let $\gamma$ be a curve in $\mathcal{A}_a'$ that is monotone after projection to $\beta_a$, such as a geodesic in $\mathcal{A}_a'$.
We define the \textit{synthetic rotation} $|r|_{\sigma}(\gamma) \in \frac{1}{2}\mathbb{N}$ of $\gamma$ with respect to $\sigma$ to be half the number of intersections of $\gamma$ and $\sigma, \sigma'$.
Synthetic rotation is well-defined to within $1/2$, regardless of the choice of $\sigma$.

A basic bound on rotation from hyperbolic geometry that we shall use is:
\begin{lemma}[Thick Twisting Bound]\label{lemma-twist-bd}
	Fix $\delta_* >0 $. Then for all $\varepsilon < \delta_*$ sufficiently small, there exists a $C = C(\varepsilon, \delta_*) < \infty$ with the following property.
	
	Let $\alpha$ be any geodesic segment in $\mathcal{C}_{a}'$ with $0 \leq a \leq \varepsilon$ entirely contained in the complement of $\mathcal{C}_{a}$, i.e. $\alpha$ is contained in the part of $\mathcal{C}_{a}'$ foliated by hypercycles of length at least $\delta_*$.
	Then the basic rotation of $\alpha$ is no more than $C$.
\end{lemma}

We remark that $C(\varepsilon, \delta_*)$ is not large in practice for large $\delta_*$.
For instance, in the hyperbolic cusp $\mathcal{C}_0'$ with $\delta_* = 1$, the sharp constant is $C=2$.

\begin{proof}
	It is helpful to lift this setup to $\bbH^2$.
	One sees that the basic rotation of $\alpha$ is always finite because the geodesic containing $\alpha$ does not share an endpoint with the lift of the core geodesic $\beta_a$ (or the cusp endpoint if $a = 0$).
	Let $\alpha(0)$ be a lift of the point in $\alpha$ that makes the deepest excursion into $\mathcal{C}_a$, and let $\sigma$ be the geodesic in $\bbH^2$ orthogonal to $\beta_\alpha$ and containing $\alpha(0)$.
	By cutting $\alpha$ at $\alpha(0)$, it suffices to show the claim for geodesic rays that start at their deepest excursion into $\mathcal{C}_a'$.
	
	For these rays, working in $\bbH^2$ also shows that the total rotation is monotone increasing when $\alpha$ is pushed deeper into $\mathcal{C}_a$ by a loxodromic transformation preserving $\sigma$.
	Since $\alpha(0)$ is the point of deepest excursion of $\alpha$, the segments $\alpha$ and $\sigma$ have an oriented angle of no more than $\pi/2$ pointing in the direction away from $\beta_a$ on $\sigma$.
	Basic hyperbolic geometry also shows that the total rotation of a geodesic ray through $\alpha(0)$ and with such an angle in $[0, \pi/2]$ is monotone increasing under increasing the angle up to $\pi/2$ with an elliptic transformation that holds $\alpha(0)$ fixed.
	
	So it suffices to prove the claim for geodesic rays $\alpha$ tangent to the hypercycle $\nu_{\delta_*}$ of length $\delta_*$, and to prove this it suffices to consider complete geodesics $\alpha$ tangent to $\nu_{\delta_*}$.
	The synthetic rotation of $\alpha$ is then given by the number of orbits of $0 \in \partial \bbH^2$ under powers of $M(a)$ that are contained in $[-1,1]$, where $M(a)$ is the holonomy of the core curve of $\mathcal{C}_a$, normalized so that the tangent vector to $\nu_{\delta_*}$ at $\alpha(0)$ is $(i,1)$ in the upper half-plane model of $\bbH^2$.
	This is in fact the very same M\"obius transformation that appeared in Lemma \ref{lem:const-geo-curv-fine} of Part I and is written in Eq.~\ref{eq-holonomy-good}).
	The claim follows from the expression of Eq.~(\ref{eq-holonomy-good}) and the unifom comparability of basic rotation and synthetic rotation.
\end{proof}

\subsection{Spiraling estimates}\label{ss-combo-spiraling}
We now turn to estimating lengths of geodesics that spiral around the core of a hyperbolic annulus.

\subsubsection{Non-crossing estimates}
We begin with geodesic segments in cut-off cylinders.
We fix a $\delta_*  \in (0,1]$ and take cut-off cylinders and annuli to have hypercycle boundaries of length $\delta_*$. 

The object of interest to us here is the quantity, for a geodesic $\alpha_{a,t}$ with boundary on the long boundary component of $\mathcal{C}_a$ and with total basic rotation $t$, $$\mathcal{L}_a(t) = \ell(\alpha_{a,t}) - ta.$$
The goal of this paragraph is to get a useful bound on $\mathcal{L}_a(t)$ (Prop.~\ref{prop-residue-bound}).

There is a maximal $s(a,t) > 0$ so that $\alpha_{a,t}$ intersects the hypercycle $\eta_{s(a,t)}$.
Define the \textit{excursion depth} $\mathcal{E}(a,t)$ of $\alpha_{a,t}$ to be the distance between $\eta_{s(a,t)}$ and the long boundary component of $\mathcal{C}_{a}$.
Recall that $R_a$ is the distance between the boundary components of $\mathcal{C}_a$ and $R_a = \arccosh(\delta_*/a)$.

\begin{figure}
	\begin{center}
		\includegraphics[scale=0.43]{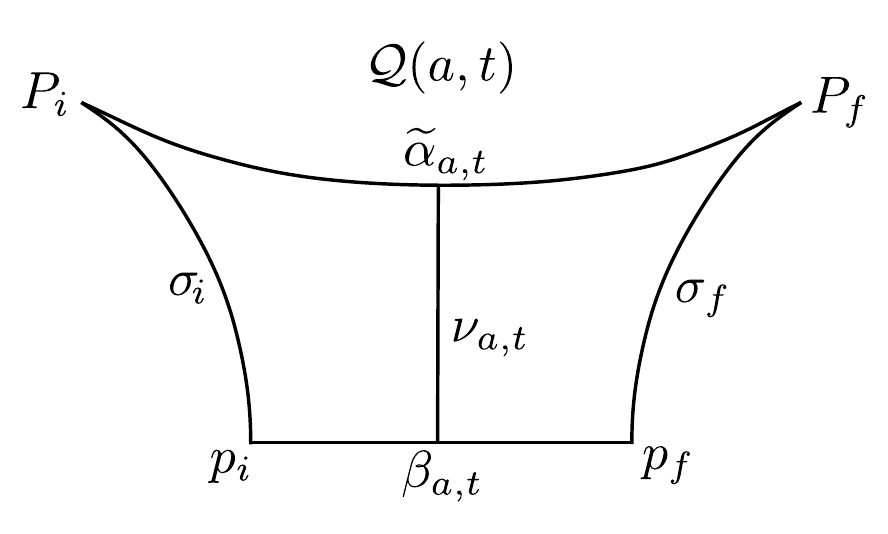}
	\end{center}
	\caption{The basic set-up after lifting to the universal cover.}\label{fig-lift-setup}
\end{figure}

We estimate $\mathcal{E}(a,t)$ through the following construction.
Lift to the universal covering of $\mathcal{C}_{a}$, then place two geodesics $\sigma_i$ and $\sigma_f$ orthogonal to the lift $\widetilde{\beta}_t$ of the core geodesic spaced $ta$ apart.
Call $\beta_{a,t}$ the cut out segment of $\widetilde{\beta}_{t}$.
See Figure \ref{fig-lift-setup}.
The resulting shape has four corners, cyclically ordered as $p_i, p_f, P_f, P_i$ with $p_i$ and $p_f$ on the lift of the core geodesic.
Then a lift $\widetilde{\alpha}_{a,t}$ of $\alpha_{a,t}$ is the geodesic between $P_i$ and $P_f$.
Denote by $\mathcal{Q}(a,t)$ the hyperbolic quadrilateral with sides $\beta_{a,t}, \widetilde{\alpha}_{a,t}, \sigma_i,$ and $\sigma_f$.
Denote by $\mathcal{P}(a,t)$ the quadrilateral obtained by splitting $\mathcal{Q}(a,t)$ in half by a geodesic segment $\nu_{a,t}$ orthogonal to the midpoint $b$ of $\beta_{t}$.
Then $\mathcal{E}(a,t) = R_a - \ell(\nu_{a,t})$.

It seems worthwhile to note that this construction gives an analytic expression for $\mathcal{E}(a,t)$, as $\mathcal{P}(a,t)$ is a hyperbolic quadrilateral with three right angles.
Such quadrilaterals are known as \textit{Lambert quadrilaterals}, and have well-understood analytic expresions for side lengths, e.g.~\cite{beardon1995geometry}.
From these, one sees:
\begin{align}\label{explicit-expression}
	\mathcal{E}(a,t) = R_a - \arctanh \left( \frac{\tanh(R_a)}{\cosh(ta)} \right).
\end{align}
In particular, $\mathcal{E}(a,t)$ is monotone increasing in $t$ and real-analytic.

We have found coarse estimates from $\delta$-hyperbolicity to generally be more useful than Eq. (\ref{explicit-expression}) in proving the estimates we shall need below.
Let $\delta_\bbH$ be the optimal Gromov hyperbolicity constant of $\bbH^2$.
The basic estimates here are: 

\newC{good-len-approx}
\newC{C-excursion-error}

\begin{lemma}[Excursion Estimates]\label{lemma-excursion-ests}
	Fix $\delta_* >0$ and let $\mathcal{O} > 0$ be given.
	Then there are finite positive constants $\useC{good-len-approx}$ and $\useC{C-excursion-error}$ with only $ \useC{C-excursion-error}$ depending on $\mathcal{O}$ so that for any $0 < a < \delta_*$:
	\begin{enumerate} 
		\item\label{additive-coarse-excursion-bound} $|2 \mathcal{E}(a,t)  - \mathcal{L}_a(t)| \leq \useC{good-len-approx}$ for all $t \geq 0$,
		\item\label{claim-E-approx} $|\mathcal{E}(a,t) - \log t| \leq \useC{C-excursion-error}$ for $1 \leq t \leq 4\delta_\bbH/a + \mathcal{O}$,
		\item\label{claim-deep-E} $R_a - 2\delta_\bbH  \leq \mathcal{E}(a,t) \leq R_a$ for $t \geq 4\delta_\bbH/a$.
	\end{enumerate}
\end{lemma}

\newC{C-triangle}
\begin{proof}
	We retain the above notation.
	For Claim (\ref{additive-coarse-excursion-bound}), consider the two right triangles formed in $\mathcal{P}(a,t)$ after adding the geodesic segment $\zeta_{a,t}$ from $P_i$ to the midpoint $b$ of $\beta_{a,t}$.
	The two triangles meet along the long side $\zeta_{a,t}$.
	A basic fact from hyperbolic trigonometry (which follows, for instance, from the law of cosines) is that the length of the long side of a right hyperbolic triangle is within a uniform constant $\useC{C-triangle}$ of the sum of the shorter side lengths.
	Because $\ell(\nu_{a,t}) = R_a - \mathcal{E}(a,t)$, we then have
	\begin{align*}
		|\mathcal{E}(a,t) - \mathcal{L}_a(t)/2| &= | R_a - \ell(\nu_{a,t}) - (\ell(\widetilde{\alpha}_{a,t}) - ta)/2 | \\ &\leq | R_a + ta/2 - \ell(\zeta_{a,t})| + | \ell(\zeta_{a,t}) - (\ell(\widetilde{\alpha}_{a,t})/2 + \ell(\nu_{a,t}))| \\
		&\leq 2\useC{C-triangle}.
	\end{align*}
	This proves (\ref{additive-coarse-excursion-bound}).

	Claim (\ref{claim-E-approx}) is the most involved.
	We prove it by showing that the strategy of going down to a hypercycle of fixed length, crossing between $\sigma_i$ and $\sigma_f$ on that hypercycle, and going back up to the top is not very inefficient in comparison to $\widetilde{\alpha}_{a,t}$ in this time range.
	The main point is getting a uniform upper bound on $\ell(\eta_{s(a,t)})$ in this time range, which then determines a good length of hypercycle to go down to.
	The argument is a bit complicated by this strategy not being well-defined for times that are too large or too small relative to $a$, which makes the edges of the time range here need to be handled with separate considerations.
		\newC{C-hypercycle-bound}
		
	Let us begin.
	We claim that there is a $\useC{C-hypercycle-bound}$ so that $\ell(\eta_{s(a,t)}) \leq \useC{C-hypercycle-bound}$ for all $a < \delta_*$ and $1 \leq t \leq 4\delta_\bbH/a + \mathcal{O}$.
	We do case analysis on whether or not $\mathcal{E}(a,t) \leq R_a - 2\delta_\bbH$.
	
	Note first that if $\mathcal{E}(a,t) \geq R_a - 2\delta_{\bbH}$ and $ta \leq 4\delta_\bbH + \mathcal{O}a$ then $\ell(\eta_{s(a,t)}) \leq (4\delta_\bbH + \mathcal{O}\delta_*) \cosh(2\delta_\bbH)$.
	This handles the case $\mathcal{E}(a,t) \leq R_a - 2\delta_\bbH$.

	So assume that $\mathcal{E}(a,t) < R_a - 2\delta_\bbH$.
	Let $p_m$ be the midpoint of $\widetilde{\alpha}_{a,t}$ and let $r_m$ be the intersection point of $\sigma_i$ and $\eta_{s(a,t)}$.
	Let $\mathcal{I}(a,t)$ be the all-right quadrilaterial with three geodesic sides formed by $\nu_{a,t}$, arcs of $\beta_{a,t}$ and $\sigma_i$, and half of $\eta_{s(a,t)}$.
	
	We prove the desired bound by applying Gau\ss-Bonnet to $\mathcal{I}(a,t)$.
		As $\mathcal{I}(a,t) \subset \mathcal{P}(a,t)$ and $\mathcal{P}(a,t)$ is a geodesic quadrilaterial with three right angles, the area $A(a,t)$ of $\mathcal{I}(a,t)$ is bounded above by $\pi/2$.
	As $\mathcal{E}(a,t) < R_a - 2\delta_\bbH$, the geodesic curvature $\kappa(a,t)$ of $\eta_{s(a,t)}$ satisfies $\kappa(a,t) \geq \tanh(2\delta_\bbH)$.
	Now because $\mathcal{I}(a,t)$ has four right angles and three geodesic sides, Gau\ss-Bonnet gives
	$$ -A(a,t) = - \frac{1}{2}\int_{\eta_{s(a,t)}} \kappa(a,t). $$
	The bounds on these quantities then show that $\pi/\tanh(2\delta_\bbH)$ is an upper bound for $\ell(\eta_{s(a,t)})$ under these hypotheses.
	So with $\useC{C-hypercycle-bound} = \max\{\pi/\tanh(2\delta_\bbH), (4\delta_\bbH + \mathcal{O}\delta_*) \cosh(2\delta_\bbH)\}$, we have $\ell(\eta_{s(a,t)}) \leq \useC{C-hypercycle-bound}$ for all $t \in [1, t\delta_\bbH/a + \mathcal{O}]$.

	We next observe that as $ta \leq 4\delta_\bbH + \mathcal{O}\delta_*$, Claim (\ref{additive-coarse-excursion-bound}) shows that $\ell(\widetilde{\alpha}_{a,t})/2$ and $\mathcal{E}(a,t)$ differ by at most an $(a,t)$-independent additive constant, and so it suffices to prove the claim for $\ell(\widetilde{\alpha}_{a,t})/2$.
	Finally, we note that using the very coarse estimate $0 \leq \ell(\widetilde{\alpha}_{a,t}) \leq t\delta_*$, we can adjust $\useC{C-excursion-error}$ to make the desired bound trivial for $t < \useC{C-hypercycle-bound}/\delta_*$, and henceforth assume $t \geq \useC{C-hypercycle-bound}/\delta_*$.
	
	Write the depth of the hypercycle with boundary length $\useC{C-hypercycle-bound}$ as $D(a,t)$.
	Using the curve obtained by going down to this hypercycle, going over, then going back up as a competitor for $\widetilde{\alpha}_{a,t}$ gives the bound $\ell(\widetilde{\alpha}_{a,t})/2 \leq D(a,t) + \useC{C-hypercycle-bound}/2$. 
	On the other hand, as $\useC{C-hypercycle-bound}$ bounds $\ell(\eta_{s(a,t)})$ above and the length of hypercycles is monotone in distance from the core, we have $\mathcal{E}(a,t) \geq D(a,t)$.
	Comparing the estimates this gives on $d(P_i, \beta_{a,t})$ and $d(p_m, \beta_{a,t})$, we then get $\ell(\widetilde{\alpha}_{a,t})/2 \geq D(a,t)$.
	So $|D(a,t) - \ell(\widetilde{\alpha}_{a,t})/2| \leq \useC{C-hypercycle-bound}/2$.
	
	We now estimate $D(a,t)$.
	It is a technically useful remark that  $\useC{C-hypercycle-bound} \geq 4\delta_\bbH + \mathcal{O}\delta_*$, though it would be trivial to arrange this if it were not already the case.
	By definition, $\useC{C-hypercycle-bound} = \cosh ( R_a - D(a,t))ta$, so that $D(a,t) = R_a - \arccosh(\useC{C-hypercycle-bound} /ta)$.
	This expression is valid for all $t$ under consideration because $\useC{C-hypercycle-bound} \geq 4\delta_\bbH + \mathcal{O}\delta_*$.
	Now, using the definition of $R_a$ and the bound $|\arccosh(x) - \log(x)| \in [0,\log2)$ for $x \in [1,\infty)$,
	\begin{align*}
		|D(a,t) - \log t| &= |\arccosh (\delta_*/a) - \arccosh(\useC{C-hypercycle-bound}/ta) - \log t|  \\
		&\leq |\log(\delta_*/a) - \log(\useC{C-hypercycle-bound}/ta) - \log t| + \log 4 \\
		&\leq | \log (\delta_*/\useC{C-hypercycle-bound})| + \log 4,
	\end{align*}
	and Claim (\ref{claim-E-approx}) follows.
	
	Finally, for Claim (\ref{claim-deep-E}), observe that in this case $\ell(\beta_{a,t}) \geq 4\delta_\bbH$.
	As $\mathcal{Q}(a,t)$ is a quadrilateral, every point $p \in \partial \mathcal{Q}(a,t)$ is within the $2\delta_\bbH$-neighborhood of the sides of $\mathcal{Q}(a,t)$ not containing $p$.
	So $p_m$ is within a $2\delta_\bbH$-neighborhood of $\beta_{a,t}$, and the claim follows.
\end{proof}

The main proposition of this subsection is now:

\begin{proposition}\label{prop-residue-bound}
	Let $\delta_* > 0$ be fixed, and let $\tau > 0$ and $C > 1$ be given.
	Then there are positive constants $\varepsilon_* $  and $T_0$ so that for all $a_1 \leq a_2 < \varepsilon_*, t_1 > T_0$, and $|t_2 - t_1| \leq \tau$,
	$$C^{-1} \leq \frac{\mathcal{L}_{a_1}(t_1)}{\mathcal{L}_{a_2}(t_2)} \leq  C \frac{\log(a_1)}{\log(a_2)} . $$
\end{proposition}

\begin{proof}
	What happens is that $\mathcal{L}_{a_1}(t_1)$ and $\mathcal{L}_{a_2}(t_2)$ are approximately the same size until $\mathcal{L}_{a_2}(t_2)$ approaches its supremum, which is around $\log a_2$.
	Then $\mathcal{L}_{a_1}(t_1)$ keeps growing until it approaches its supremum, which is around $\log a_1$.
	Decreasing $\varepsilon_*$ and increasing $T_0$ makes these estimates become arbitrarily multiplicatively accurate.
	
	Let us now be precise.
	Let $\tau > 0$ and $\delta_* \in (0,1]$ be fixed.
	We begin by arranging $T_0$ and $\varepsilon_*$ so that the necessary estimates can be made below.
	So take $T_0$ and $\varepsilon_*$ so that with the constants suppiled by Lemma \ref{lemma-excursion-ests} with $\mathcal{O} = 2 \tau$:
	\begin{enumerate}
		\item \label{assumption-T0-big-init} $C^{-1/6} \log t < \log t - \max( \useC{C-excursion-error}, \useC{good-len-approx}/2)$ and $C^{1/6}\log t > \log t +  \max( \useC{C-excursion-error}, \useC{good-len-approx}/2)$ for all $t \geq T_0 - \tau$,
		\item \label{assumption-adjust-fine} $C^{-1/6} \log t \leq  \log(t - \tau) \leq \log(t + \tau) \leq C^{1/6} \log t$ for all $t > T_0 - \tau$,
		\item \label{assumption-swap-R} $C^{-1/6} |\log a| \leq R_a \leq C^{1/6}| \log a|$ for all $0 < a < \varepsilon_*$,
		\item \label{assumption-monotonicity-swap}  $(R_a + x)/R_a \in [C^{-1/6}, C^{1/6}]$ for all $0 < a < \varepsilon_*$ and $|x| < 2\delta_\bbH$,
		\item \label{assumption-yet-another-one} $(|\log a| + x)/|\log a| \in [C^{-1/6}, C^{1/6}]$ for all $0 < a < \varepsilon_*$ and $|x| < \log (4\delta_\bbH)$. 
	\end{enumerate}
	
	Now fix $a_1 \leq a_2 < \varepsilon_*$ and let $t_1, t_2$ satisfy the hypotheses.
	Note first that that using (\ref{assumption-T0-big-init}), Lemma \ref{lemma-excursion-ests}, and monotonicity of $\mathcal{E}(a, t)$ in $t$ shows that for $i = 1,2$,
	\begin{align}
		C^{-1/6}2\mathcal{E}(a_i,t_i) \leq \mathcal{L}_{a_i}(t_i) \leq C^{1/6}2\mathcal{E}(a_i,t). \label{eq-init-comparison}
	\end{align}   
	
	We now do case analysis on the size of $t_1$ to apply the various estimates in the regimes of behaviors in Lemma \ref{lemma-excursion-ests} to $\mathcal{E}(a,t)$.
	Suppose first that $t_1 \leq 4\delta_\bbH/a_2 + \tau$.
	Then using (\ref{assumption-T0-big-init}) above and Lemma \ref{lemma-excursion-ests}.(\ref{claim-E-approx}) shows for $i = 1,2$,
	\begin{align*}
		C^{-1/6} \log (t_i) \leq \log(t_i) - \useC{C-excursion-error} \leq \mathcal{E}(a_i,t_i) \leq \log(t_i) + \useC{C-excursion-error}  \leq C^{1/6} \log(t_i).
	\end{align*}
	Combining this with (\ref{assumption-adjust-fine}) and Eq.~ (\ref{eq-init-comparison}), in this time range
	\begin{align*}
		C^{-1} \leq C^{-2/3} \frac{\log(t_1) }{ \log(t_2)} \leq \frac{\mathcal{L}_{a_1}(t_1)}{\mathcal{L}_{a_2}(t_2)} \leq C^{2/3}  \frac{\log(t_1) }{ \log(t_2)} \leq C.
	\end{align*}
	
	Next suppose that $4\delta_\bbH/a_2 + \tau \leq t_1$, so that $t_2 \geq 4\delta_\bbH/a_2$.
	Then by Lemma \ref{lemma-excursion-ests}.(\ref{claim-deep-E}) and assumptions (\ref{assumption-swap-R}) and (\ref{assumption-monotonicity-swap}) we have the following estimate on $\mathcal{E}(a_2, t_2)$: 	\begin{align}
		C^{1/6} |\log a_2|  \geq R_{a_2} \geq \mathcal{E}(a_2,t_2) \geq R_{a_2} -2\delta_\bbH \geq C^{-1/6} R_{a_2} \geq C^{-1/3} |\log a_2|. \label{eq-intermediate-times}
	\end{align}
	
	For $(a_1, t_1)$, we break into two cases.
	If $t_1 \geq 4\delta_\bbH/a_1$, then arguing as in Eq.~(\ref{eq-init-comparison}) above $C^{1/6} |\log a_1 | \geq \mathcal{E}(a_1, t_1) \geq C^{-1/6} |\log a_1|$.
	Next, suppose $t_1 \leq 4\delta_\bbH/a_1$. 
	Then arguing as in Eq.~\ref{eq-intermediate-times}) above and using the constraint $t_1 \geq 4\delta_\bbH/a_2$,
	\begin{align*}
		C^{1/6} |\log a_1| \geq R_{a_1} \geq \mathcal{E}(a_1, t_1) \geq C^{-1/6} \log t_1 \geq C^{-1/6}(|\log a_2| + \log (4 \delta_\bbH)) \geq C^{-1/6} |\log a_2|.
	\end{align*}
	Combining these estimates, in all cases
	\begin{align*}
		C^{-1} = C^{-1} \frac{\log a_2}{\log a_2} \leq \frac{\mathcal{L}_{a_1}(t_1)}{\mathcal{L}_{a_2}(t_2)} \leq C \frac{\log a_1}{\log a_2}.
	\end{align*}
\end{proof}

\subsubsection{Crossing Estimates}
We now estimate the length of geodesics in $\mathcal{A}_a'$ that cross through the cylinder.
The estimates are much simpler than in the non-crossing case.

In this paragraph, change notation so that $\alpha_{a,t}$ is a geodesic in $\mathcal{A}_a'$ between the boundary components of basic rotation $t$.
In analogy to before, define
\begin{align}
	\mathcal{L}'_a(t) = \ell(\alpha_{a,t}) - ta.
\end{align}

\newC{C-crossing}
In this setting, something stronger than Prop.~\ref{prop-residue-bound} is true:
\begin{lemma}\label{lemma-crossing-estimate}
	There is a constant $\useC{C-crossing}$ so that $2R_a - \useC{C-crossing} \leq  \mathcal{L}'_a(t) \leq 2R_a$ for all $t \geq 0$ and $\varepsilon > 0$.
\end{lemma}

\begin{proof}
	Once again lift to the universal covering of $\mathcal{A}'_a$ and examine a lift $\widetilde{\alpha}_{a,t}$ of $\alpha_{a,t}$ in $\mathbb{H}^2$.
	Let $\beta_{a,t}$ be a segment of the lift of the core curve of length $ta$, with orthogonal geodesic segments $\sigma_i, \sigma_f$ of length $2R_a$ at its endpoints.
	Then $\widetilde{\alpha}_{a,t}$ connects one side of $\sigma_i$ to the other side of $\sigma_f$, and by symmetry meets $\beta_{a,t}$ at its midpoint $b$.
	The curve $\widetilde{\alpha}_{a,t}$ is the union of the long sides of two right geodesic triangles formed with halves of $\sigma_i, \sigma_f$, halves of $\beta_{a,t}$, and halves of $\widetilde{\alpha}_{a,t}$.
	
	Because the long side of a right triangle in $\bbH^2$ is uniformly close in length to the sum of the other sides, we obtain a bound of $\ell(\alpha_{a,t})\geq 2R_a + ta -\useC{C-crossing}$ for a uniform $\useC{C-crossing}$.
	Using the curve that goes down along $\sigma_i$ to $\beta_{a,t}$ then across $\beta_{a,t}$, then down $\sigma_f$ as a competitor for $\widetilde{\alpha}_{a,t}$ gives the complementary bound $\ell(\alpha_{a,t}) \leq 2R_a + ta$.  
\end{proof}

\section{Length Combinatorialization}\label{s-len-combo}

We now turn to combinatorializing the lengths of curves.
The final statement is Prop.~\ref{prop-combo}, which has slightly more information than Thm.~\ref{theo-distortion}.
We begin by introducing the families of hyperbolic surfaces and decompositions we consider (\S\ref{ss-combo-basics}).
We then describe appropriately invariant combinatorial data induced by isotopy classes of curves with respect to this decomposition, and show that these combinatorial data determine the basic shape of a closed geodesic in the family of surfaces we consider to within certain bounds (\S\ref{ss-combo-combo}).
Finally, we prove our desired combinatorialization (\S\ref{ss-combo-length-ests}).

\subsection{Basic Constructions}\label{ss-combo-basics}

We begin by setting notation for cutting surfaces into hexagons, and giving reminders on how to systematically encode geodesics in pairs of pants.

\subsubsection{Hexagon Decompositions}
We recall the relevant systems of curves.
Let us begin by setting terminology on pants decompositions. Let $\mathscr{P} = \{\gamma_1, ..., \gamma_{3g-3}\}$ be the curves in pants decomposition of $S$.
Write the corresponding pants as $\mathcal{P} = \{P_1, ..., P_{2g-2}\}$.

\begin{definition}[Untwisted Locus]
	Let $\mathscr{U}(\mathscr{P})$ be the subset of $\mathcal{T}(S)$ consisting of hyperbolic structures that have $0$ twist parameters along all curves in $\mathscr{P}$.
	
	For $\varepsilon > 0$ let $\mathscr{U_\varepsilon(P)}$ be the collection of $X \in \mathscr{U(P)}$ so that all length parameters are less than $\varepsilon$, and for $\varepsilon > 0$ let $\mathcal{P}_\varepsilon$ be the collection of (marked) isometry classes of hyperbolic pants with geodesic boundaries of lengths not exceeding $\varepsilon$.
\end{definition}

Each pair of pants $P_j$ has three distinct boundary curves and is uniquely split by three geodesic arcs $\nu_{i,j}$ ($i=1,2,3$) orthogonal to the boundary components into two right hyperbolic hexagons $H^\pm_j$.
For $X \in \mathscr{U}(\mathscr{P})$, the curves $\nu_{i,j}$ glue together into simple closed geodesics $\nu_1, ...,\nu_l$ on $X$.
Denote this collection of simple geodesics by $\mathscr{H}$.

 The complementary regions for $\mathscr{P} \cup \mathscr{H}$ in general $X$ are hexagons; for $X \in \mathscr{U}(\mathscr{P})$ the complementary regions are right hexagons.
 Accordingly, let us define:
\begin{definition}
	For $X \in \mathscr{U(\mathscr{P})}$, the pair $(\mathscr{P}, \mathscr{H})$ is a {\rm{hexagon system}} on $X$.
\end{definition}

To simplify statements of invariance properties of combinatorial data induced by a hexagon system, it is useful to restrict to pants decompositions with the following property.

\begin{definition}\label{def-convenient}
	A pants decomposition $\mathscr{P}$ is {\rm{convenient}} if it is embedded and every curve $\nu \in \mathscr{H}$ has connected intersection with each pair of pants.
\end{definition}

It is an exercise to show that convenient pants decompositions exist.
We henceforth only consider convenient pants decompositions.

We shall often fix a $\delta_*>0$ chosen to be sufficiently small that the following holds.
For any $P \in \mathcal{P}_{\delta_*}$, the hypercycle arcs of length $\delta_*$ around each boundary component of $P$ are embedded, and cutting off the boundary components along these arcs leaves a topological pair of pants $P_{\delta_*}$.
Following Part I, we call $P_{\delta_*}$ the \textit{$\delta_*$-shorts} in $P$.

\subsubsection{Geodesics in pairs of pants}
We have found the technical framework used by Chas-McMullen-Phillips in \cite{chasMcMullenPhillips2019almost} useful for encoding closed geodesics on a pair of pants, and describe it in this paragraph.
The case of non-closed geodesics requires only small modification.

Let $\widetilde{G} = \langle x, y,z \mid x^2 = y^2 = z^2 = e \rangle$, and let $G$ be the subgroup generated by $(a,b,c) = (xy, yz, zx)$.
It is a free group on two generators, with presentation $\langle a,b,c\mid abc =e\rangle$.

Every word in $\widetilde{G}$ may be uniquely put into \textit{reduced form}, i.e. with no letter appearing twice consecutively.
For $w \in \widetilde{G}$ (or $G$), write $\len(w)$ as the length of the reduced form of $w$ in $\widetilde{G}$.
A word is \textit{cyclically reduced} if $w$ has different starting and ending letters.
One says that two cyclically reduced words $w_1, w_2$ are \textit{equivalent} if $w_1$ is conjugate to $w_2$ and in this case writes $w_1 \sim w_2$.
A \textit{run of length $n$} in a word $w \in G$ is an appearance of $(xy)^{\pm n}$, $(yz)^{\pm n}$ or $(zx)^{\pm n}$ in a reduced form of $w$ (or cyclically reduced form, depending on the context). 

Let $P$ be a pair of pants.
Label the three seams of $P$ by $x,y$ and $z$.
Denote by $\mathcal{G}_o$ the collection of closed (not necessarily primitive) oriented geodesics in $P$ that are not peripheral.
Denote by $\mathcal{G}_-$ the collection of geodesic segments between fixed reference points in opposite right hexagons in $P$ and $\mathcal{G}_\upsilon$ the collection of geodesic segments starting and ending at fixed reference points in the same hexagon.
By writing down the ordered intersections of such segments with seams, we get maps from these collections of geodesics to reduced words in the cases of $\mathcal{G}_\upsilon$ and $\mathcal{G}_-$ and to equivalence classes of cyclically reduced words in the case of $\mathcal{G}_o$.
It follows from basic hyperbolic geometry (lift to the universal cover) that these maps are bijections after basic constraints are adjoined:
\begin{align*}
	\mathcal{G}_o &\cong  \left\{ {\text{cyclically reduced } w \in \widetilde{G}} \text{ so }w\not\sim a^n,b^n,c^n\,\, (n \in \bbZ)  \right\} /\sim, \\
	\mathcal{G}_{\upsilon} &\cong \{\text{reduced } w \in \widetilde{G} \text{ with } \mathrm{len}(w) \equiv 0 \pmod2 \}, \\
	\mathcal{G}_- & \cong \{\text{reduced } w \in \widetilde{G} \text{ with } \mathrm{len}(w) \equiv 1 \pmod 2 \}.
\end{align*}

\subsection{Combinatorial Data from Curves}\label{ss-combo-combo}
Our task in this subsection is to give concrete senses in which geodesics look similar in their realizations across similarly shaped hyperbolic surfaces.
The source of invariance that we use is the existence of boundary conjugations of $\partial \bbH^2$ between Fuchsian groups.
In analogy to the classification of simple closed curves \cite{fathiLaudenbachPoenaru2012thurston}, our data for closed curves is obtained in terms of intersections with geodesics adapted to a pants decomposition.

The basic initial observation is:

\begin{lemma}[Simple Intersection Patterns]\label{lemma-simple-inter-pattern}
	Let $\mathscr{S} = \{\gamma_1, ...,\gamma_n\}$ be a finite collection of disjointly realizable simple closed curves and let $\gamma$ be an oriented closed curve in $S$.
	Let $X_1$ and $X_2$ be two marked hyperbolic structures on $S$ and let $\gamma_i^{X_j}$ and $\gamma^{X_j}$ be the geodesic realizations of these curves on $X_j$ ($j = 1,2$).
	Let $(c_k)^{X_j}_{k\in I_j}$ (with $I_j = \emptyset$ or $\mathbb{Z}$) be the sequence of indices of $\gamma_{c_k}^{X_j}$ that $\gamma^{X_j}$ intersects in order, with respect to a base-point on $\gamma^{X_j}$.
	
	Then there is an order-preserving bijection $I_1 \to I_2$ so that after applying this bijection, $(c_k)^{X_1} = (c_k)^{X_2}$. 
\end{lemma}

\begin{proof}
	Fix $j \in \{1,2\}$ and take geodesic representatives on $X_j$, then pick a lift of $\gamma^{X_j}$ to $\mathbb{H}^2$ with endpoints $(x_i, x_f) \in \partial \bbH^2$.
	Let $x_f$ be the positive-time endpoint of $\gamma^{X_j}$.
	One obtains orientations on each connected component $U_1, U_2$ of $\partial \bbH^2 - \{x_i,x_f\}$ by pointing the sides towards $x_f$.
	Note that a geodesic $\eta$ in $\bbH^2$ intersects $\gamma^{X_j}$ if and only if the endpoints of $\eta$ are in opposite connected components of $\partial \bbH^2 - \{x_i,x_f\}$.
	
	Let ${\mathscr{S}}_j'$ be the collection of lifts of elements of $\mathscr{S}_i$ in $X_j$ to $\bbH^2$ whose endpoints are on different connected components of $\partial \bbH^2 - \{x_i,x_f\}$.
	Note that ${\mathscr{S}}_j'$ is countable and obtains a total order from the orientation of $U_i$.
	This order separates points because $\mathscr{S}$ consists of closed curves and $S$ is closed.
	That the geodesics in $\mathscr{S}$ are disjointly realizable implies that this ordering is independent of the choice of connected component $U_i$.
	
	Note that this construction is independent of the choice of lift of $\gamma^{X_j}$ in the sense that changing the lift changes all relevant constructions by a covering transformation of $\bbH^2$, which induces an order-preserving bijection between the relevant sets ${\mathscr{S}}'_j$.
	
	Now if $i \neq j$ there is an equivariant homeomorphism $\xi: \partial \bbH^2 \to \partial \bbH^2$ conjugating uniformizing groups for $X_j$ and $X_i$ and taking the endpoints of the lift of $\gamma^{X_j}$ to the endpoints of a lift of $\gamma^{X_i}$.
	As the boundary map is a homeomorphism and takes endpoints of lifts of elements of $\mathscr{S}$ with respect to $X_j$ to endpoints of lifts of elements of $\mathscr{S}$ with respect to $X_i$, we obtain an order-preserving bijection between ${\mathscr{S}}_j'$ and ${\mathscr{S}}_i'$.
	One sees that the sequences $(c_k)^{X_j}$ are obtained from ${\mathscr{S}}_{j}'$ by listing the curves that elements of ${\mathscr{S}}_{j}'$ are lifts of in order, so the claim is proved.
\end{proof}

Lemma \ref{lemma-simple-inter-pattern} makes the following well-defined.

\begin{definition}
	Let $\gamma$ be a closed geodesic in $X$.
	For $\mathscr{S} = \{\gamma_1, ..., \gamma_n \}$ a family of disjointly realizable simple closed curves on $S$, let $\mathscr{C}_\mathscr{S}$ denote the circularly ordered list of elements of $\mathscr{S}$ that are obtained by intersections with $\gamma$ on one traversal of $\gamma$.
	Let $\mathscr{C}_\mathscr{S}'$ be the totally ordered collection of elements of $\mathscr{S}$ obtained by traversing $\gamma$ any number of times.
\end{definition}

Let us conclude by setting notation for intersection numbers (e.g. \cite{fathiLaudenbachPoenaru2012thurston}).

\begin{definition}
	For $\gamma$ an isotopy class of closed curves and $\mathscr{S}$ a finite collection of simple closed curves in $S$, the {\rm{intersection number}} $i(\gamma, \mathscr{S})$ is $\# \mathscr{C_S}$.
\end{definition}

\subsubsection{Crossing Data for Hexagons}\label{sss-crossing-data}
Fix a hexagon system $(\mathscr{P}, \mathscr{H})$ on $S$ and a hyperbolic structure $X \in \mathcal{T}(S)$, and let $\gamma$ be an oriented closed geodesic in $X$.
Note that no curves in $\mathscr{C'_P}$ or $\mathscr{C'_H}$ share any endpoint, since all involved curves are simple, closed, and distinct.
Because curves in $\mathscr{P}$ and $\mathscr{H}$ intersect each other, there is not a canonical circular order on $\mathscr{C_P} \cup \mathscr{C_H}$.
However, we will see that the failure of Lemma \ref{lemma-simple-inter-pattern} can be well-understood in this setting, and that this failure encodes structure of $\gamma$.

To begin, fix $X_1, X_2 \in \mathcal{T}(S)$ and a lift of $\gamma \subset X_1$ to $\mathbb{H}^2$.
We abuse notation and also denote this lift by $\gamma$.
Let $x_i,x_f$ be its endpoints.
From the orientation of $S$ there is an associated orientation of $\partial \bbH^2$ that is preserved by the relevant equivariant map $\xi: \partial \bbH^2 \to \partial \bbH^2$.
Let us disambiguate connected components $U_1$ and $U_2$ of $\partial \bbH^2 - \{x_i,x_f\}$ so that the orientation of $U_1$ pointing from $x_i$ to $x_f$ agrees with the background orientation of $\partial \bbH^2$, and the two relevant orientations on $U_2$ disagree.

\begin{definition}
	Fix $x_0 \in \gamma$.
	For $i = 1, 2$, let $\iota_i:  \mathscr{C'_P} \cup \mathscr{C}'_H \to \bbZ$ be the bijection obtained by the inclusions of $\mathscr{C'_P}$ and $\mathscr{C'_H}$ into $U_i$ given by taking endpoints.
	Let $\prec_i$ denote the order on $\mathscr{C'_P} \cup \mathscr{C}'_H $ induced by $\iota_i$.
\end{definition}

This structure is well-defined:

\begin{lemma} The orders $\prec_i$ induced on $\mathscr{C'_P} \cup \mathscr{C'_H}$ by $\iota_1$ and $\iota_2$ are well-defined independent of the reference marked hyperbolic structure $X \in \mathcal{T}(S)$ and base-point $x_0 \in \gamma$. 
\end{lemma}

\begin{proof}
	Independence of the choice of base-point $x_0 \in \gamma$ is clear from construction, while independence of the hyperbolic structure comes from the equivariant map $\xi : \partial \bbH^2 \to \partial \bbH^2$ conjugating the uniformizing group of $X_1$ to that of $X_2$ being an orientation-preserving homeomorphism.
\end{proof}

These orders encode a great deal of data.
A first useful observation is that the map $\iota_1 \circ \iota_2^{-1}$ is order-preserving on a pair $\eta, \nu \in \mathscr{C'_P} \cup \mathscr{C'_H}$ if and only if $\nu$ and $\eta$ are disjoint.
This follows immediately from examining crossings in $\bbH^2$.
To record it: 
\begin{lemma}\label{lemma-crossing-swap}
	Let $\nu \in \mathscr{C_P}'$ and $\eta \in \mathscr{C_H}'$.
	Then $\eta$ and $\nu$ intersect if and only if $\eta \prec_i \nu$ and $\nu \prec_j \eta$ for some $i \neq j$.
\end{lemma}

A corollary, which follows from the non-existence of geodesic bigons in $\bbH^2$, is:
\begin{corollary}\label{cor-actual-inters}
	If $\eta, \nu \in \mathscr{C'_P} \cup \mathscr{C'_H}$ and $\eta \prec_i \nu$ for $i = 1, 2$, then for all $X \in \mathcal{T}(S)$, the intersection in $X$ of $\gamma$ corresponding to $\eta$ appears before that corresponding to $\nu$.
\end{corollary}

Directions of crossings are also encoded.
Every curve $\gamma_i \in \mathscr{P}$ bounds two pairs of pants; the orientation of $\gamma_i$ distinguishes the position of lifts of these pairs of pants to each side of $\gamma_i$ in $\bbH^2$.

\begin{definition}
	For a crossing $c \in \mathscr{C'_P}$ corresponding to a curve $\gamma_c$, the geodesic $\gamma$ crosses $\gamma_c$ {\rm{upwards}} if the quadruple $(\gamma_c^+,\gamma^+,\gamma_c^-,\gamma^-)$ is positively oriented in $\partial\bbH^2$ and {\rm{downwards}} if $(\gamma_c^+,\gamma^+,\gamma_c^-,\gamma^-)$ is negatively oriented.
\end{definition}

This is well-defined independent of $X \in \mathcal{T}(S)$, and encodes the order in which $\gamma$ traverses pairs of pants in $\mathcal{P}$.
So, the following is well-defined.

\begin{lemma}[Pants Sequence]
	Write $\mathscr{C_P'} = (c_1, ..., c_k)$.
	There is a sequence $(P_1, ..., P_k)$ of pants in $\mathcal{P}$ so that for any $X \in \mathcal{T}(S)$, the geodesic realization of $\gamma$ is entirely in $P_i$ in between the crossings $c_i$ and $c_{i+1}$. Denote this sequence by $\mathcal{PS}(\gamma)$.
\end{lemma}

The following dichotomy is quite useful:

\begin{definition}
	$d \in \mathscr{C'_H}$ is {\rm{crossing}} if there exists $c \in \mathscr{C'_P}$ and $ i \neq j$ so that $d \prec_i c$ and $c \prec_j d$.
	Otherwise, $d$ is {\rm{internal}}.
\end{definition}

\begin{definition}
	Let $\mathscr{C_P} = (c_1, ..., c_k)$.
	For all $i = 0,...,n-1$ let $s_i$ be the sequence of internal $d \in \mathscr{C'_H}$ so $c_i \prec_k d \prec_k c_{i+1}$ for all $k \in \{1,2\}$.
	Call $s_0, ..., s_{n-1}$ the {\rm{internal words}} of $\gamma$.
\end{definition}

\subsubsection{Twists around Pants Curves}\label{sss-twisting-data}
We now explain how under the additional hypothesis that $X \in \mathscr{U}(\mathscr{P})$, the orders $\prec_i$ encode how many times $\gamma$ wraps around a geodesic it crosses, up to uniform error.

\begin{definition}[Rotation about a Crossing]\label{def-crossing-rot}
	Let $\gamma : \bbR \to X$ be a geodesic in a hyperbolic surface $X$, let $\beta$ be a simple closed geodesic in $X$, and let $\gamma$ transversely intersect $\beta$ at $t \in \bbR$. Denote this crossing by $c$.
	
	Let $\mathcal{A}(\beta)'$ be the annular cover of $S$ with an embedded core geodesic $\widetilde{\beta}$ projecting to $\beta$ and let $\widetilde{\gamma}$ be the lift of $\gamma$ to $\mathcal{A}(\beta)'$ that places $\gamma(t)$ on $\widetilde{\beta}$.
	Let $\sigma$ be a geodesic in $\mathcal{A}(\beta)'$ that is perpendicular to $\widetilde{\beta}$.
	
	The {\rm{unsigned synthetic rotation}} ${|r|}(\gamma,c)$ of $\gamma$ about the crossing $c$ is $|r|_\sigma(\widetilde{\gamma})$, evaluated on this annulus.
	If $\beta$ is oriented, define the {\rm{synthetic rotation}} $r(\gamma,c)$ to be $|r|(\gamma,c)$ if the projection of $\gamma$ to $\widetilde{\beta}$ agrees with the orientation of $\widetilde{\beta}$ and $-|r|(\gamma,c)$ otherwise.
	
	Both $r(\gamma, c)$ and $|r|(\gamma,c)$ are well-defined up to addition of $\pm 1/2$.
\end{definition}

When the hyperbolic surface is not clear from context and $\gamma$ is a closed geodesic, we write the synthetic rotation of an intersection by ${|r|}(\gamma, c,X)$. 
Note that ${|r|}(\gamma,c)$ is finite for all transverse crossings, as a consequence of Lemma \ref{lemma-twist-bd}.
Indeed, the basic rotation of any geodesic segment entirely on the outside of a neighborhood of fixed radius of the core geodesic in $\mathcal{A}_a'$ is uniformly bounded.

Note that by our assumption that $\mathscr{P}$ is convenient, for $X \in \mathscr{U}(\mathscr{P})$, every curve $\gamma_i \in \mathscr{P}$ has exactly two transverse curves $x,y$ in $\mathscr{H}$ that meet $\gamma_i$ transversely and exactly a half-rotation apart.
These curves then distinguish a choice of orthogonal geodesic in Definition \ref{def-crossing-rot}, and so remove the arbitrary choice in the definition of $|r|(\gamma,c,X)$.

\begin{convention}
	We always make this choice for $X \in \mathscr{U(P)}$, so that $|r|(\gamma, c, X)$ is well-defined on $\mathscr{U(P)}$.
\end{convention}

	\begin{figure}
	\includegraphics[scale=0.43]{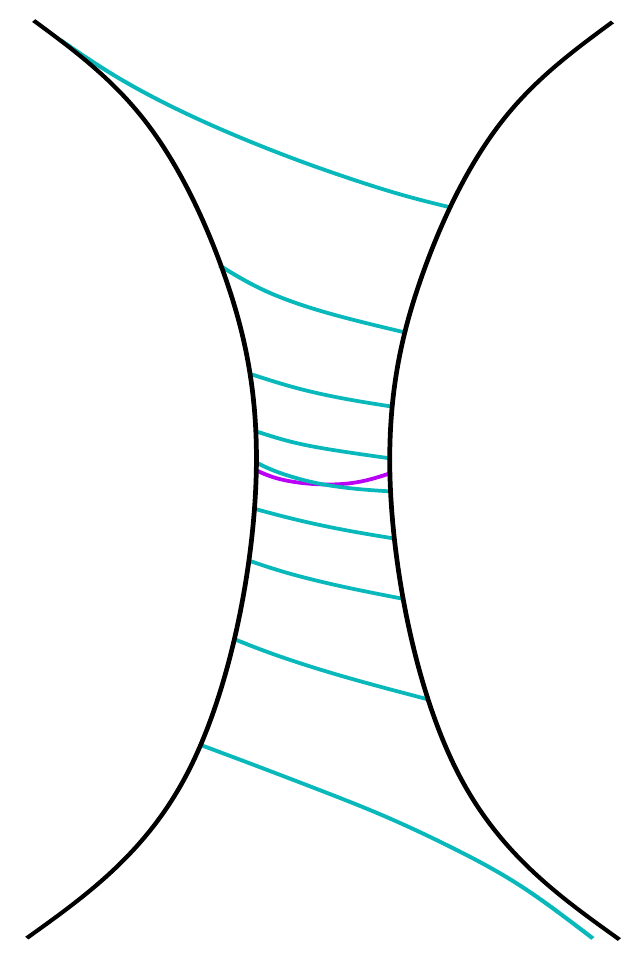} \hspace{1.2cm} \includegraphics[scale=0.43]{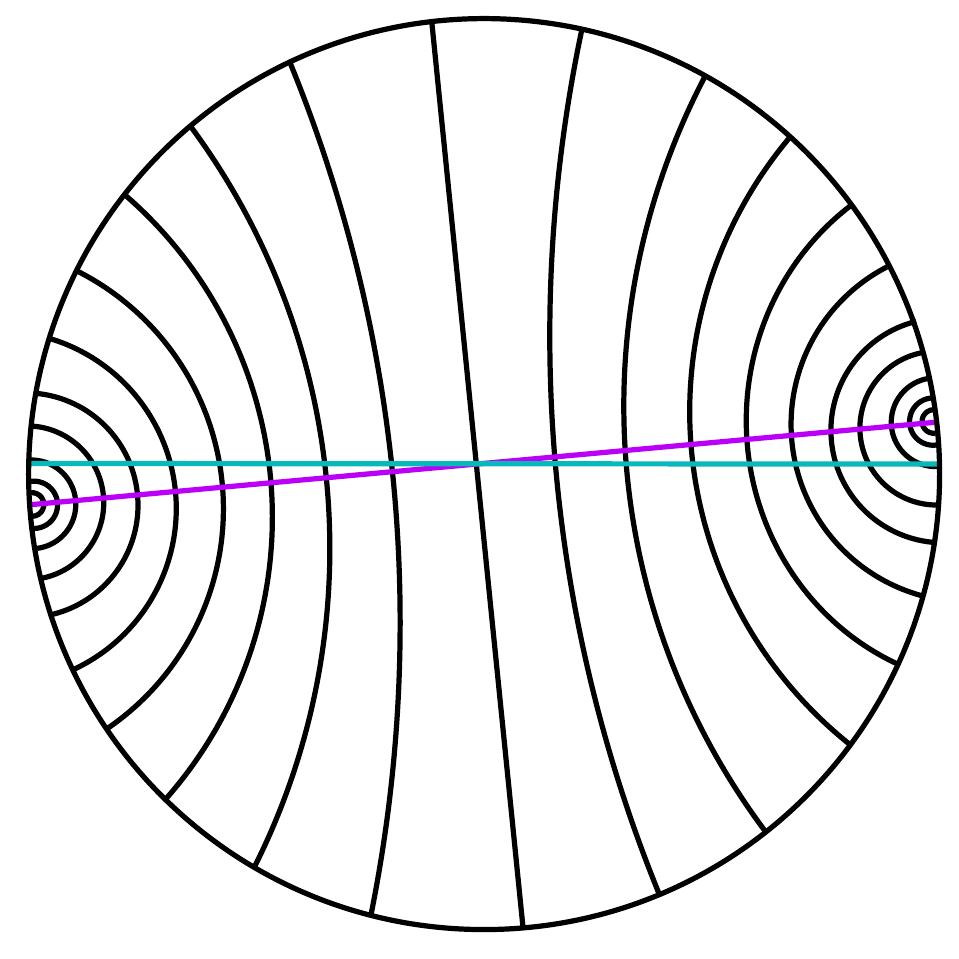}
	\caption{Twisting in a cylinder and its lift to $\bbH^2$.} \label{fig-twist}
\end{figure}

This choice allows for a clean characterization of synthetic rotation in terms of $\prec_i$.

\begin{lemma}[Order and Rotations]\label{order-rotations}
	Let $X \in \mathscr{U(P)}$ and let $\gamma$ be a closed oriented geodesic in $X$.
	For every $c \in \mathscr{C'_P}$, let $n_c(\gamma)$ be half the maximum number of $d \in \mathscr{C_H'}$ so that there is a choice of $i,j$ so that $c\prec_i d$ and $d \prec_j c$.
	Then $|{r}|(\gamma, c, X) = n_c(\gamma)$.
\end{lemma}

See Figure \ref{fig-twist}.

\begin{proof}
	Let $\gamma_c$ be the pants curve corresponding to $c$.
	Begin by lifting $\gamma$ to the annular covering of $S$ corresponding to $\gamma_c$.
	Then $\mathscr{C_H}'$ is in canonical order-preserving bijection with the list of intersections of $\gamma$ with lifts of elements of $\mathscr{H}$ as they occur in this covering.
	From the definition of synthetic rotation, the value $l = |r|(\gamma, c,X)$ is exactly half the number of elements $\nu_1, ..., \nu_k$ of $\mathscr{C_H}'$ whose representatives intersect the core curve $\gamma_c$ in this covering.
	After further lifting to $\bbH^2$, these are exactly the elements of $\mathscr{C_H}'$ that intersect the corresponding lift of $\gamma_c$.
	The Lemma now follows from Lemma \ref{lemma-crossing-swap}.
\end{proof}

We next note that the sign of the synthetic rotation is also well-defined on $\mathscr{U}(\mathscr{P})$.
To see this, as a convention begin by orienting the curves $\eta, \nu \in \mathscr{H}$ that cross $\gamma_i$ so that they cross $\gamma_i$ upwards.
After lifting to the universal cover, the triples $(\gamma_i^+,\gamma^+,\eta_j^+)$ and $(\gamma_i^+,\gamma^+,\nu_j^+)$ are either all positively oriented or all negatively oriented for $\nu_j$ and $\eta_j$ lifts of $\nu$ and $\eta$ in $\mathscr{C'_H}$ that intersect the lift of $\gamma_i$.
One sees that ${r}(\gamma,c_i) > 0$ if and only if $|r|(\gamma,c_i) > 0$ and these triples are positively oriented.
In the downward crossing case, ${r}(\gamma,c_i) > 0$ if and only if $|r|(\gamma,c_i) > 0$ and these triples are negatively oriented.

As this condition is independent of $X \in \mathscr{U}(\mathscr{P})$ and $n_c(\gamma)$ depends only on $\prec_1$ and $\prec_2$,
\begin{corollary}[Rotation Well-Definition]\label{lemma-rot-well-def}
	For all $X,Y \in \mathscr{U(P)}$, $\gamma$ a free isotopy class of closed geodesics on $S$ and $c \in \mathscr{C'_P}$ for $\gamma$, we have ${r}(\gamma,c,X) = {r}(\gamma,c,Y)$.
\end{corollary}

\subsubsection{Rotation tracking}\label{sss-rot-track}
With this structure established on rotations about crossings, we next establish basic control on combinatorial rotation (Def.~\ref{def-combo-rot}).

If a pair $w_iw_{i+1}$ of letters in $\mathscr{C_H}'$ and pants curve $\gamma_i$ satisfy condition (\ref{condition-counting}) in the definition of $n(\gamma, \gamma_i)$ (Def.~\ref{def-combo-rot}), we say that $w_i w_{i+1}$ \textit{counts} for $r(\gamma, \gamma_i)$.

\begin{lemma}[Rotation control]\label{lem:rot-control} Let $\gamma$ be a closed curve and write $\mathscr{C_H} = w$ as $w_1 ... w_n$, up to circular permutation.
	\begin{enumerate}
		\item \label{claim-total-rot-well-defd} $r(\gamma, \gamma_i)$ is well-defined independent of $X \in \mathcal{T}(S)$,
		\item \label{claim-count-exclusivity} Every pair $w_i w_{i+1}$ counts for at most one $r(\gamma, \gamma_k)$,
		\item \label{claim-count-crossing} If $w_i$ and $w_{i+1}$ are both crossing for $c \in \mathscr{C_P}$ corresponding to $\gamma_k \in \mathscr{P}$, then $w_i w_{i+1}$ counts for $r(\gamma, \gamma_k)$,
		\item \label{claim-count-internal} If $w_i w_{i+1}$ are both internal for the same element of the pants sequence $\mathcal{PS}(\gamma)$, then $w_i w_{i+1}$ counts for $r(\gamma, \gamma_k)$ for the corresponding pants curve $\gamma_k$,
		\item \label{claim-most-rotation-understood} Let $M(\gamma)$ be the number of consecutive pairs $w_i w_{i+1}$ in $w$ that are neither crossing for the same $c \in \mathscr{C_P}$ nor internal to the same element of the pants sequence.
		Then $M(\gamma) \leq 2 i(\gamma, \mathscr{P})$. 
	\end{enumerate}
\end{lemma}

\begin{proof}
	For the first two claims, it is useful to lift to the universal cover, and lift a given pair of consecutive symbols $w_{i}w_{i+1}$ in $\mathscr{C_H}$ to $\mathscr{C_H}'$.
	Let $\nu_i,$ and $\nu_{i+1}$ be the curves in $\mathscr{H}$ corresponding to $w_{i}$ and $w_{i+1}$ and let $\widetilde{\nu}_i$ and $\widetilde{\nu}_{i+1}$ be the corresponding lifts to $\bbH^2$ for $w_i$ and $w_{i+1}$.
	From the definition, $w_i w_{i+1}$ is the lift of a pair that counts for $r(\gamma, \gamma_k)$ if and only if there is a lift of $\gamma_k$ that intersects both $\widetilde{\nu}_i$ and $\widetilde{\nu}_{i+1}$.
	This is encoded in the order of endpoints of the relevant curves in $\partial \Gamma$, and so is invariant of $X \in \mathcal{T}(S)$.
	This proves (\ref{claim-total-rot-well-defd}).
	
	Next, note that convenience of $\mathscr{P}$ implies the following. Let $\eta$ be a lift of a pants curve that intersects $\widetilde{\nu_i}$ and $\widetilde{\nu}_{i+1}$, and let $\eta'$ be a distinct lift of a pants curve that intersects $\widetilde{\nu_i}$ and some lift $\widetilde{\nu}_{i+1}'$ of $\nu_{i+1}$.
	Then there is a lift $\eta''$ of a different pants curve that separates $\eta'$ from $\eta$ and so that $\eta''$ intersects $\widetilde{\nu}_i$ and is disjoint from $\widetilde{\nu}_{i+1}$, so that $\widetilde{\nu}_{i+1}' \neq \widetilde{\nu}_{i+1}$.
	This shows that at most one pants curve has a lift that intersects both $\widetilde{\nu}_i$ and $\widetilde{\nu}_{i+1}$.
	This proves (\ref{claim-count-exclusivity}).
	
	Claims (\ref{claim-count-crossing}) and (\ref{claim-count-internal}) follow directly from definitions, the analysis of Lemma \ref{lemma-rot-well-def}, and Lemma \ref{cor-actual-inters}.	
	For the final Claim (\ref{claim-most-rotation-understood}), the orders $\prec_i$ on $\mathscr{C_P}' \cup \mathscr{C_H}'$ split $w$ into a circularly ordered sequence of possibly empty words that are crossing for $c \in \mathscr{C_P}$ or internal to an element of the pants sequence.
	The pairs of letters under consideration appear only in the transitions between these regions, with at most two such pairs per crossing $c \in \mathscr{C_P}$.
\end{proof}

\subsubsection{Initial Decomposition}\label{sss-init-decomposition}
We continue.
In this paragraph we examine the behavior of geodesics near the end of groups of crossing intersections, and use this structure to partition geodesics in a way that isolates their behavior around crossings.
Let $\gamma$ be a closed geodesic in $X \in \mathscr{U}(\mathscr{P})$.

\begin{lemma}[Transitions]\label{lemma-transitions}
	Let $c \in \mathscr{C_P}'$ be a crossing into $P \in \mathcal{P}$ across $\gamma_k$.
	Let $w = (b_{0}, ..., b_{m}) \subset \mathscr{C_H}$ be the crossing indicies corresponding to $c$.
	Let $x,y$ denote the elements of $\mathscr{H}$ that intersect $\gamma_k$ and let $z \in \mathscr{H}$ denote the curve corresponding to the other seam of $P$.
	
\begin{enumerate}
	\item \label{claim-w-crossing-alternates} Then $w$ is an alternating sequence of $x$ and $y$.
	\item \label{claim-w-crossing-other-letter-end} 	Furthermore, if the internal word $s_i$ corresponding to $c$ is nonempty, then the first letter $b_{m+1}$ of $s_i$ is $z$.
	Symmetrically, if $z'$ is the curve corresponding to the other side of the other pair of pants $P'$ bounded by $\gamma_k$ and $s_{i-1} \neq \emptyset$ then $b_{i_0-1} =z'$. 
\end{enumerate}	
\end{lemma}

\begin{proof}
	Let $w_*$ be a possibly empty word so that $c \prec_i w_*$ and $w_* \prec_j c$ for some $i,j \in \{1,2\}$.
	The elements of $\mathscr{H}$ appearing as letters in $w_*$ must intersect $\gamma_k$ by Lemma \ref{lemma-crossing-swap}.
	So $w_*$ is a sequence of the two letters $x, y$ corresponding to the seams of $P$ that intersect $\gamma_k$; one sees that the letters alternate by examining the annular cover corresponding to $\gamma_k$.
	This proves (\ref{claim-w-crossing-alternates}).
	
	Now suppose $s_i$ is nonempty, let $b_{i_m+1}$ denote the first letter of $s_i$, and let $w$ be the word of crossing indicies for $c$.
	Let $\gamma_c$ denote the part in $P$ beginning with the crossing $c$ of the segment in $\gamma_c$ corresponding to the maximal sequence $w'$ of alternating $x,y$ in $\mathscr{C_H}'$ containing $w$ and not exceeding the greatest element of $s_i$.
	Then $\gamma_c$ is entirely contained in the annulus $A_z = P - z$ because $s_i \prec c'$ for all $c' \succ c$ in $\mathscr{C_P}'$ and all crossings in $\mathscr{H}$ in $\gamma_c$ are of the symbols $x,y$.
	
	Since $A_z$ is an annulus with fundamental group generated by $\gamma_k$, Lemma \ref{lemma-crossing-swap} again shows that all intersections in $w'$ are crossing for $c$.
	So $w'$ is disjoint from $s_i$.
	The definition of $w'$ implies that the first letter of $s_i$ is neither $x$ nor $y$, and proves the desired claim for $b_{i_m+1}$.
	The analogue for $b_{i_0-1}$ is symmetric.
\end{proof}

We now begin adding more geometric constraints to the surfaces under consideration and using pinchedness to constrain the shapes of geodesics.

\begin{definition}
	A thickness parameter $\delta_* >0$ is {\rm{combinatorially isolating}} if there exists a $\delta > \delta_*$ so that:
	\begin{enumerate}
		\item For all $P  \in \mathcal{P}_{\delta_*}$, the $\delta$-shorts in $P$ are an embedded pair of topological pants,
		\item For all $a \leq \delta_*$ and bi-infinite geodesics $\gamma$ in $\mathcal{C}_{a}'$ that intersect $\mathcal{C}_{a, \delta_*}$, the intersection of $\gamma$ with $\mathcal{C}_{a, \delta}$ has basic rotation at least $1$.
	\end{enumerate}
\end{definition}
We emphasize that in this definition, the bi-infinite geodesics do not cross the core of $\mathcal{C}_{a}'$, and the rotation does not have to happen in between the $\delta_0$-thick and $\delta$-thick parts of $\mathcal{C}_{a}'$.
It is an exercise to show that all sufficiently small $\delta_*$ are combinatorially isolating.

Combinatorial isolation can be used to constrain the location of spirals.

\begin{lemma}[Spiral Placements]\label{cor-spiral-placements}
	Write  $\mathscr{C_P} = (c_1, ..., c_k)$, and choose pants $P_{c_i}$ so that $\gamma$ crosses from $P_{c_{i-1}}$ to $P_{c_i}$ at $c_i$.
	For all $i$, let $w_i \subset \mathscr{C_H}$ be the crossing intersections for $c_i$.
	Then:
	\begin{enumerate}
		\item \label{claim-almost-no-crossing-overlaps} For all $i$, $\# (w_i \cap w_{i+1}) \leq 1$,
		\item \label{claim-crossings-happen} For all $i$ and $X \in \mathscr{U(P)}$, the segment $\eta_{c_i}$ of $\gamma$ in $P_{c_i} \cup P_{c_{i-1}}$ corresponding to $c_i$ contains all but possibly the first and last intersections in $w_i$,
		\item \label{claim-thick-spiral-end} Suppose, furthermore, $\delta_*$ is combinatorially isolating, $X \in \mathscr{U_{\delta_*}(P)}$, and $\# w_i \geq 2$.
		Let $\eta_{c_i}^o$ be the maximal segment in $\eta_{c_i}$ that contains the crossing $c$ and has connected intersection with the $\delta_*$-shorts of both $P_{c_i}$ and $P_{c_{i-1}}$.
		Then $\eta_{c_i}^o$ contains the crossings corresponding to the sub-word $w_{i}^o$ of $w_i$ that omits the first and last letters of $w_i$.
	\end{enumerate}	
\end{lemma}

\begin{proof}
	Claim (\ref{claim-almost-no-crossing-overlaps}) is because of Lemma \ref{lemma-transitions}.(\ref{claim-w-crossing-alternates}) and the consequence of convenience that within each individual $P \in \mathcal{P}$, the boundary components are distinguished by the elements of $\mathscr{H}$ that they intersect.
	
	If $l$ is a letter in $w$ so that $c_{i-1} \prec_k l \prec_k c_{i+1}$ for all $k \in \{1,2\}$ then the crossing corresponding to $l$ occurs in $\eta_c$ by Cor.~\ref{cor-actual-inters}.
	So Claim (\ref{claim-crossings-happen}) follows from (\ref{claim-almost-no-crossing-overlaps}).
	
	We now turn to the final claim.
	We assume $\# w_c > 2$ so that the claim is not vacuous, and do case-work on whether or not $s_i$ and $s_{i-1}$ are nonempty.
	
	If $s_{i}$ is nonempty, then the letter $z$ appears after the end of $w_i$ with corresponding crossing in $P_i$ by Lemma \ref{lemma-transitions}.
	Note that this occurs, at the latest, in the first time $\eta_{c_i}$ enters the $\delta_*$-collars of $P_i$ from $P_{i,\delta_*}$ by combinatorial isolation and the characterization of crossing words in Lemma \ref{lemma-transitions}.
	And so $\eta_{c_i}^o$ ends in this direction after the end of $w_i^o$.
	The case in which $s_{i-1}$ is nonempty is similar.

	For the other case, suppose $s_i$ is empty.
	Note that this implies $\eta_c$ exits $P_i$ from a different boundary component than the curve $\gamma_i$ it entered across.
	Let $\gamma_e$ denote the exiting pants curve.
	Claim (\ref{claim-crossings-happen}) shows that all but the first and last intersections in $w_i$ occur in $\eta_c$, and the proof showed that the final crossings may fail to occur in $\eta_c$ only if they correspond to an element of $\mathscr{H}$ that intersects $\gamma_e$.
	Without loss of generality, assume that the letters $x,z$ intersect $\gamma_e$ and $x,y$ for $\gamma_i$.
	From combinatorial isolation of $\delta_*$ and characterizations of crossing intersection in Lemma \ref{lemma-transitions}, only at most one appearance of $x$ may appear in the part of $w_i$ coming from the portion of $\eta_{c_i}$ after leaving the $\delta_*$-shorts of $P_i$.
	Because $w_i$ is formed by an alternating sequence of $x$ and $y$, we conclude that both of these sources of ambiguity together only allow at most one letter from the corresponding end of $w_i$ to not appear in $\eta_{c_i}^o$.
	This gives the claim.
\end{proof}

The following decomposition combines what we have proved in this section so far.

\begin{proposition}[Path tracking]\label{lemma-path-tracking}
	Let $\delta_*$ be combinatorially isolating and let $\gamma$ be an oriented closed geodesic in $X \in \mathscr{U}_{\delta_*}(\mathscr{P})$ with $\mathscr{C_P} = (c_1, ...,c_k)$ with rotation numbers $r(\gamma, c_i)$ and internal words $s_0,...,s_{k-1}$.
	
	If $k = 0$, then $\gamma$ is the geodesic in the free homotopy class of $\mathcal{G}$ represented by $s_0$.
	
	If $k \neq0$, then for any $Y \in \mathscr{U}_{\delta_*}(\mathscr{P})$, the geodesic representative of $\gamma$ has a partition into segments $\alpha_1, ..., \alpha_k$ and $\beta_1,..., \beta_k$ so that:
	\begin{enumerate}
		\item \label{claim-twisting-part} For all $i$, the segment $\beta_i$ passes through the crossing $c_i$, has basic rotation number within $1$ of ${r}(\gamma,c_i)$, has connected intersection with $P_{i, \delta_*}$, and the letters in $\mathscr{C_H} \cup \mathscr{C_P} $ occuring in $\beta_i$ are one appearance of $c_i$ and crossing intersections in $\mathscr{C_H}$ for $c_i$.
		\item \label{claim-internal-words-OK} The segment $\alpha_i$ is in $P_i$, begins and ends in $P_{i,\delta_*}$ and has crossings corresponding to the reduced word $s_i \in \widetilde{G}$, adjusted by at most one letter on the start and on the end.
	\end{enumerate}
\end{proposition}

\begin{proof}
	The $k = 0$ case is trivial; assume $k \neq 0$.
	We describe where to split $\gamma$.
	We will give one curve $\alpha_i$ in $P_i$ per crossing $c_i$, and the rest of $\gamma$ will be the union of the segments $\beta_i$.
	
	Fix $i \in \{1, ..., k\}$ and write the word in $\widetilde{G}$ corresponding to $s_i$ by $w = l_1...l_m$.
	If $w$ is empty then the entry and exit curves $\gamma_i$ and $\gamma_{i+1}$ for $\gamma$ in $P_i$ cannot be the same as a consequence of Lemma \ref{lemma-transitions}, and so $\gamma$ has to cross the $\delta_*$-shorts of $P_i$.
	In this case, using Cor. \ref{cor-spiral-placements}.(\ref{claim-thick-spiral-end}), split $\gamma$ while it is in the $\delta_*$-shorts of $P_i$, after at most all but the last letter of the word $w_c$ of crossing intersections for $c$ that are not crossing for $c_{i+1}$.
	Take $\alpha_i$ to be trivial, in this case.
	
	Next, assume $w$ is non-empty.
	By Lemma \ref{lemma-transitions}, the first and last letters $l_1, l_m$ of $w$ correspond to elements of $\mathscr{H}$ that do not meet $\gamma_i$ and $\gamma_{i+1}$, respectively.
	If $\gamma$ leaves the $\delta_*$-shorts in $P_i$ before the corresponding crossing, start $\alpha_i$ as $\eta$ leaves the $\delta_*$-shorts; Lemma \ref{cor-spiral-placements}.(\ref{claim-thick-spiral-end}) ensures that at most one extra crossing is added on to $s_i$ by this choice.
	If $\eta$ does not leave the $\delta_*$-shorts before the corresponding crossing, start $\alpha_i$ briefly before the first crossing appearing in $s_i$.
	The curves $\alpha_i$ are constructed to satisfy property (\ref{claim-internal-words-OK}).
	
	Let $\beta_i$ be the segment in the complement of the union of the $\alpha_i$ containing the crossing $c_i$.
	The curves $\beta_i$ are arranged using Lemma \ref{cor-spiral-placements} to include all but at most $2$ of their intersections in $\mathscr{C_H}$ that are crossing for $c_i$, and this gives Claim (\ref{claim-twisting-part}).
\end{proof}

\subsection{Combinatorialization}\label{ss-combo-length-ests}
We now relate the lengths of the curves appearing in the decomposition of Prop. \ref{lemma-path-tracking} to their combinatorial data.
In the background, fix a combinatorially isolating $\delta_0 > 0$, restrict throughout the remainder of the subsection to $X \in \mathscr{U}_{\delta_0}(\mathscr{P})$, and take decompositions of curves $\gamma$ by applying Prop. \ref{lemma-path-tracking} with this $\delta_0$.
There will be other thickness parameters appearing throughout the subsection, but this $\delta_0$ never changes.

We begin by treating the crossing curves $\beta_i$ before turning to the general combinatorialization (\ref{sss-crossing-combo}).
In handling estimates for internal curves $\alpha_i$ (for $\gamma$ not contained in a single pair of pants), the general estimate borrows from the length gaurantees on the beginning and ending segments $\beta_i$.
Because of this, it is convenient to assemble tools to control the curves $\alpha_i$ (\S \ref{sss-internal-curve-control}) and then prove the general combinatorialization (\S \ref{sss-general-combinatorialization}).

\subsubsection{Crossing Curves}\label{sss-crossing-combo}
We analyze the lengths of the curves $\beta_i$.
Let $\delta_0$ be as above and let $\varepsilon < \delta_0$.

For $X \in \mathscr{U}_{\varepsilon}(\mathscr{P})$ and $P_i \in \mathcal{P}$ pick arbitrary points $x_1, x_2$ in the $\delta_0$-shorts of $P_i$ on each of the two orthogonal hexagons glued together to make $P_i$.
For $w \in \widetilde{G}$ reduced, denote by $\ell(w)$ the length of the unqiue geodesic in $P_i$ with intersection pattern $w$ between $x_1$ and $p \in \{x_1, x_2\}$.
We adopt the convention that whenever $\ell(w)$ is used, points are implicitly fixed so as to make it a well-defined number.

\newD{D-diam-delta}
Let us now introduce some constants.
For a fixed $\varepsilon < \delta_0$, let $\useD{D-diam-delta} < \infty$ be the supremum of the diameter of the $\delta_0$-shorts of pairs of pants $P \in \mathcal{P}_{\varepsilon}$.
Then $\useD{D-diam-delta}$ decreases as we decrease $\varepsilon$, and so a global choice works regardless of if we restrict to pants in $\mathcal{P}_\varepsilon$ for a smaller $\varepsilon$ later.

For any reduced word $w \in \widetilde{G}$ and geodesic segments $\alpha, \alpha'$ starting and ending in $P_{\delta_0}$ representing $w$ in a fixed pair of pants $P \in \mathcal{P}_\varepsilon$,
we have $|\ell(\alpha) - \ell(\alpha')| \leq  2\useD{D-diam-delta}$.

\newe{epsilon-init}
Let us now fix $\varepsilon_{init} < \delta_0$ as an upper bound to the thickness of pants in the following, so that Lemma \ref{lemma-twist-bd} gives a bound $C(\varepsilon_{init}, \delta_0)$.
The way in which we use this is as follows.

\newD{D-thick-tw}
\newC{C-th-tw-bd}
\begin{corollary}\label{cor-thick-tw-length}
	There are constants $\useD{D-thick-tw}$ and $\useC{C-th-tw-bd}$ so that for any $P \in \mathcal{P}_{\varepsilon_{init}}$ and geodesic segment $\gamma$ contained in $P_{\delta_0}$ with reduced word in $w \in \widetilde{G}$ using only two letters, $\ell(\gamma) < \useD{D-thick-tw}$ and $\# w \leq \useC{C-th-tw-bd}$.
\end{corollary}

\begin{proof}
	The hypothesis that $w$ uses only two letters makes $\alpha$ lift to a geodesic segment in a collar $\mathcal{C}_a'$ with $a < \delta_0$ that is entirely contained in the $\delta_0$-thick part of $\mathcal{C}_a'$ and has synthetic rotation at least $(\# w - 2)/2$.
	Lemma \ref{lemma-twist-bd} then places a uniform upper bound $L$ on $\# w$, and containment of $\alpha$ in $P_{\delta_0}$ then forces $\ell(\alpha) \leq (L+1) \useD{D-diam-delta}$. 
\end{proof}
\newC{initial-combo}
The following is the estimate we use on the lengths of $\beta_i$.

\begin{proposition}[Crossing Estimate]\label{prop-initial-combo}
    Let $\varepsilon > 0$ be sufficiently small. Then there is a $\useC{initial-combo}$ so that for all geodesics $\gamma$ in $X \in \mathscr{U}_{\varepsilon}(\mathscr{P})$ and $c \in \mathscr{C_P}$ corresponding to a crossing of $\gamma_i$ with $\ell_X(\gamma_i) = a_i$, the curve $\beta_c$ of the decomposition of Prop.~\ref{lemma-path-tracking} has length estimated by
    \begin{align*}
        \big|2|\log a_i| + |r|(\gamma,c_i) a_i  - \ell(\beta_c)\big|\leq \useC{initial-combo}.
    \end{align*}
\end{proposition}

The point here is that $\useC{initial-combo}$ is independent of the curve $\gamma$ and surface $X \in  \mathscr{U}_{\varepsilon}(\mathscr{P})$.

\begin{proof}
	By Lemma \ref{lemma-path-tracking}.(\ref{claim-twisting-part}), the curve $\beta_c$ crosses $c$, has basic rotation within $1$ of $r(\gamma, c)$, and has connected intersection with $P_{i-1,\delta_0}$ and $P_{i,\delta_0}$.
	Let $\beta_c'$ be the part of $\beta_c$ contained in the $\delta_0$-collar around $\gamma_c$.
	Then Lemma \ref{lemma-twist-bd} and Cor.~\ref{cor-thick-tw-length} show that $|\ell(\beta_c) - \ell(\beta_c')| \leq 2\useD{D-thick-tw}$, and $\beta_c'$ has unsigned total rotation $|r|(\beta_c')$ satisfying $| |r|(\beta_c') - |r|(\gamma, c,X)| \leq \useC{C-th-tw-bd} + 2$.
	Applying Lemma \ref{lemma-crossing-estimate}, taking $\varepsilon < 1$, and noting that $$|R_a - |\log a_i|| = |\arccosh (\delta_0/a_i) - |\log(a_i)|| \leq \log2\delta_0 $$ gives the claim.
\end{proof}

\subsubsection{Controlling Internal Curves}\label{sss-internal-curve-control}
What is left is now to estimate the lengths of the curves $\alpha_i$ appearing in Prop.~\ref{lemma-path-tracking}.
Properly handling excursions in a way that is uniform across curves makes this more intricate than the combinatorialization for $\beta_i$.

The ambuguity introduced by choices of reference points in $\ell(w)$ is encoded in the following initial observation: 
\newC{C-ai-error}

\begin{lemma}\label{lemma-a-init}
	Take the notation of Lemma \ref{lemma-path-tracking}.
	In particular, $\mathscr{C_P} = (c_1, ..., c_k)$ with corresponding internal words $s_1, ..., s_{k}$.
	Then for all $X \in \mathscr{U}_{\delta_0}(\mathscr{P})$, isotopy classes of curves $\gamma$, and $i = 1, ..., k$,
	$|\ell_X(w_i) - \ell(\alpha_i)| \leq 4\useD{D-diam-delta}$.
\end{lemma}

\begin{proof}
		This is an immediate corollary of Prop.~\ref{lemma-path-tracking}: $2\useD{D-diam-delta}$ come from the ambiguity of $\ell(w)$ and $2\useD{D-diam-delta}$ come from the potential difference of up to $2$ letters from $s_i$ and the word corresponding to crossings of $\alpha_i$.
\end{proof}

We now turn to directly constraining the shapes of geodesics in pants.
We will use the following consequence of combinatorial isolation.
\begin{lemma}\label{lemma-transition-crossings}
    Let $\delta_*$ be combinatorially isolating and $P \in \mathcal{P}_{\delta_*}$.
    Let $\gamma$ be a geodesic segment in $P$ with corresponding word $w \in \widetilde{G}$.
    Any sub-segment of $\gamma$ corresponding to a three-letter sub-word $abc$ of $w$ using all letters $x,y,z$ must intersect the $\delta_*$-shorts of $P_i$.
\end{lemma}

\begin{proof}
Suppose without loss of generality that the sub-word is $xyz$.
If the segment of $\alpha$ corresponding to $xy$ does not intersect the $\delta_*$-shorts of $P_i$, then this segment is contained in the annulus $\mathcal{A}_{xy}$ in $P_i - P_{i,\delta_*}$ whose fundamental group is generated by $(xy) \in G$.
Then $\gamma$ has to cross the $\delta_*$-thick part to finish the word $xyz$ as $\mathcal{A}_{xy}$ is disjoint from $z$.
\end{proof}

The key point on the segments $\alpha_i$ is the following lemma, which allows us to break up curves in $P$ into parts that stay in uniformly thick peices, and well-controlled parts that pick up a gauranteed amount of length of our choosing.
We call these length gaurantees \textit{depth parameters} in the below.

\begin{lemma}[Depth Parameters]\label{lemma-internal-split-ests}
    Let $D, T_0 < \infty$ be given.
    Then there is an $\varepsilon_0 > 0$ and constants $T \in (T_0, \infty)$ and $\delta_1 > 0$ so that for every $P \in \mathcal{P}_{\varepsilon_0}$,
    \begin{enumerate}
        \item \label{claim-capris} Every geodesic segment $\alpha$ starting and ending in $P_{\delta_0}$ corresponding to a word $w \in \widetilde{G}$ with no run of length more than $T$ is entirely contained in $P_{\delta_1}$,
        \item \label{claim-undertainty-bound} Any geodesic segment $\alpha$ beginning and ending in $P_{\delta_0}$ corresponding to a run $r$ of length $N \geq T - (2\useC{C-th-tw-bd}+1)$ is the union of two segments of lengths in $[0,\useD{D-thick-tw}]$ and a geodesic $\alpha_{a_i,t}$ in a $\delta_0$-collar $\mathcal{C}_{a_i}$ of a boundary component $\gamma_i$ of $P$ with $N - (2\useC{C-th-tw-bd} +1) \leq t \leq N + 1$,
        \item \label{claim-depth} $\mathcal{L}_{a}(t) \geq D$ for all $a <\varepsilon_0$ and $t \geq T - (2\useC{C-th-tw-bd}+1)$.
    \end{enumerate}
\end{lemma}

\begin{proof}
	We start by taking $\varepsilon_0' < \varepsilon_{init}$, and decrease it later as needed.
	Let $\alpha$ be a geodesic segment beginning and ending in $P_{\delta_0}$ with intersection pattern corresponding to a word $w$.
    We begin by arranging for (\ref{claim-undertainty-bound}) to hold, which does not place much constraint on $T$.
    
    Let $P \in \mathcal{P}_{\varepsilon_0'}$ and let $\alpha$ be a geodesic segment in $P$ starting and ending in $P_{\delta_0}$ with corresponding word $w$ a run of length $N$ on the letters $x,y$.
    Let $a_i$ denote the length of the corresponding boundary component and let $\mathcal{C}_{a_i}$ be the corresponding $\delta_0$-collar.
    We note that for $N > 1$, $\alpha$ cannot intersect the other $\delta_0$-collars in $P$ because of combinatorial isolation of $\delta_0$, and by definition $\alpha$ cannot intersect $z$.
	
	So for $N > 2$, by Cor.~\ref{cor-thick-tw-length}, $\alpha$ consists of two intervals in $P_{\delta_0}$ with corresponding words of length less than $\useC{C-th-tw-bd}$ and segments of lengths in $[0,\useD{D-thick-tw}]$ and a segment $\alpha_{a_i,t}$ in $\mathcal{C}_{a_i,\delta_0}$ with synthetic rotation between $N - 2\useC{C-th-tw-bd}$ and $N$.
	The statement is then well-defined and the desired bound on $t$ then follows for all $T > 2\useC{C-th-tw-bd} + 3$ from the uniform comparability of basic rotation and synthetic rotation.
	So take $T$ satisfying this.
	Note that this argument shows that we may increase $T$ while still maintaining (\ref{claim-undertainty-bound}).
	Do this so that $T > T_0$.
	
	By applying Lemma \ref{lemma-excursion-ests} with $\delta_* = \delta_0$, after further decreasing $\varepsilon_0'$ and increasing $T$, we may arrange for (\ref{claim-depth}) to hold.
	This will be the final value of $T$.
	
	We finally address (\ref{claim-capris}).
	Any segment $\alpha_{a_i, t}$ in a $\delta_0$-collar $\mathcal{C}_{a_i}$ with basic rotation $t$ corresponds to a run in $w$ of length at least $t - 2$, and so we obtain the bound $t \leq T +2$ under the hypothesis of (\ref{claim-capris}).
	Lemma \ref{lemma-excursion-ests} gives a uniform upper bound ${E}$ on $\mathcal{E}(a,t)$ for all $a < \varepsilon_0$ and $t < T + 2$.
	By decreasing $\varepsilon_0'$, we may ensure that $R_a > {E}$ for all $a < \varepsilon_0'$.
	
	Let $\delta_{a, {E}}$ denote the length of the hypercycle in $\mathcal{C}_a$ at distance ${E}$ from the long boundary component.
	Then $\delta_{a,{E}}/\delta_0$ is given by $\cosh(R_a - {E})/\cosh(R_a) > \exp(-{E})/4$.
	After further decreasing $\varepsilon_0'$ to a final value $\varepsilon_0$, arrange for $\varepsilon_0 < \delta_0\exp(-{E})/4$, and take $\delta_1 =\delta_0\exp(-{E})/4$. 
	Then all claims hold with these values of $\varepsilon_0, \delta_1, T$.
\end{proof}

To handle the parts of curves in $\delta_1$-shorts of pants given in the decoposition of Lemma \ref{lemma-internal-split-ests}, we use Lemma \ref{lemma-compact-part-estimate} from Part I.

\begin{lemma}\label{cor-compact-est}
    Let $C > 1$ be given. Let $\varepsilon > 0$ and $T' < \infty$ have the property that any geodesic segment in $P \in \mathcal{P}_{\varepsilon}$ that starts and ends in $P_{\delta_0}$ and has crossing pattern a word with no run of length exceeding $T'$ is contained in $P_{\delta_1}$.
    
    Then there is an $\varepsilon_0 < \varepsilon$ with the additional following properties.
    \begin{enumerate}
    	\item     For all $P_1, P_2 \in \mathcal{P}_{\varepsilon_0}$, words $w \in \widetilde{G}$ with no run of length exceeding $T'$, and geodesics $\alpha_i$ in $P_i$ with crossings given by $w$ starting and ending in $P_{i,\delta_0}$,
    	\begin{align*}
    		C^{-1} \ell_{P_1}(\alpha)  - 6 \useD{D-diam-delta}  \leq \ell_{P_2}(\alpha_2) \leq C \ell_{P_1}(\alpha_1) + 6\useD{D-diam-delta}.
    	\end{align*}
    	Furthermore, $\ell_{P_i}(\alpha_i) \geq \delta_1 (\len(w) -1)$ ($i = 1, 2$).
    	\item \label{claim-closed-detail} For all $P_1, P_2 \in \mathcal{P}_{\varepsilon_0}$ and cyclically reduced words $w\in G$ with no (cyclic) run of length exceeding $T'$ and closed geodesics $\alpha_i$ in $P_i$ with crossings given by $w$, $\ell_{P_1}(\alpha_1) \in [C^{-1} \ell_{P_2}(\alpha _2), C \ell_{P_2}(\alpha_2)]$. Furthermore, $\ell_{P_i}(\alpha_i) \geq \delta_1\len(w)$ ($i = 1,2$).
    \end{enumerate}

\end{lemma}

\begin{proof}
	The length lower bounds are both direct estimates from injectivity radii of $\delta_1$-shorts.
	For the comparisons, we may as well assume $C < 2$ and take $\varepsilon$ and $C$ to be small enough so that for any $P_1 \in \mathcal{P}_\varepsilon$ there is a point $p \in P_1$ so that for any $P_2 \in \mathcal{P}_\varepsilon$ and $C$-Lipschitz map $f: P_1 \to P_2$ satisfying the conclusions of Lemma \ref{lemma-compact-part-estimate}, $f(p)\in P_{2, \delta_0}$.

    Now take $\varepsilon_0 < \varepsilon$ as produced by Lemma \ref{lemma-compact-part-estimate} applied with the given $C$.
    Fix points $p_i \in P_{i,\delta_0}$ so that the Lipschitz maps produced by Lemma \ref{lemma-compact-part-estimate} map $p_i$ into $P_{j,\delta_0}$ for $j \neq i$.
    Modify $\alpha_1$ to start and end at the point $p_1$ in Lemma \ref{lemma-compact-part-estimate}, adding at most $2\useD{D-diam-delta}$.

    Then, modifying the image with two segments of length not exceeding $\useD{D-diam-delta}$ Lemma \ref{lemma-compact-part-estimate} gives a curve in $P_2$ containing a curve of the same isotopy class rel. endpoints as $\alpha_2$ and with the same start and end-points with length at most $C\ell_{P_1}(\alpha_1) + 6\useD{D-diam-delta}$.
    This gives the bound $\ell(\alpha_2) \leq C \ell_{P_1}(\alpha_1) + 6\useD{D-diam-delta}$.
    The other bound follows by applying this bound with the roles of $P_1$ and $P_2$ interchanged.
    
    In the closed curve case, there is no estimation to do on error from choices of endpoints, and the claim follows directly from the properties of the maps produced by Lemma \ref{lemma-compact-part-estimate}.
\end{proof}

Lemma \ref{lemma-internal-split-ests} splits consideration of parts of curves in pants into portions in shorts and portions of quantitative commitment to spiraling around a collar.
Accordingly, it is useful to record the invariant data that encodes when such committed spiralling may occur.

\begin{definition}
    Let $T \geq 2$ in $\frac{1}{2}\mathbb{N}$.
    For a geodesic $\gamma$ in $X \in \mathscr{U(P)}$, let $e_1(\gamma,T), ..., e_m(\gamma,T)$ denote the cyclically ordered maximal runs appearing in the internal words $s_i$ in $\mathscr{C_H}$ of lengths no less than $T$ and let $|r|(\alpha, e_i,T)$ be half their corresponding (word) lengths.
    Call $e_1(\gamma,T), ..., e_{m}(\gamma,T)$ the {\rm{$T$-internal runs}} of $\gamma$.
\end{definition}

As the internal words $s_i$ are invariant of $X \in \mathscr{U(P)}$, we have

\begin{lemma}\label{lemma-internal-runs-well-defd}
    The run sequence $e_1(\gamma,T),...,e_m(\gamma,T)$ and lengths $|r|(\gamma,e_i,T)$ are independent of $X \in \mathscr{U}(\mathscr{P})$. 
\end{lemma}

  Note that the $T$-internal runs $e_{i}(\gamma, T)$ are disjoint except possibly at endpoints for conecutive $e_{i}(\gamma, T), e_{i+1}(\gamma, T)$ contained in the same internal word.

\subsubsection{General Combinatorialization}\label{sss-general-combinatorialization}
We are now ready to prove the main result of this section:
We retain the notations and fixed $\delta_0, \varepsilon_{init}$ of the previous paragraphs in this subsection.

\newC{C-combo-K}
\newC{C-combo-L}
\newC{C-leftovers}
\begin{proposition}\label{prop-combo}
    Fix a convenient pants decomposition $\mathscr{P}$.
    Let $C > 1$ and $B < \infty$ be given.
    Then there exist $\varepsilon >0$ so that the following hold.
        
    \begin{enumerate}
        \item \label{claim-nontwisting-big} For all $X\in \mathscr{U}_{\varepsilon}(\mathscr{P})$ and $\gamma$ neither trivial nor in $\mathscr{P}$, we have $\mathscr{L}(\gamma,X) > 0$ and $\mathscr{L}(\alpha,X) \geq B i(\gamma, \mathscr{P})$.
        \item \label{claim-distortion} For all $X,Y \in \mathscr{U}_{\varepsilon}(\mathscr{P})$ and nontrivial $\gamma$ not contained in $\mathscr{P}$,
        \begin{align}
          \min_{\gamma_i \in \mathscr{P}, \sigma \in \{\pm 1\}} C^{-1} \left(\frac{\log(\ell_X(\gamma_i))}{\log(\ell_Y(\gamma_i))}\right)^{\sigma}  \leq   \frac{\mathscr{L}(\gamma,X)}{\mathscr{L}(\gamma,Y)} \leq \max_{\gamma_i \in \mathscr{P}, \sigma \in \{\pm 1\}} C \left(\frac{\log(\ell_X(\gamma_i))}{\log(\ell_Y(\gamma_i))}\right)^{\sigma}. \label{eq-comparison}
        \end{align}
    \end{enumerate}
\end{proposition}

We remark that the thresholds $\varepsilon$ produced by the proof will be extremely small in practice, because they are chosen so that $-\log \varepsilon$ is quite large.

\begin{proof}
	The proof is by applying estimates we have made to the decomposition of Prop. \ref{lemma-path-tracking}, and by absorbing additive uncertainties in lengths into gaurantees of lengths during runs or crossings.
	 We split into cases on whether or not the curve $\gamma$ is entirely contained in a single $P \in \mathcal{P}$.
	We give a detailed proof for $\gamma$ not contained in a single $P \in \mathcal{P}$.
	The other case is similar but slightly distinct; we conclude by remarking on how the proof is adjusted here.

	We begin.
	Determine constants $D, T, \delta_1$, and $\varepsilon$ by the following procedure:
	\begin{enumerate}
		\item \label{init-consts} Apply Lemma \ref{prop-residue-bound} with fixed value $\delta_*= \delta_0$ to obtain values $T_0$ and $\varepsilon_0$ for which the conclusion of Lemma \ref{prop-residue-bound} holds with values $C^{1/2}$ and $\tau = 4 \useC{C-th-tw-bd} +2$,
		\item \label{giant-D} Let $T, \delta_1, \varepsilon_1$ be obtained by applying Lemma \ref{lemma-internal-split-ests} with inputs $D$ and $T_0$, where $T_0$ is as in (\ref{init-consts}) and $D$ satisfies that $(x + R)/x \in [C^{-1/4}, C^{1/4}]$ for all $x > D$ and $|R| < 2\useD{D-thick-tw} + 6 C^{1/4}\useD{D-diam-delta} + (2\useC{C-th-tw-bd}+1)\delta_0 + \varepsilon_0/2$.
		\item Then $\delta_1, T$, and $\varepsilon_2 = \min \{\varepsilon_0, \varepsilon_1\}$ satisfy the hypotheses of Lemma \ref{cor-compact-est}.
		Apply Lemma \ref{cor-compact-est} with the multiplicative constant $C^{1/2}$ to obtain an $\varepsilon_3 < \varepsilon_2$ for which the conclusion of Lemma \ref{cor-compact-est} holds.
		\item \label{last-assumption} Take $\varepsilon_4 < \varepsilon_3$ to be small enough that furthermore $(|\log a| + R)/|\log a| \in [C^{-1/2}, C^{1/2}]$ for all $a < \varepsilon$, and $|R|< \useC{initial-combo} + (4 + 6C^{1/4}) \useD{D-diam-delta} + \varepsilon_3$,
		\item \label{real-last-assumption} Take $\varepsilon_5 < \varepsilon_4$ so that also $C^{-1/2} |\log a| \geq B/2$ for all $a < \varepsilon_5$,
		\item \label{but-wait-theres-more} Take $\varepsilon < \varepsilon_5$ so that furthermore $1 +x \in [C^{-1/4}, C^{1/4}]$ for all $|x| < \varepsilon/2\delta_1$. 
	\end{enumerate}
	
    Let $X \in \mathscr{U}_{\varepsilon}(\mathscr{P})$ and $\gamma$ be neither trivial, in $\mathscr{P}$, nor contained in an individual $P \in \mathcal{P}$.
    Then $\mathscr{C_P} = (c_1, ..., c_k)$ is nonempty.
    Let $\gamma_{c_i} \in \mathscr{P}$ be the corresponding curves of lengths $a_X(c_i)$ (resp. $a_Y(c_i)$).
    Let $s_1, ..., s_m$ be the words corresponding to internal sequences in $\gamma$.
    The estimates of Prop.~\ref{prop-initial-combo} and Lemma \ref{lemma-a-init} give
    \begin{align}
        \ell_X(\gamma) &= \sum_{i=1}^l (2|\log a_X(c_i)| + |r|(\gamma, c_i) a_X(c_i) + \ell_X(s_i) + A_i) \qquad (|A_i| \leq \useC{initial-combo} + 4 \useD{D-diam-delta}).\label{eq-first-len-decomp}
    \end{align}

   For each $i=1, ..., l$, use the threshold $T$ to split (up to length-$1$ overlaps at endpoints) $s_i$ into words $f_{i, 0}(\gamma,T), e_{i, 1}(\gamma,T), ..., e_{i, k_i}(\gamma,T), f_{i, k_i}(\gamma,T)$ with $e_{i, 1}(\gamma,T), ..., e_{i,k_i}(\gamma,T)$ the $T$-internal runs inside $s_i$ around curves $\gamma(e_{i,j})$ and $f_{i,0}, ..., f_{i, k_i}$ possibly empty maximal words disjoint from and in-between the words $e_{i,j }(\gamma, T)$.
   
   Let $\upsilon_i$ be a geodesic starting and ending in $P_{i,\delta_0}$ and representing the combinatorics of the word $s_i$.
   Then Lemma \ref{lemma-transition-crossings} shows that $\upsilon_i$ can be partitioned into disjoint segments $\eta_{i, 1}, ..., \eta_{i, k_i}$ representing the $T$-internal words $e_{i, j}(\gamma, T)$ in $s_i$ up to potentially cutting off one letter from endpoints and segments $\nu_{i,0}, ..., \nu_{i, k_i}$ representing $f_{i, j}(\gamma, T)$.
    Then, of course,
    $$ \ell_X(s_i) = \sum_{j=1}^{k_i} \ell_X(\eta_{i, j}) + \sum_{j=0}^{k_i} \ell_X(\nu_{i, j}).$$

   From Prop.~\ref{lemma-internal-split-ests},
    \begin{align}
        \sum_{i=1}^l \sum_{j = 1}^{k_i} \ell_X(\eta_{i, j}) = \sum_{i=1}^l \sum_{j = 1}^{k_i} [\mathcal{L}_{\ell_X(\gamma(e_{i,j}))}(t_{i,j}) + |r|(\alpha, e_{i,j},T)\ell_X(\gamma(e_{i, j}))+ C_{i, j}]\label{eq-len-sub-decomp}
    \end{align}
    with $|t_{i,j} - |r|(\alpha, e_{i,j},T)| \leq 2\useC{C-th-tw-bd}+1$ and $|C_{i,j}|\leq 2\useD{D-thick-tw} + (2\useC{C-th-tw-bd}+1)\delta_0$.
   
   	Lemma \ref{lem:rot-control} shows that the contributions to $\mathscr{R}(\gamma,X)$ that do not appear in Eqs.~(\ref{eq-first-len-decomp}) and (\ref{eq-len-sub-decomp}) consist of an amiguity from transitions between $\beta_i$ and $\alpha_i$ of size at most $\varepsilon i(\gamma, \mathscr{P})$, contributions from letters in each $f_{i,j}(\gamma,T)$ of size no more than $\varepsilon(\len(f_{i,j}(\gamma,T)) -1)/2,$ and terms of size at most $\varepsilon/2$ for each transition between $T$-internal runs and internal words $f_{i,j}(\gamma, T)$.
   	Grouping the first contributions and the constants $A_i$ into constants $A_i'$, the second contributions into constants $E_{i,j}$, and the third contributions together with $C_{i,j}$ into constants $C_{i,j}'$, we have
    \begin{align*}
        \mathscr{L}(\gamma, X) = \sum_{i=1}^l 2 [|\log(a_X(\gamma_{c_i}))| + A_i] +  \sum_{i=1}^l \left[ \sum_{j = 1}^{k_i} [\mathcal{L}_{\ell_X(\gamma(e_{i,j}))}(t_{i,j}) + C_{i, j}'] + \sum_{j=0}^{k_i} [\ell_X(\nu_{i,j}) +E_{i,j}] \right]
    \end{align*}
	with $|A_i'| \leq |\useC{initial-combo}| + 4 \useD{D-diam-delta} + \varepsilon$, $|C_{i,j}'| < 2\useD{D-thick-tw} + (2\useC{C-th-tw-bd}+1)\delta_0 + \varepsilon/2$, and $|E_{i,j}| \leq \varepsilon (\len (f_{i,j}(\gamma,T)) - 1)/2$.
    Note that positivity of $\mathscr{L}(\gamma,X)$ follows at this point in the argument.

    Fom (\ref{but-wait-theres-more}) and Lemma \ref{cor-compact-est}, there are constants $J_{i,j} \in [C^{-1/4} ,C^{1/4}]$ so that $\ell_X(\nu_{i,j}) + E_{i,j} = J_{i,j}\ell_X(\nu_{i,j})$.
    Now let $Y \in \mathscr{U}_\varepsilon(\mathscr{P})$ be given, and write the constants corresponding to $J_{i,j}$ by following the same analysis in $Y$ as $J_{i,j}'$.
    Write by $\nu_{i,j}^Y$ the segments in $Y$ corresponding to $\nu_{i, j}$ in $X$.
    Using Lemma \ref{cor-compact-est}, there are numbers $D_{i, j}$ with $|D_{i,j}| \leq C^{1/4}6 \useD{D-diam-delta}$ so that $J_{i,j}\ell_X(\nu_{i,j}) + D_{i,j} =  F_{i,j} \ell_Y(\nu_{i,j}^Y)$ with $F_{i, j} \in [C^{-3/4},C^{3/4}]$.
    
    Next, for each $i = 1, ...,m$, absorb all but $D_{i, 0}$ into the constants $C_{i,j}'$ to obtain new constants $C_{i,j}''$ with $|C_{i,j}''| \leq 2\useD{D-thick-tw} + C^{1/4}6 \useD{D-diam-delta}+\varepsilon/2$.
    Absorb the remaining constants $D_{i,0}$ into the corresponding constants $A_i'$ to obtain new constants $A_i''$ with $|A_i''| \leq  \useC{initial-combo} + (4+6C^{1/4})\useD{D-diam-delta} + \varepsilon$.
    
    Next, note that the twisting terms $t_{i,j}'$ in $Y$ differ from $t_{i,j}$ by no more than $4 \useC{C-th-tw-bd} +2$.
	Absorb the constants $A_i''$ into multiplicative constants $G_i \in [C^{-1/2}, C^{1/2}]$ on the $|\log(a_X(\gamma_{c_i}))|$-terms (resp. $G_i' \in [C^{-1/2}, C^{1/2}]$ on the $|\log(a_Y(\gamma_{c_i}))|$-terms) using (4), and absorb the constants $C_{i,j}''$ into multiplicative constants $H_{i,j} \in [C^{-1/4}, C^{1/4}]$ (resp. $H_{i,j}' \in [C^{-1/4}, C^{1/4}]$) on $\mathcal{L}_{\ell_Y(\gamma_i)}(t_{i,j}')$-terms (resp. $\mathcal{L}_{\ell_X(\gamma_i)}(t_{i,j})$-terms) using (\ref{giant-D}). This gives
    \begin{align}
        \mathscr{L}(\alpha, X) &= \sum_{i=1}^l \left(2 G_i |\log(a_X(\gamma_{c_i}))| +  \sum_{j = 1}^{k_i} H_{i, j} \mathcal{L}_{\ell_X(\gamma_i)}(t_{i,j}) + \sum_{j=0}^{k_i} F_{i,j} \ell_Y(\nu_{i, j}^Y) \right) , \label{eq-length-expanded} \\
        \mathscr{L}(\alpha, Y) & = \sum_{i=1}^l \left( 2 G_i' |\log(a_Y(\gamma_{c_i}))| +  \sum_{j = 1}^{k_i} H_{i,j}' \mathcal{L}_{\ell_Y(\gamma_i)}(t_{i,j}') + \sum_{j=0}^{k_i} J_{i,j}'\ell_Y(\nu_{i,j}^Y)\right) \nonumber
    \end{align}
    
    Condition (\ref{init-consts}) and our bounds on $G_i, G_i', H_{i,j}, H_{i,j}'$, and $F_{i,j}$ ensure that the ratios of each pair of terms in these sums across $X$ and $Y$ under the canonical bijection of entries are bounded by the desired expression in Eq.~ (\ref{eq-comparison}).
    As all terms in each sum are positive, this gives the desired estimate for (\ref{claim-distortion}) in the case where $\gamma$ is not contained in a single pair of pants.
    The estimate in Claim (\ref{claim-nontwisting-big}) also follows immediately from Eq. (\ref{eq-length-expanded}) and condition (\ref{real-last-assumption}).
    
    In the case where $\gamma$ is contained in a single pair of pants, make the following modifications.
    Keep the same constants $D, T, \delta_1$, and $\varepsilon$.
    For $\gamma$ that do not make a run of length $T$, the purely multiplicative estimate in Lemma \ref{cor-compact-est}.(\ref{claim-closed-detail}) gives the claim immediately.
    For $\gamma$ that do make a run of length $T$, use the same decomposition of $\gamma$ induced by the Depth Parameters Lemma \ref{lemma-internal-split-ests}.
    Starting the decomposition at points in the $\delta_0$-shorts of the pairs of pants immediately before a run of length $T$ eliminates the extra segment $\nu_0$.
    Then the same estimation scheme as in the multiple-pants case works.
\end{proof}

\section{Limit Cones}\label{s-limit-cones}

In this section, we deduce polyhedrality of limit cones of examples from the length combinatorialization of Prop. \ref{prop-combo}.
Let $\mathscr{P}$ be a convenient pants decomposition.

\begin{definition}
	Let $L = [(a_{ij})_{j=1, ..., n}]_{i =1, ..., 3g-3}$ be a sequence of $3g-3$ points in $(\bbR^{+})^n$ with all $a_{ij} < 1$ and let $C > 1$.
	Let $\mathcal{C}_L$ be the cone over the convex hull of $L$.
	Let 
	\begin{align*} 
		\mathcal{H}_L(C) &=  \left\{(b_1, ..., b_n) \in (\bbR^{+})^n  : \frac{b_i}{b_j} \leq C \max_{\substack{l=1, ..., 3g-3, \\ 1 \leq m , k \leq n}}  \frac{\log a_{lm}}{\log a_{lk}} , 1 \leq i, j \leq n \right\}, \\
			\mathcal{H}(C) &= \left\{ (b_1, ..., b_n) \in (\bbR^{+})^n :  \frac{b_i}{b_j} \leq C , 1 \leq i,  j \leq n \right\}
	\end{align*}

	For such an $L$, let $\rho_j \in \mathcal{T}(S)$ be given by $(a_{1j}, ..., a_{(3g-3)j}, 0, ..., 0)$ in Fenchel-Neilsen coordinates with respect to $\mathscr{P}$ for all $j$ and let $\rho_L = (\rho_1, ..., \rho_n) \in \mathcal{T}(S)^n$ be the direct sum representation.
	Call $\rho_L$ {{representation}} {\rm{associated to}} $L$. 
\end{definition}

Observe that $\mathcal{C}_L$ is invariant under simultaneously scaling the entries of $L$.
For $a \in \mathbb{R}^+$, write this rescaling by $aL$.

\begin{lemma}\label{lemma-cones-shrink}
	Let $L$ be so that $\mathcal{C}_L$ contains $(1, ..., 1)$.
	Then there exists a $C > 1$ and $\varepsilon > 0$ so that $\mathcal{H}_{aL}(C) - \{0\}$ is contained in the interior of $\mathcal{C}_L = \mathcal{C}_{aL}$ for all $0< a < \varepsilon$.
\end{lemma}

\begin{proof}
	For any cone $\mathcal{C}$ containing $(1, ..., 1)$ in its interior, the cones $\mathcal{H}(C)$ are contained in the interior of $\mathcal{C}$ for all $C$ sufficiently close to $1$.
	Take such a $C$ for $\mathcal{C}_L$.
		For $a < 1$,
	$$\mathcal{H}_{aL}(C) = \left\{ (b_1, ..., b_n) \in (\bbR^{+})^n : \frac{b_i}{b_j} \leq C \max_{\substack{l=1, ..., 3g-3, \\ 1 \leq m , k \leq n}} \frac{\log a_{lm} + \log a }{\log a_{lk} + \log a}, 1 \leq i,j \leq n  \right\},$$
	with the ratios of logarithm terms all tending to $1$ as $a \to 0$.
	So $\mathcal{H}_{aL}(C)$ converges (after projectivization) in the Hausdorff topology to $\mathcal{H}(C)$, and the claim follows.
\end{proof}

\begin{definition}
	For $\rho \in \mathcal{T}(S)^n$ and $\mathscr{P}$ a pants decomposition of $S$, denote the vector of length coordinates of $\rho$ as $L_\rho = (\ell_{ij})_{j=1, ..., n}^{i=1, ..., 3g-3}$.
	
	For $\rho \in \mathcal{T}(S)^n$, $T \geq 2$, and $\gamma$ a homotopy class of curves, write $\mathscr{R}_{\rho}(\gamma) =(\mathscr{R}(\gamma, \rho_i))_{i=1, ..., n}$ and $\mathscr{L}_{\rho}(\gamma) = (\mathscr{L}(\gamma, \rho_i))_{i=1, ..., n}$.
\end{definition}

Recall that the limit cone of a representation $\rho$ is denoted by $\mathcal{L}(\rho(\Gamma))$, and Prop. \ref{prop-Z-closures} and Benoist's limit cone theorem \cite[\S 1.2]{benoist1997proprietes} show $\mathcal{L}(\rho(\Gamma))$ is convex for all $\rho \in \mathcal{T}(S)^n$ with pairwise distinct projections.

\begin{theorem}\label{thm-cones}
	Suppose $\mathcal{C}_{L}$ contains $(1, ..., 1)$ in its interior.
	Then there is an $\varepsilon > 0$ so that for all $a < \varepsilon$, there is a neighborhood $U_{aL} \subset \mathcal{T}(S)^n$ of the representation $\rho_{aL}$ corresponding to $aL$ so that for all $\rho \in U_{aL}$, we have $\mathcal{C}_{L_{\rho}} = \mathcal{L}(\rho(\Gamma))$.
\end{theorem}

Note that the twist parameters do not appear in the limit cone, so that limit cones in these neighborhoods are constant along twisting, simultaneous scaling of entries, or adjusting Jordan projections of pants curves that land in $\text{Int}(\mathcal{C}_L)$.
This implies Cor. \ref{cor-cones-nonrigid}.

\begin{proof}
	Note that for any $L$, the cone $\mathcal{C}_{L}$ is the convex hull of the Jordan projections for $\rho_{L}$ of the curves in the pants decomposition $\mathscr{P},$ so that $\mathcal{C}_{L} \subset \mathcal{L}(\rho_L(\gamma))$.

	Now let $L$ be so that $\mathcal{C}_L$ contains $(1, ..., 1)$ in its interior.
	Take $B > 0$ as large as we wish and take $C$ sufficiently close to $1$ so that $\mathcal{H}(C) \subset \text{Int}( \mathcal{C}_L)$, let $\varepsilon > 0$ be as supplied by Prop. \ref{prop-combo} for these $C$ and $B$, and take $\varepsilon_0 \leq \varepsilon$ so that for all $a < \varepsilon_0$, every coordinate of $aL$ is less than $\varepsilon/2$.
	
	Let $a < \varepsilon_0/2$.
	Let us begin by proving that the claim holds for all $\rho$ with $0$ twist coordinates in a neighborhood $U_{aL}$ of $\rho_{aL}$, taken so that the following two properties hold.
	First, all $\rho \in U_{aL}$ have all their length coordinates less than $\varepsilon_0$.
	Second, $\mathcal{H}_{L_\rho}(C) $ is contained in the interior of $\mathcal{C}_{L_\rho}$ for all $\rho \in U_{aL}$.
	This can be arranged with Lemma \ref{lemma-cones-shrink}.
	
	So let $\rho$ be as specified.
	For the remaining containment, let $\gamma$ be a closed curve and write the Jordan projection of $\rho(\gamma)$ as $(\mathscr{R}(\gamma, \rho_i)+ \mathscr{L}(\gamma, \rho_i))_{i=1,..., n}$.
	Then by our arrangement of constants, Prop. \ref{prop-combo} shows $(\mathscr{L}(\gamma, \rho_i))_{i=1, ..., n} \in \mathcal{H}_{L_\rho}(C)$, which is contained in $\mathcal{C}_{L_\rho}$.
	Furthermore, by definition $\mathscr{R}(\rho, \gamma)$ is in $\mathcal{C}_{L_\rho}$.
	Because cones are semigroups and $\mathcal{H}_{L_\rho}(C) \subset \mathcal{C}_{L_\rho}$, the Jordan projection of $\gamma$ is in $\mathcal{C}_L$.
	As $\gamma$ was arbitrary, $\mathcal{C}_{L_\rho} = \mathcal{L}(\rho(\Gamma))$.
	
	We now handle twisting.
	Fix a $\rho$ as above and let $v = (v_{ij})_{i=1, ..., 3g-3}^{j=1,...,n} \in \bbR^{n(3g-3)}$ have all coordinates no more than $1$ in size.
	Let $\rho_{tv}$ be the path in $\mathcal{T}(S)^n$ with constant length coordinates and twist coordinate $tv$.
	Denote the factor representations of $\rho_{tv}$ by $\rho_{tv,i}$.
	
	Let $\gamma$ be any isotopy class of closed curves in $S$.
	From Wolpert's Cosine Formula \cite{kerckhoff1983nielsen,wolpert1981elementary}, we have $\frac{d}{dt} \ell_{\rho_{tv,i}}(\gamma) \leq \varepsilon_0 i(\gamma,\mathscr{P})$ for all $t$ and $v$.
	As $\mathscr{L}(\gamma, X) \geq  B i (\gamma, \mathscr{P})$, there is a $T_0$ depending only on the distance of $\mathcal{H}_{aL}(C)$ from $\partial \mathcal{C}_{L_\rho}$ when intersected with the hyperplane $||x||_{\mathrm{L}^1} = 1$ in $\bbR^n$ so that $\mathscr{L}(\gamma, \rho_{tv}) \in \mathcal{C}_{L_\rho}$ for all $t \leq T_0$.
	So for all $t \leq T_0$, we have $\mathcal{C}_{L_\rho} = \mathcal{L}(\rho_{tv}(\Gamma))$.
	
	As the separation of $\mathcal{H}_{L_\rho}(C)$ from $\partial \mathcal{C}_{L_\rho}$ is locally uniform in $\rho$, the claim follows.
\end{proof}

The results on limit cones claimed in the introduction follow immediately as special cases.

\subsection{Finite type surfaces}\label{ss-finite-type-limit-cones}
We briefly describe analogues in the finite-type case.
As in Part I, the statements to Thms.~\ref{theo:main-cones} and \ref{theo-designer-cone} are modified for surfaces $S_{g,p,b}$ with $p$ punctures, $b$ boundary components, and genus $g$ by replacing appearances of $3g-3$ with the maximum number $\kappa(S_{g,p,b}) = 3g-3 + 2b +p$ of disjoint non-peripheral simple curves in $S_{g,p,b}$.
Cor.~\ref{cor-cones-nonrigid} has the modification that the dimension count becomes $2n\kappa(S_{g,p,b}) - (nk+b) +1$.

For surfaces $S_{g,0,b}$ with boundary, these results are deduced by isometrically embedding $S_{g,0,b}$ into a closed surface and applying Prop.~\ref{prop-combo} there, as in the proof of Thm.~\ref{thm-cones}.

For surfaces with punctures, it is important to note that the proof of Prop.~\ref{prop-combo} does not directly work due to cusps interfering with the combinatorial structure we use.
However, there is a work-around to this for the desired statements: open up punctures to boundary components, apply Prop.~\ref{prop-combo}, and take a limit as the boundaries are pinched.

We now sketch the analogue of Thm.~\ref{theo:main-cones} in the same notation as its proof; the others are similar.
For untwisted $\rho_t \in \mathcal{T}(S_{g,p,b})^n$ under consideration that are deformations of a given $\rho$, realize $\rho_t$ as a limit of $\rho_{t,k} \in \mathcal{T}(S_{g,0,p+b})$, with all Fenchel-Nielsen coordinates not corresponding to the opened punctures the same as those of $\rho_t$.
Take a $C$ so that $\mathcal{H}_{L_\rho}(C)$ is contained in the interior of $\mathcal{C}_{L_\rho}$.
From continuity of Jordan projections, $\lambda(\rho_t(\gamma)) = \lim_{k\to\infty} \lambda(\rho_{t,k}(\gamma))$ for any $\gamma \in \pi_1(S_{g,p,b})$.

Now analyze $\lambda(\rho_t(\gamma))$ for a fixed $\gamma$.
Apply Prop.~\ref{prop-combo} to $\lambda(\rho_{t,k}(\gamma))$, note that contributions from rotation around the opened cusps become small compared to the size of $\lambda(\rho_{t,k}(\gamma))$ as $k$ grows large, and take the limit.
That Prop.~\ref{prop-combo} is uniform across all $X,Y \in \mathscr{U}_\varepsilon(\mathscr{P})$ and curves $\gamma$ implies that there is a $\gamma$-independent $a > 0$ so that $\lambda(\rho_{t}(\gamma))$ is the sum of an element of $\mathcal{C}_{L_\rho}$ and an element of $\mathcal{H}_{L_\rho}(C)$ for all $t < a$ and $\gamma \in \pi_1(S_{g,p,b})$.
The point is that the values of $k$ needed for a good approximation of this statement depend on $\gamma$, but the limiting statement of interest on $\lambda(\rho_{t}(\gamma))$ is $\gamma$-independent. 
This now implies $\mathcal{L}(\rho_t(\pi_1(S_{g,p,b}))) = \mathcal{C}_{L_\rho}$ for all $t < a$, as desired.
Twisting is handled as in the closed case.

\bibliographystyle{plain}
\bibliography{refs}

@article {bessonCourtoisGallot1995entropies,
	AUTHOR = {Besson, G. and Courtois, G. and Gallot, S.},
	TITLE = {Entropies et rigidit\'es des espaces localement sym\'etriques
	de courbure strictement n\'egative},
	JOURNAL = {Geom. Funct. Anal.},
	FJOURNAL = {Geometric and Functional Analysis},
	VOLUME = {5},
	YEAR = {1995},
	NUMBER = {5},
	PAGES = {731--799}
}

@article {hutchingsMorganRitpreRos2002proof,
	AUTHOR = {Hutchings, Michael and Morgan, Frank and Ritor\'e, Manuel and
	Ros, Antonio},
	TITLE = {Proof of the double bubble conjecture},
	JOURNAL = {Ann. of Math. (2)},
	FJOURNAL = {Annals of Mathematics. Second Series},
	VOLUME = {155},
	YEAR = {2002},
	NUMBER = {2},
	PAGES = {459--489}
}

@article {viazovska2017sphere,
	AUTHOR = {Viazovska, Maryna S.},
	TITLE = {The sphere packing problem in dimension 8},
	JOURNAL = {Ann. of Math. (2)},
	FJOURNAL = {Annals of Mathematics. Second Series},
	VOLUME = {185},
	YEAR = {2017},
	NUMBER = {3},
	PAGES = {991--1015}
}

@article {LangSchroeder1997kirszbraun,
    AUTHOR = {Lang, Urs and Schroeder, Viktor},
     TITLE = {Kirszbraun's theorem and metric spaces of bounded curvature},
   JOURNAL = {Geom. Funct. Anal.},
  FJOURNAL = {Geometric and Functional Analysis},
    VOLUME = {7},
      YEAR = {1997},
    NUMBER = {3},
     PAGES = {535--560}
}

@article {kerckhoff1980asymptotic,
    AUTHOR = {Kerckhoff, Steven P.},
     TITLE = {The asymptotic geometry of {T}eichm\"uller space},
   JOURNAL = {Topology},
  FJOURNAL = {Topology. An International Journal of Mathematics},
    VOLUME = {19},
      YEAR = {1980},
    NUMBER = {1},
     PAGES = {23--41},
      ISSN = {0040-9383}
}

@inproceedings {alexanderKapovitchPetrunin2011Alexandrov,
    AUTHOR = {Alexander, Stephanie and Kapovitch, Vitali and Petrunin,
              Anton},
     TITLE = {Alexandrov meets {K}irszbraun},
 BOOKTITLE = {Proceedings of the {G}\"{o}kova {G}eometry-{T}opology
              {C}onference 2010},
     PAGES = {88--109},
 PUBLISHER = {Int. Press, Somerville, MA},
      YEAR = {2011}
}

@article{Kirszbraun1934Uber,
  title={{\"U}ber die zusammenziehende und {L}ipschitzsche Transformationen},
  author={Kirszbraun, Moj{\.z}esz David},
  journal={Fundamenta Mathematicae},
  year={1934},
  volume={22},
  pages={77-108}
}

@article {gueritaudKassel2017maximally,
	AUTHOR = {Gu{\'e}ritaud, Fran{\c c}ois and Kassel, Fanny},
	TITLE = {Maximally stretched laminations on geometrically finite
	hyperbolic manifolds},
	JOURNAL = {Geom. Topol.},
	FJOURNAL = {Geometry \& Topology},
	VOLUME = {21},
	YEAR = {2017},
	NUMBER = {2},
	PAGES = {693--840},
	ISSN = {1465-3060,1364-0380}
}

@misc{dancigerGueritaudKassel2025limit,
	author = {Danciger, Jeffrey and Gu{\' e}ritaud, Fran{\c c}ois and Kassel, Fanny},
	title = {Limit cones of multi-{F}uchsian representations},
	year = {2025},
	note = {Preprint}
}

@article {Valentine1945lipschitz,
    AUTHOR = {Valentine, Frederick Albert},
     TITLE = {A {L}ipschitz condition preserving extension for a vector
              function},
   JOURNAL = {Amer. J. Math.},
  FJOURNAL = {American Journal of Mathematics},
    VOLUME = {67},
      YEAR = {1945},
     PAGES = {83--93}
}

@misc{tsouvalas2023holder,
    author = {Tsouvalas, Konstantinos},
    title  = {The {H\"o}lder exponent of {A}nosov limit maps},
    year   = {2023},
    note   = {\href{https://arxiv.org/abs/2306.15823}{arxiv:2306.15823}}
}

@article {sambarino2015orbital,
    AUTHOR = {Sambarino, Andr\'{e}s},
     TITLE = {The orbital counting problem for hyperconvex representations},
   JOURNAL = {Ann. Inst. Fourier (Grenoble)},
  FJOURNAL = {Universit\'{e} de Grenoble. Annales de l'Institut Fourier},
    VOLUME = {65},
      YEAR = {2015},
    NUMBER = {4},
     PAGES = {1755--1797}
}

@incollection {benoist2004convexesI,
    AUTHOR = {Benoist, Yves},
     TITLE = {Convexes divisibles. {I}},
 BOOKTITLE = {Algebraic groups and arithmetic},
     PAGES = {339--374},
 PUBLISHER = {Tata Inst. Fund. Res., Mumbai},
      YEAR = {2004}
}

@article{mirzakhani2008growth,
	AUTHOR = {Mirzakhani, Maryam},
	TITLE = {Growth of the number of simple closed geodesics on hyperbolic
	surfaces},
	JOURNAL = {Ann. of Math. (2)},
	FJOURNAL = {Annals of Mathematics. Second Series},
	VOLUME = {168},
	YEAR = {2008},
	NUMBER = {1},
	PAGES = {97--125}
}

@article {guichard2005regularity,
    AUTHOR = {Guichard, Olivier},
     TITLE = {Sur la r\'{e}gularit\'{e} {H}\"{o}lder des convexes divisibles},
   JOURNAL = {Ergodic Theory Dynam. Systems},
  FJOURNAL = {Ergodic Theory and Dynamical Systems},
    VOLUME = {25},
      YEAR = {2005},
    NUMBER = {6},
     PAGES = {1857--1880}
}

@misc{calderonTao2025deflating,
	title={Deflating hyperbolic surfaces and the shapes of optimal {L}ipschitz maps}, 
	author={Aaron Calderon and Jing Tao},
	year={2025},
	note={\href{https://arxiv.org/abs/2510.19930}{arxiv:2510.19930}}
}

@article {zhangZimmer2024regularity,
	AUTHOR = {Zhang, Tengren and Zimmer, Andrew},
	TITLE = {Regularity of limit sets of {A}nosov representations},
	JOURNAL = {J. Topol.},
	FJOURNAL = {Journal of Topology},
	VOLUME = {17},
	YEAR = {2024},
	NUMBER = {3},
	PAGES = {Paper No. e12355, 72},
	ISSN = {1753-8416,1753-8424}
}

@misc{kasselTholozan2024sharpness,
	title={Sharpness of proper and cocompact actions on reductive homogeneous spaces}, 
	author={Fanny Kassel and Nicolas Tholozan},
	year={2024},
	note={\href{https://arxiv.org/abs/2410.08179}{arxiv:2410.08179}}
}

@article {lenzhenRafiTao2012bounded,
	AUTHOR = {Lenzhen, Anna and Rafi, Kasra and Tao, Jing},
	TITLE = {Bounded combinatorics and the {L}ipschitz metric on
	{T}eichm{\"u}ller space},
	JOURNAL = {Geom. Dedicata},
	FJOURNAL = {Geometriae Dedicata},
	VOLUME = {159},
	YEAR = {2012},
	PAGES = {353--371},
	ISSN = {0046-5755,1572-9168}
}

@article {lenzhenRafiTao2015shadow,
	AUTHOR = {Lenzhen, Anna and Rafi, Kasra and Tao, Jing},
	TITLE = {The shadow of a {T}hurston geodesic to the curve graph},
	JOURNAL = {J. Topol.},
	FJOURNAL = {Journal of Topology},
	VOLUME = {8},
	YEAR = {2015},
	NUMBER = {4},
	PAGES = {1085--1118},
	ISSN = {1753-8416,1753-8424}
}

@article {dumasLenzhenRafiTao2020coarse,
	AUTHOR = {Dumas, David and Lenzhen, Anna and Rafi, Kasra and Tao, Jing},
	TITLE = {Coarse and fine geometry of the {T}hurston metric},
	JOURNAL = {Forum Math. Sigma},
	FJOURNAL = {Forum of Mathematics. Sigma},
	VOLUME = {8},
	YEAR = {2020},
	PAGES = {Paper No. e28, 58},
	ISSN = {2050-5094}
}

@article {pan2023local,
	AUTHOR = {Pan, Huiping},
	TITLE = {Local rigidity of the {T}eichm\"uller space with the
	{T}hurston metric},
	JOURNAL = {Sci. China Math.},
	FJOURNAL = {Science China. Mathematics},
	VOLUME = {66},
	YEAR = {2023},
	NUMBER = {8},
	PAGES = {1751--1766},
	ISSN = {1674-7283,1869-1862}
}

@misc{barnatanOhshikaPapadopoulos2025convex,
	title={Convex structures of the unit tangent spheres in {T}eichm{\"u}ller space}, 
	author={Assaf Bar-Natan and Ken'Ichi Ohshika and Athanase Papadopoulos},
	year={2025},
	note={\href{https://arxiv.org/abs/2503.20404}{arxiv:2503.20404}},
}

@misc{panWolf2024envelopes,
	title={Envelopes of the {T}hurston metric on {T}eichm\"uller space}, 
	author={Huiping Pan and Michael Wolf},
	year={2024},
	note={\href{https://arxiv.org/abs/2401.06607}{arxiv:2401.06607}}, 
}

@book {beardon1995geometry,
	AUTHOR = {Beardon, Alan F.},
	TITLE = {The geometry of discrete groups},
	SERIES = {Graduate Texts in Mathematics},
	VOLUME = {91},
	NOTE = {Corrected reprint of the 1983 original},
	PUBLISHER = {Springer-Verlag, New York},
	YEAR = {1995},
	PAGES = {xii+337},
	ISBN = {0-387-90788-2}
}

@article {chasMcMullenPhillips2019almost,
	AUTHOR = {Chas, Moira and McMullen, Curtis T. and Phillips, Anthony},
	TITLE = {Almost simple geodesics on the triply-punctured sphere},
	JOURNAL = {Math. Z.},
	FJOURNAL = {Mathematische Zeitschrift},
	VOLUME = {291},
	YEAR = {2019},
	NUMBER = {3-4},
	PAGES = {1175--1196},
	ISSN = {0025-5874,1432-1823}
}

@article {contreras2016ground,
	AUTHOR = {Contreras, Gonzalo},
	TITLE = {Ground states are generically a periodic orbit},
	JOURNAL = {Invent. Math.},
	FJOURNAL = {Inventiones Mathematicae},
	VOLUME = {205},
	YEAR = {2016},
	NUMBER = {2},
	PAGES = {383--412},
	ISSN = {0020-9910,1432-1297}
}

@inbook{bochi2018ergodic,
	author = {Jairo Bochi},
	title = {Ergodic optimization of {B}irkhoff averages and {L}yapunov exponents},
	booktitle = {Proceedings of the International Congress of Mathematicians (ICM 2018)},
	chapter = {},
    year = {2018},
	pages = {1825-1846},
	doi = {10.1142/9789813272880_0119}
}

@article{huntOtt1996optimal,
	title = {Optimal periodic orbits of chaotic systems occur at low period},
	author = {Hunt, Brian R. and Ott, Edward},
	journal = {Phys. Rev. E},
	volume = {54},
	issue = {1},
	pages = {328--337},
	year = {1996}
}

@article{deyOh2025deformations,
	author = {Dey, Subhadip and Oh, Hee},
	title = {Deformations of {A}nosov subgroups: Limit cones and growth indicators},
	journal = {J. Lond. Math. Soc.},
	volume = {112},
	number = {3},
	pages = {e70280},
    year = {2025}
}

@article {huangOhshikaPapadopoulos2025infinitesimal,
	AUTHOR = {Huang, Yi and Ohshika, Kenichi and Papadopoulos, Athanase},
	TITLE = {The infinitesimal and global {T}hurston geometry of
	{T}eichm\"uller space},
	JOURNAL = {J. Differential Geom.},
	FJOURNAL = {Journal of Differential Geometry},
	VOLUME = {131},
	YEAR = {2025},
	NUMBER = {2},
	PAGES = {--},
	ISSN = {0022-040X,1945-743X}
}

@article {minsky1996extremal,
	AUTHOR = {Minsky, Yair N.},
	TITLE = {Extremal length estimates and product regions in
	{T}eichm\"uller space},
	JOURNAL = {Duke Math. J.},
	FJOURNAL = {Duke Mathematical Journal},
	VOLUME = {83},
	YEAR = {1996},
	NUMBER = {2},
	PAGES = {249--286},
	ISSN = {0012-7094,1547-7398}
}

@misc{davaloRiestenberg2024finitesided,
	title={Finite-sided {D}irichlet domains and {A}nosov subgroups}, 
	author={Colin Davalo and J. Maxwell Riestenberg},
	year={2024},
	note={\href{https://arxiv.org/abs/2308.08607}{arxiv:2402.06408}},
}

@misc{panWolf2025ray,
	title = {Ray structures on {T}eichm{\"u}ller space},
	author = {Pan, Huiping and Wolf, Michael},
	note = {To appear in \textit{Acta Math.}}
}

@article {daskalopoulosUhlenbeck2024transverse,
	AUTHOR = {Daskalopoulos, Georgios and Uhlenbeck, Karen},
	TITLE = {Transverse measures and best {L}ipschitz and least gradient
	maps},
	JOURNAL = {J. Differential Geom.},
	FJOURNAL = {Journal of Differential Geometry},
	VOLUME = {127},
	YEAR = {2024},
	NUMBER = {3},
	PAGES = {969--1018},
	ISSN = {0022-040X,1945-743X}
}

@book{fathiLaudenbachPoenaru2012thurston,
	author = {Albert Fathi and François Laudenbach and Valentin Poénaru},
	publisher = {Princeton University Press},
	title = {Thurston's Work on Surfaces},
	volume = {48},
	year = {2012},
	note = {Translated by Djun M. Kim and Dan Margalit}
}

@article{alessandriniDisarlo2022generalizing,
	AUTHOR = {Alessandrini, Daniele and Disarlo, Valentina},
	TITLE = {Generalizing stretch lines for surfaces with boundary},
	JOURNAL = {Int. Math. Res. Not. IMRN},
	FJOURNAL = {International Mathematics Research Notices. IMRN},
	YEAR = {2022},
	NUMBER = {23},
	PAGES = {18919--18991},
	ISSN = {1073-7928,1687-0247}
}

@misc{cantrellCertReyes2024jointtranslation,
	title={The joint translation spectrum and {M}anhattan manifolds}, 
	author={Stephen Cantrell and Eduardo Reyes and Cagri Sert},
	year={2024},
	note={\href{https://arxiv.org/abs/2411.06375}{arxiv:2411.06375}} 
}

@article {theret2010divergence,
	AUTHOR = {Th\'eret, Guillaume},
	TITLE = {Divergence et parall\'elisme des rayons d'\'etirement
	cylindriques},
	JOURNAL = {Algebr. Geom. Topol.},
	FJOURNAL = {Algebraic \& Geometric Topology},
	VOLUME = {10},
	YEAR = {2010},
	NUMBER = {4},
	PAGES = {2451--2468},
	ISSN = {1472-2747,1472-2739}
}

@incollection {su2016problems,
	AUTHOR = {Su, Weixu},
	TITLE = {Problems on the {T}hurston metric},
	BOOKTITLE = {Handbook of {T}eichm\"uller theory. {V}ol. {V}},
	SERIES = {IRMA Lect. Math. Theor. Phys.},
	VOLUME = {26},
	PAGES = {55--72},
	PUBLISHER = {Eur. Math. Soc., Z\"urich},
	YEAR = {2016},
	ISBN = {978-3-03719-160-6}
}

@article {dalBoKim2000criterion,
	AUTHOR = {Dal'Bo, Fran{\c c}oise and Kim, Inkang},
	TITLE = {A criterion of conjugacy for {Z}ariski dense subgroups},
	JOURNAL = {C. R. Acad. Sci. Paris S\'er. I Math.},
	FJOURNAL = {Comptes Rendus de l'Acad{\'e}mie des Sciences. S{\'e}rie I.
	Math{\'e}matique},
	VOLUME = {330},
	YEAR = {2000},
	NUMBER = {8},
	PAGES = {647--650},
	ISSN = {0764-4442}
}

@article {borel1960density,
	AUTHOR = {Borel, Armand},
	TITLE = {Density properties for certain subgroups of semi-simple groups
	without compact components},
	JOURNAL = {Ann. of Math. (2)},
	FJOURNAL = {Annals of Mathematics. Second Series},
	VOLUME = {72},
	YEAR = {1960},
	PAGES = {179--188},
	ISSN = {0003-486X}
}

@article {papadopoulosTheret2012someLipschitz,
	AUTHOR = {Papadopoulos, Athanase and Th\'eret, Guillaume},
	TITLE = {Some {L}ipschitz maps between hyperbolic surfaces with
	applications to {T}eichm\"uller theory},
	JOURNAL = {Geom. Dedicata},
	FJOURNAL = {Geometriae Dedicata},
	VOLUME = {161},
	YEAR = {2012},
	PAGES = {63--83},
	ISSN = {0046-5755,1572-9168}
}

@article {papadopoulosYamada2017deforming,
	AUTHOR = {Papadopoulos, Athanase and Yamada, Sumio},
	TITLE = {Deforming hyperbolic hexagons with applications to the arc and
	the {T}hurston metrics on {T}eichm\"uller spaces},
	JOURNAL = {Monatsh. Math.},
	FJOURNAL = {Monatshefte f\"ur Mathematik},
	VOLUME = {182},
	YEAR = {2017},
	NUMBER = {4},
	PAGES = {913--939},
	ISSN = {0026-9255,1436-5081}
}

@misc{thurston1998minimal,
    title = {Minimal stretch maps between hyperbolic surfaces},
    author = {Thurston, William P.},
    year = {1998},
    note = {\href{https://arxiv.org/abs/math/9801039}{arxiv:9801039} (digitized version of unpublished 1986 preprint)}
}

@article {wolpert1981elementary,
	AUTHOR = {Wolpert, Scott},
	TITLE = {An elementary formula for the {F}enchel-{N}ielsen twist},
	JOURNAL = {Comment. Math. Helv.},
	FJOURNAL = {Commentarii Mathematici Helvetici},
	VOLUME = {56},
	YEAR = {1981},
	NUMBER = {1},
	PAGES = {132--135},
	ISSN = {0010-2571,1420-8946}
}

@article {kerckhoff1983nielsen,
	AUTHOR = {Kerckhoff, Steven P.},
	TITLE = {The {N}ielsen realization problem},
	JOURNAL = {Ann. of Math. (2)},
	FJOURNAL = {Annals of Mathematics. Second Series},
	VOLUME = {117},
	YEAR = {1983},
	NUMBER = {2},
	PAGES = {235--265},
	ISSN = {0003-486X,1939-8980}
}

@article{leeZhang2017collar,
    title = {Collar lemma for {H}itchin representations},
    author = {Lee, Gye-Seon and Zhang, Tengren},
    journal = {Geom. Topol.},
    volume = {21},
    year = {2017},
    pages = {2243-2280}
}

@article{benoist1997proprietes,
    AUTHOR = {Benoist, Yves},
     TITLE = {Propri\'{e}t\'{e}s asymptotiques des groupes lin\'{e}aires},
   JOURNAL = {Geom. Funct. Anal.},
  FJOURNAL = {Geometric and Functional Analysis},
    VOLUME = {7},
      YEAR = {1997},
    NUMBER = {1},
     PAGES = {1--47}
}

\end{document}